\documentclass{amsart}

\usepackage{amssymb}
\usepackage{eepic}
\usepackage{color}
\usepackage{mathtools,xparse}
\usepackage[linktocpage=true]{hyperref}
\usepackage[utf8x]{inputenc}

\allowdisplaybreaks

\numberwithin{equation}{section}

\newtheorem{theorem}{Theorem}[section]
\newtheorem{lemma}[theorem]{Lemma}
\newtheorem{proposition}[theorem]{Proposition}
\newtheorem{corollary}[theorem]{Corollary}
\newtheorem{remark}[theorem]{Remark}
\newtheorem{problem}[theorem]{Problem}

\newtheorem{TheoA}{Theorem A}
\newtheorem{TheoB}{Theorem B}
\newtheorem{TheoLocal}{Local form of Theorem A}

\newcommand{\N}{\mathbf{N}}
\newcommand{\Z}{\mathbf{Z}}
\newcommand{\R}{\mathbf{R}}
\newcommand{\C}{\mathbf{C}}

\newcommand{\SL}{\mathrm{SL}_n(\R)}
\newcommand{\GL}{\mathrm{GL}_n(\R)}
\newcommand{\SLk}{\mathrm{SL}_d(\R)}
\newcommand{\sphere}{\ensuremath{\mathbf{S}}}
\newcommand{\summ}{\sum\nolimits}

\newcommand{\weight}[1]{\left\bracevert \hskip-3pt #1  \hskip-3pt \right\bracevert}

\def\G{\mathrm{G}}
\def\1{\mathbf{1}}
\def\H{\mathcal{H}}
\def\Q{\mathcal{Q}}
\def\M{\mathcal{M}}

\def\RR{\mathcal{R}_\Sigma}

\def\V{\mathrm{\mathcal{L}(G)}}

\newcommand{\dem}{\noindent {\bf Proof. }}

\newcommand{\demAGral}{\noindent {\bf Proof of Theorem A. }}

\newcommand{\demB}{\noindent {\bf Proof of Theorem B. }}
\newcommand{\demC}{\noindent {\bf Proof of Theorem \ref{thm:multipliers_SO(n,1)}. }}
\newcommand{\fin}{\hspace*{\fill} $\square$ \vskip0.2cm}

\def\esssup{\mathop{\mathrm{ess \, sup}}}

\begin{document}

\null

\vskip-50pt

\null

\title[Fourier multipliers in $\SL$]{Fourier multipliers in $\SL$}
\author[J. Parcet, \'E. Ricard, M. de la Salle]{Javier Parcet, \'Eric Ricard \\ and  Mikael de la Salle}

\maketitle

\null

\vskip-50pt

\null

\begin{abstract}
We establish precise regularity conditions for $L_p$-boundedness of Fourier multipliers in the group algebra of $\SL$. Our main result is inspired by the H\"ormander-Mikhlin criterion from classical harmonic analysis, although it is substantially and necessarily different. Locally, we get sharp growth rates of Lie derivatives around the singularity and nearly optimal regularity. The asymptotics also match Mikhlin formula for an exponentially growing weight with respect to the word length. Additional decay comes imposed by this growth and Mikhlin condition for high order terms. Lafforgue/de la Salle's rigidity theorem fits here. The proof includes a new relation between Fourier and Schur $L_p$-multipliers for nonamenable groups. 
By transference, matters are reduced to a rather nontrivial $RC_p$-inequality for $\SL$-twisted forms of Riesz transforms associated to fractional Laplacians. 

\vskip3pt

Our second result gives a new and much stronger rigidity theorem for radial multipliers in $\SL$. More precisely, additional regularity and Mikhlin type conditions are proved to be necessary up to an order $\sim |\frac12 - \frac1p| (n-1)$ for large enough $n$ in terms of $p$. Locally, necessary and sufficient growth rates match up to that order. Asymptotically, extra decay for the symbol and its derivatives imposes more accurate and additional rigidity in a wider range of $L_p$-spaces. This rigidity increases with the rank, so we can construct radial generating functions satisfying our H\"ormander-Mikhlin sufficient conditions in a given rank $n$ and failing the rigidity conditions for ranks $m >> n$. We also prove automatic regularity and rigidity estimates for first and higher order derivatives of $\mathrm{K}$-biinvariant multipliers in the rank 1 groups $\mathrm{SO}(n,1)$. 
\end{abstract}

\addtolength{\parskip}{+1ex}

\vskip10pt

\section*{\bf Introduction}

We study the relation between regularity and $L_p$-boundedness for multipliers in the group algebra of $\SL$. In Euclidean harmonic analysis, this central topic orbits around the H\"ormander-Mikhlin fundamental condition \cite{Ho,Mi}. It defines a large class of Fourier multipliers |including Riesz transforms and Littlewood-Paley partitions of unity| which are crucial in Fourier summability or pseudodifferential operator theory. Given a measurable function $m: \R^n \to \C$, its Fourier multiplier is the linear map determined by $$\widehat{T_mf}(\xi) = m(\xi) \widehat{f}(\xi).$$ Then, $T_m$ is $L_p$-bounded on $\R^n$ for $1 < p < \infty$ whenever 
\begin{equation} \tag{HM} \label{Eq-HM}
\big| \partial_\xi^\gamma m(\xi) \big| \, \lesssim \, |\xi|^{-|\gamma|} \quad \mbox{for all} \quad 0 \le |\gamma| \le \Big[ \frac{n}{2} \Big] + 1.
\end{equation}
This condition imposes $m$ to be bounded, smooth over $\R^n \setminus \{0\}$ and to satisfy certain local and asymptotic behavior. Locally, $m$ admits a singularity at $0$ with a mild control of derivatives around it up to order $[\frac{n}{2}] +1$. This singularity links to deep concepts in harmonic analysis and justifies the key role of the H\"ormander-Mikhlin theorem in Fourier multiplier $L_p$-theory. Asymptotically, the same derivatives decay faster and faster to $0$, at a polynomial rate given by the differentiation order. It is optimal in the sense that we may not consider less classical derivatives |a Sobolev type formulation (recalled below in this paper) admits differentiability orders up to $\frac{n}{2} + \varepsilon$| or larger upper bounds for them. 

The H\"ormander-Mikhlin theorem has been investigated during the last decades for nilpotent groups by Christ, Cowling, M\"uller, Ricci or Stein among others. In the context of semisimple Lie groups, these questions have only been considered under additional symmetry assumptions. Let $\G$ be a real semisimple (noncompact and connected) Lie group with finite center and $\mathrm{K}$ be maximal compact in $\G$. Consider the Riemannian symmetric space $\G/\mathrm{K}$ equipped with its $\G$-invariant measure under left multiplication. If $\mathfrak{a}$ is a Cartan subalgebra for $(\G,\mathrm{K})$, the $\G$-invariant maps on $L_2(\G/\mathrm{K})$ can be identified via the spherical transform with Weyl-group-invariant elements in $L_\infty(\mathfrak{a}^*)$. H\"ormander-Mikhlin criteria for this class of multipliers were first considered by Clerc and Stein \cite{CSt}, which established a necessary analyticity condition and a weak form of Mikhlin sufficient condition. Stanton and Tomas \cite{ST} obtained nearly optimal results in rank one based on precise local/asymptotic expansion formulae for spherical functions on $\G/\mathrm{K}$. Anker \cite{A} finally discovered satisfactory Mikhlin conditions in high ranks. We refer to \cite{RW} for an interesting generalization in $\mathrm{SL}_2(\R)$ and to \cite{AL,I1,I2,LM} for related results.

In this paper we work with the full semisimple Lie group $\SL$ and place it in the frequency side. Its dual is no longer a group and it is described as a group von Neumann algebra, a key model of quantum (nonclassical) group. The interest of Fourier multipliers over group algebras was early recognized in the pioneering work of Haagerup \cite{H}, as well as in the research carried out thereafter in the context of approximation properties \cite{DCH,CDSW,CH,H2}. The corresponding theory of Fourier $L_p$-multipliers is basic in noncommutative harmonic analysis, with potential applications in geometric group theory and operator algebra. It has recently gained a considerable momentum \cite{GJP,JMP1,JMP2,JR,dLdlS,LdlS,MR,PRo}. 

If $(\mathcal M,\tau_{\mathcal M})$ is a semi-finite von Neumann algebra, $L_p(\mathcal M)$ denotes the associated noncommutative $L_p$-space (we skip the reference to the trace as it will always be natural). We refer to \cite{PX} for precise definitions. When dealing with $\mathcal B(\H)$, the bounded linear operators on a Hilbert space $\H$ with its canonical trace, the resulting noncommutative $L_p$-spaces are exactly the Schatten $p$-classes $S_p(\H)$. When $\H = \C^n$, we simply write $S_p^n$ when $p<\infty$ and $M_n$ for $p=\infty$. Given $(\mathcal{M},\tau_{\mathcal{M}})$, the canonical basis of $\C^n$ yields an identification $M_n(\mathcal{M})=M_n\otimes \mathcal{M}$. Using the tensor product trace at the $L_p$-level, we can also identify $L_p(M_n(\mathcal{M}))$ with $S_p^n\otimes L_p(\mathcal{M})$.
Given another semi-finite algebra and $(\mathcal{N},\tau_{\mathcal{N}})$ and a map $T: L_p(\mathcal{M})\to L_p(\mathcal{N})$, we say that it is completely bounded if  
$$ \|T\|_{\mathrm{cb}}=\sup_{n\geq 1}  \big\| \mathrm{Id}_{S_p^n}\otimes T :  L_p(M_n(\mathcal{M}))\to L_p(M_n(\mathcal{N})) \big\| <\infty.$$
Using completely bounded maps rather than bounded maps is often a price to pay 
to deal with noncommutative integration. We may use  $\leq_{\mathrm{cb}}$ to say that an inequality 
also holds for matrix ampliations with the same constants. For instance, if $T_i:L_p(\mathcal{M})\to L_p(\mathcal{N})$
are completely bounded maps 
$$\Big\| \Big(\sum_{i=1}^d |T_i(f_i)|^2\Big)^{\frac12}\Big\|_p \leq_{\mathrm{cb}} C \Big\| \Big(\sum_{i=1}^d |f_i|^2\Big)^{\frac12}\Big\|_p$$
means that for any $n$ and any $F_i\in L_p(M_n(\mathcal{M}))$, 
$$\Big\| \Big(\sum_{i=1}^d \big| (\mathrm{Id}_{S_p^n}\otimes T_i) (F_i) \big|^2\Big)^{\frac12}\Big\|_{L_p(M_n(\mathcal{N}))} \leq_{\mathrm{cb}} C \Big\| \Big(\sum_{i=1}^d |F_i|^2\Big)^{\frac12}\Big\|_{L_p(M_n(\mathcal{M}))}.$$

Given a locally compact unimodular group $\G$ with left regular representation $\lambda$, its group von Neumann algebra $\V$ is the weak-$*$ closure in $\mathcal{B}(L_2(\G))$ of $\mathrm{span}(\lambda(\G))$. If $\mu$ denotes the Haar measure of $\G$, we may approximate every element affiliated to $\V$ by operators of the form $$f \, = \, \int_\G \widehat{f}(g) \lambda(g) \, d\mu(g)$$ for smooth enough $\widehat{f}$ (say $\widehat{f}\in \mathcal{C}_c(\G)*\mathcal{C}_c(\G)$).
If $e$ is the unit in $\G$, the Haar trace $\tau$ is then the normal semifinite trace determined by $\tau(f) = \widehat{f}(e)$ for every such $f$. As usual we call $\widehat f$ the Fourier transform of $f$, it gives an isometry between $L_2(\V)$ and $L_2(\G)$. Given a bounded measurable symbol $m: \G \to \C$, its associated Fourier multiplier is the map defined on $L_2(\V)$ as follows 
$$T_m\Big(\int_\G \widehat{f}(g) \lambda(g) \, d\mu(g)\Big)=\int_\G m(g)\widehat{f}(g) \lambda(g) \, d\mu(g).$$
It intertwines pointwise multiplication with Fourier transform $\widehat{T_m f}(g) =m(g) \widehat{f}(g)$.   
We are interested to know when these maps extend from $L_p(\V) \cap L_2(\V)$ (or equivalently from smooth functions) to a bounded map on  $L_p(\V)$.

To such a symbol $m: \G \to \C$, one can also  associate Herz-Schur multipliers (or simply Schur multipliers). At the $L_2$-level,
any element $A \in S_2(L_2(\G))$ has a kernel $(A_{gh})_{g,h\in \G}$ which sits in $L_2(\G\times \G)$.
The Schur multiplier $S_m$ is the map $S_2(L_2(\G))\to S_2(L_2(\G))$ given by $$S_m \big( (A_{gh})_{g,h\in \G} \big) = \big( m(gh^{-1})A_{gh} \big)_{g,h\in \G}.$$ We aim to know if $S_m$ extends to $S_p(L_2(\G))$. Actually, we are mostly interested in restrictions of $S_m$ to $S_p(L_2(\Sigma))$ where $\Sigma$ is a measurable set in $\G$. To avoid heavy notation, we will 
still write $S_m$ for its restrictions but we will always make them precise in the norm by writing $ \| S_m \|_{\mathrm{cb}(S_p(L_2(\Sigma)))}$ or $\|S_m\|_{\mathcal B(S_p(L_2(\Sigma)))}$.

There are close links between Fourier and Herz-Schur multipliers, which we call transference \cite{CS,LdlS, NR}. The most classical instance of this relation is that the complete $L_p$-boundedness of a Fourier multiplier with symbol $m$  on $\G$ ensures that $S_m$ is completely $S_p(L_2(\G))$-bounded. Known results in the opposite 
direction require $\G$ to be amenable. To illustrate this, consider $\G=\Z$, 
any $f \in \mathcal L (\Z)=L_\infty(\mathbf T)$ is represented in $\mathcal B(\ell_2(\Z))$ by a matrix $(\widehat f (i-j))_{i,j\in \Z}$. Thus formally we have 
$$ S_m(f) = \big( T_m(f)(i-j) \big)_{i,j\in \Z}.$$

\noindent \textbf{I. Main results.} The rigidity theorems in \cite{dLdlS,LdlS} establish the failure of the completely bounded approximation property (CBAP) in the noncommutative $L_p$ space over the group von Neumann algebra of any lattice in $\SL$, for certain values of $p$. Roughly speaking, what is behind in harmonic analysis words is that Fourier $L_p$-summability fails in the group algebra of $\SL$, when $|1/p-1/2|$ is large enough in terms of the rank. This constitutes an $L_p$-refinement of Haagerup's theorem on the failure of weak amenability for higher rank simple Lie groups \cite{CDSW,H2}. As far as these group algebras are concerned, Fourier multiplier theory has been limited so far to rigidity theorems and the search of \emph{necessary} conditions to this end. In this paper, we provide  \emph{sufficient} conditions for $L_p$ boundedness of Fourier multipliers in the group algebra of $\SL$. 

The aforementioned rigidity gets in conflict with the classical H\"ormander-Mikhlin criterion. Indeed, condition \eqref{Eq-HM} certainly includes $\mathcal{C}_0$-functions with arbitrarily slow decay. On the contrary, any K-biinvariant symbol $m \in \mathcal{C}_0(\mathrm{SL}_3(\R))$ satisfies the asymptotic rigidity estimate 
\begin{equation} \tag{$\mathrm{AR}$} \label{Eq-Rig-p}
\left| m \left( \begin{array}{ccc} e^s & 0 & 0 \\ 0 & 1 & 0 \\ 0 & 0 & e^{-s} \end{array} \right) \right| \, \le \, C e^{-\delta |s|} \big\| T_m \big\|_{\mathrm{cb}(L_p(\mathcal{L}(\mathrm{SL}_3(\R))))}
\end{equation}
for all $p > 4$ and any $\delta < 1/2 - 2/p$. This is the key inequality in \cite{LdlS}. Considerably stronger rigidity estimates for radial Fourier multipliers are proved in Theorem B below. Asymptotic rigidity is not witnessed by H\"ormander-Mikhlin conditions and sufficient conditions for $L_p$-boundedness must incorporate it. To do so, let us write $\| \hskip3pt \|$ for the operator norm in $M_n(\R)$. Define the following weight function 
\[ \weight{g} := \max\Big\{ \|g-e\|, \| g^{-1}-e \| \Big\}.\]
Observe that, on a neighborhood of the identity, the weight $\weight{g}$ is comparable to the Euclidean distance from $g$ to $e$ in $\SL$. Asymptotically, it is comparable to 
\[L(g):=\max \Big\{ \|g\|, \|g^{-1}\| \Big\}.\]  Note in passing that $(g,h) \mapsto \log L(gh^{-1})$ is a standard left-invariant pseudometric, comparable to the word length from a compact symmetric generating set in $\SL$. 

We work with the natural differential operators. Consider the left-invariant vector fields generated by an orthonormal basis $\mathrm{X}_1, \mathrm{X}_2, \ldots, \mathrm{X}_{n^2-1}$ of the Lie algebra $\mathfrak{sl}_n(\R)$. The corresponding Lie derivatives  
$$\partial_{\mathrm{X}_j} m (g) \, = \, \frac{d}{ds}\Big|_{s=0} m \big( g \exp(s \mathrm{X}_j) \big)$$
do not commute for $j \neq k$. This justifies to define the set of multi-indices $\gamma$ as ordered tuples $\gamma = (j_1, j_2, \ldots, j_k)$ with $1 \le j_i \le n^2 -1$ and $|\gamma| = k \ge 0$, which correspond to the Lie differential operators $$d_g^\gamma m (g) \, = \, \partial_{\mathrm{X}_{j_1}} \partial_{\mathrm{X}_{j_2}} \cdots \, \partial_{\mathrm{X}_{j_{|\gamma|}}} m (g) \, = \, \Big( \prod_{1 \le k \le |\gamma|}^{\rightarrow} \partial_{\mathrm{X}_{j_k}} \Big) m(g).$$ 

\begin{TheoA}
Assume that $m \in \mathcal{C}^{[\frac{n^2}{2}]+1}(\SL \setminus \{e\})$ satisfies
\begin{equation} \label{Eq-ThmA} \tag{$\star$}
\weight{g}^{|\gamma|} \big| d_g^\gamma m(g) \big| \, \le \, C_{\mathrm{hm}} \quad \mbox{for all} \quad |\gamma| \le \Big[\frac{n^2}{2} \Big] + 1.
\end{equation} \vskip-2pt
\noindent Then, $T_m$ is completely $L_p$-bounded for all $1 < p < \infty$ by a constant $C_p C_{\mathrm{hm}}$.
\end{TheoA}
We point out that \eqref{Eq-ThmA} is a way of expressing that, if we equip $\SL$ with a left-invariant Riemannian metric, the norm of $k$-th derivative in the Riemanian sense of $f$ at $g$ is bounded above by $C_{\rm hm} \!\! \weight{g}^{-k}$.

Theorem A gives the first sufficient condition for $L_p$-boundedness in the group algebra of $\SL$. Compared to \eqref{Eq-HM}, we have replaced Euclidean derivatives by Lie derivatives and the Euclidean norm $|\xi|$ by the weight $\hskip-3pt \weight{g} \hskip-3pt$:

\noindent a) \textit{Local analysis.} Both $\hskip-4pt \weight{g} \hskip-4pt$ and $d_g^\gamma$ are comparable to their Euclidean models at small distances to $e$ and Theorem A gives a satisfactory form of \eqref{Eq-HM}. The growth of derivatives around the singularity matches the sharp Euclidean estimates. As $\dim \SL = n^2 -1$, the differentiability order also nearly matches the optimal one in the Euclidean (HM). More precisely, we match optimal Euclidean order for $n$ odd and we  loose up to one derivative for $n$ even. In fact, the local form of \eqref{Eq-ThmA} can be replaced by (weaker) \emph{Sobolev conditions} of order $n^2/2 + \varepsilon$. Moreover, less regularity suffices for small $|1/p - 1/2|$ in the spirit of Calder\'on-Torchinsky \cite{CT}. In conclusion, local singularities (at the unit $e$ or anywhere else) are admissible and the regularity around them is apparently close to optimal. As far as we know, there is no result in the literature |including rigidity conditions| which gives any information on the local behavior of Fourier $L_p$-multipliers. Our local conditions are much more flexible than in \cite{A,CSt,ST}, due to the necessary analyticity there. 

\noindent b) \textit{Asymptotic analysis.} Condition \eqref{Eq-ThmA} coincides with \eqref{Eq-HM} for the given weight up to order $[n^2/2]+1$. As stated, it apparently poses a contradiction with the asymptotic rigidity in \eqref{Eq-Rig-p}. However, this is not the case in the weight we work with, since high order H\"ormander-Mikhlin conditions impose the same decay rates for lower order terms. More precisely, the exponential growth of $L$ with respect to the word length implies the following inequality for any $\phi \in \mathcal{C}^1(\SL \setminus \{e\})$ and $\beta>2$ \vskip-12pt
$$\sup_{\mathrm{X} \in \mathfrak{sl}_n(\R)} L(g)^{\beta} \big| \partial_\mathrm{X} \phi (g) \big| \le 1 \ \Rightarrow \ L(g)^{\beta} |\phi(g) - \alpha| \le C_\beta \quad \mbox{for some } \alpha \in \C,$$
where the supremum runs over all unit vectors in the Lie algebra of $\SL$. This forces the symbol $m - \alpha \hskip-2pt : \SL \to \C$ to decay at the same rate $\beta_0 = [n^2/2]+1$ as the highest order derivative in \eqref{Eq-ThmA}, see Remark \ref{Rem-Linear}. The known rigidity theorems in this context \cite{dLdlS,LdlS} are no longer in conflict with it. In fact, Weyl's integration formula implies that the critical integrability index for $L$ is $[n^2/2]$, so that $L^{-\beta_0}$ is in $L_1(\SL)$ and the asymptotic part of Theorem A reduces to the local part by a simple patching argument. In particular, the asymptotic bound $L^{1 - \beta_0} / \log^{1+\delta} L$ for $m$ and its Lie derivatives suffices. Conceivable, the logarithmic factor could be removed for $n$ even, under the smaller regularity order $[(n^2-1)/2]+1 = [n^2/2]$. Our asymptotic conditions are (necessarily) more rigid than in \cite{A,CSt,ST}. 

The classical H\"ormander-Mikhlin theorem is usually proved by interpolation from the case 
$p=2$ and an endpoint $L_\infty$-BMO estimate. In our noncommutative situation, we are not aware of any good definition of BMO on $\mathcal L(\SL)$.

After Theorem A, we have no reason to believe that our asymptotic estimates are even close to optimal, our analysis just gives a satisfactory comparison with the Euclidean setting in a natural metric for $\SL$. However, Theorem B below proves that Mikhlin type conditions are also necessary up to certain regularity order for radial multipliers. Given an open interval $J \subset \mathbf{R}$ and $\alpha > 0$, let $\mathcal{C}^\alpha(J)$ be the space of functions which admit $[\alpha]$ continuous derivatives in $J$ and such that the $[\alpha]$-th derivative of $\varphi$ is H\"older continuous of order $\alpha \hskip-1pt - \hskip-1pt [\alpha]$ on every compact subset of $J$. Given $g \in \mathrm{SL}_n(\mathbf{R})$, we use normalized Hilbert-Schmidt norms $|g|^2 = \frac1n \mathrm{tr}(g^*g)$.

\begin{TheoB} Let $\mathrm{G}=\SL$. Consider a radial $\mathrm{G}$-symbol $m(g) = \varphi(|g|)$ for $n \ge 3$. Assume that the Herz-Schur multiplier $S_m$ is $S_p(L_2(\mathrm{G}))$-bounded for some $p > 2 + \frac{2}{n-2}$ so that $\alpha_0 := (n-2)/2 - (n-1)/p > 0$. Then 
\begin{itemize}
\item[a)] $\varphi$ has a limit $\varphi_\infty$ at $\infty$,
\item[b)] if $\alpha_0$ is not integer, then $\varphi\in \mathcal{C}^{\alpha_0}(1,\infty)$.
\item[c)] if $\alpha_0$ is an integer, then $\varphi \in \mathcal{C}^{\alpha}(1,\infty)$ for every $\alpha<\alpha_0$.
\end{itemize}
Moreover, the following local/asymptotic estimates hold for the function $\varphi$:

\begin{itemize}
\item[i)] If $n>4$ and $p>2 +\frac{6}{n-4}$ then
\[\big| \varphi(x) - \varphi_\infty \big| \le C \frac{\|S_m\|_{\mathcal{B}(S_p(L_2(\mathrm{G})))}}{ x^{c_0}} \quad \mbox{where} \quad c_0 = \frac{n} {\big[\frac{3p}{p- 2}\big]},\]
if $n \leq 4$ or ($n>4$ and $p<2+\frac{6}{n-4}$) then
\[\big| \varphi(x) - \varphi_\infty \big| \le C \frac{\|S_m\|_{\mathcal{B}(S_p(L_2(\mathrm{G})))}}{x^{c_0}} \quad \mbox{where} \quad c_0 = \frac{n}{2}-\frac{n(n-1)}{p(n-2)}\]
and if $p=2+\frac{6}{n-4}$ then
\[\big| \varphi(x) - \varphi_\infty \big| \le C \frac{\|S_m\|_{\mathcal{B}(S_p(L_2(\mathrm{G})))}}{x^{c}} \quad \mbox{for every} \quad c <\frac{n}{n-2} .\]
\item[ii)] For every $x > 1$ and every integer $1 \le k <\alpha_0$,
$$\big| \partial^k \varphi (x) \big| \le C \frac{\|S_m\|_{\mathcal{B}(S_p(L_2(\mathrm{G})))}}{(x - 1)^k x^{c_k}} \quad \mbox{where} \quad c_k = \frac{n}{[\frac{2k+1}{1- \frac{2}{p}}]}.$$
\item[iii)]The H\"older constant in a neighborhood of $x$ is bounded as follows. If $\alpha_0$ is not an integer, then
$$\limsup_{y \to x} \frac{|\partial^{[\alpha_0]} \varphi(x) - \partial^{[\alpha_0]} \varphi(y)|}{|x-y|^{\alpha_0 - [\alpha_0]}} \, \le \, C \frac{\|S_m\|_{\mathcal{B}(S_p(L_2(\mathrm{G})))}}{\big( (x-1) x^{\frac{n}{n-2}} \big)^{\alpha_0}}.$$ 
If $\alpha_0$ is an integer, then for every $\alpha<\alpha_0$
$$\limsup_{y \to x} \frac{|\partial^{[\alpha]} \varphi(x) - \partial^{[\alpha]} \varphi(y)|}{|x-y|^{\alpha - [\alpha]}} \, \le \, C \frac{\|S_m\|_{\mathcal{B}(S_p(L_2(\mathrm{G})))}}{\big( (x-1) x^{\frac{n}{n-2}} \big)^\alpha}.$$ 
\end{itemize}





\end{TheoB}

\vskip3pt

Theorem B gives a major strengthening of the rigidity theorems in this context \cite{H2,dLdlS,LdlS}. By transference \cite{LdlS}, it is also valid for radial completely bounded Fourier multipliers |i.e. using the cb-norm of $T_m$| improving \eqref{Eq-Rig-p} and higher dimensional forms in various ways, notably by the conditions in Theorem B ii). The range of $p$'s for which it applies  also improves the best known results \cite{dLdlS}. Its Euclidean form for radial multipliers $m(\xi) = \varphi(|\xi|)$ in dimension $d \ge 2$ is the necessary condition
\begin{equation} \label{Eq-NecRad} \tag{TT}
p > 2 + \frac{2}{d-1} \ \Rightarrow \ |\xi|^k \big| \partial^k \varphi(\xi) \big| \, \le \, C_{p,d} \hskip1pt \big\|T_m \hskip-2pt : L_p(\R^d) \to L_p(\R^d) \big\|
\end{equation}
for $k < (d-1)/2 - d/p$. Comparing \eqref{Eq-NecRad} with \eqref{Eq-HM} enlightens the structure of radial multipliers. Indeed, taking $p$ arbitrarily large, Mikhlin conditions are necessary up to order $[d/2] - 1$ and sufficient from $[d/2] + 1$ on. Condition \eqref{Eq-NecRad} has it roots in the regularity of the Hankel transform of radial $L_p$-functions, which goes back to Schoenberg \cite{Sc} and culminated in \cite{To,Tr}. It is the most satisfactory result for radial multipliers before the celebrated characterization \cite{GS, HNS}. 

Theorem B exactly reproduces \eqref{Eq-NecRad} around the singularity ($x=1$ and $\xi=0$ respectively) when Euclidean dimension $d$ is replaced by $n-1$. Thus, Theorem B confirms that the growth rate around the singularity in Theorem A is optimal for low order derivatives. \eqref{Eq-NecRad} also suggests that the range of $p$'s in Theorem B could be best possible. Asymptotic rigidity arises from the extra decay provided by $c_0, c_1, c_2,  \ldots, c_{[\alpha]}$ and it increases with the rank. Remarkably,  there exist radial multipliers satisfying \eqref{Eq-ThmA} in rank $n$ and failing rigidity for ranks $m >> n$.

As an application of the methods of proof of Theorem B, we also prove in Theorem \ref{thm:multipliers_SO(n,1)} related rigidity estimates in $\mathrm{SO}(n,1)$ for first and higher order derivatives. This result is especially satisfactory, since this group is rank 1 and weakly amenable. It is certainly a surprise that Lafforgue's methods around strong property (T) shed some light here.

\noindent \textbf{II. Structure of the proof.} H\"ormander-Mikhlin criteria for group algebras have been recently investigated in \cite{GJP,JMP1,JMP2}. The lack of finite-dimensional orthogonal representations for $\mathrm{SL}_2(\R)$ and Kazhdan's property (T) for higher ranks implies that Euclidean geometry only mirrors the geometry of $\SL$ via nonorthogonal actions, which are beyond the scope of the above mentioned papers. Here are some benchmarks which can help the reader to follow our argument:  

\noindent a) \textit{A local measurement of nonamenability.} Almost every form of transference since the pioneer contributions of Cotlar or Calder\'on involve some kind of amenability assumption. This was dodged in \cite{JMP1,JMP2} for those nonamenable groups which act orthogonally on some finite-dimensional real Hilbert space. Among others, this excludes $\SL$. Given $\Omega, \Sigma$ relatively compact open neighborhoods of the identity in a nonamenable group $\G$, we shall introduce a constant $0 < \delta_\Sigma(\Omega) < 1$ which quantifies the \lq $\Sigma$-nonamenability\rq${}$ relative to $\Omega$. When $\G = \SL$, this is closely related to the Harish-Chandra's function \cite{CHH,Oh}. Our proof starts with a complete contraction $\iota_p \hskip-2pt : L_p(\mathcal{L}(\SL)) \to S_p(L_2(\SL))$ (denoted $j_p^\phi$ in the text), which admits a key partial converse for $p \in 2\Z_+$: $$\|f\|_p \, \le_{\mathrm{cb}} \, \frac{1}{1 - \delta_\Sigma(\Omega_p)} \|\iota_p(f)\|_p \quad \mbox{when} \quad \mathrm{supp} \widehat{f} \subset \Omega.$$ Here $\Omega_p = \Omega \Omega^{-1} \cdots \Omega^{\pm1}$ with $p/2$ terms. This is valid for all groups $\G$.  

\noindent b) \textit{Local transference Fourier \hskip-2pt $\to$ \hskip-2pt Schur \hskip-2pt $\to$ \hskip-2pt Twisted multipliers.} Generalizing previous work of Bo\.zejko/Fendler, 
the isometry between Fourier and Schur $L_p$-multipliers was proved in \cite{CS,NR} for amenable groups. In fact, Schur multiplier norms are dominated by Fourier multiplier ones even for nonamenable groups. The reverse inequality remains open. At the moment, we find no reason to believe that it holds true for nonamenable groups. However, $\Omega$--supported Fourier multipliers are dominated by $\Sigma' \times \Sigma'$--truncated Schur multipliers $(\Sigma' = \Omega^{-1} \Sigma)$, up to a constant which  equals 1 for $\G$ amenable and blows up as $\Omega \to \G$ otherwise. More precisely, if $\mathrm{supp} \hskip2pt m \subset \Omega$ and $p \in 2\Z_+$, we shall prove that 
\begin{equation} \label{Eq-LT} \tag{LT}
\big\| T_m \big\|_{\mathrm{cb}(L_p(\V))} \, \le \, \frac{1}{1 - \delta_\Sigma(\Omega_p)} \big\| S_m \big\|_{\mathrm{cb}(S_p(L_2(\Sigma')))}.
\end{equation}
Next, if $\beta: \G \to \R^k$ is a cocycle map associated to a volume-preserving action $\alpha: \G \curvearrowright \R^k$ and $\RR$ is the algebra of bounded functions $f \hskip-2pt : \R^k \to \mathcal{B}(L_2(\Sigma))$, we shall lift our multiplier $m = \dot{m} \circ \beta$ to control the right hand side of \eqref{Eq-LT} by the cb-$L_p$-norm of the \lq\lq twisted Fourier multiplier\rq\rq${}$ given by $$\widetilde{T}_{\dot{m}}: \Big( f_{gh} \Big)_{\hskip-2pt g,h \in \Sigma} \mapsto \Big( T_{\dot{m}_g}(f_{gh}) \Big)_{\hskip-2pt g,h \in \Sigma} \quad \mbox{with} \quad \widehat{T_{\dot{m}_g} (f_{gh})}(\xi) = \dot{m}(\alpha_g (\xi)) \widehat{f}_{gh}(\xi).$$

\noindent c) \hskip1pt \textit{Twisted Riesz transforms for fractional Laplacians.} Calder\'on-Zygmund and other classical methods are inefficient to bound twisted Fourier multipliers. Instead, we shall adapt a key result from \cite{JMP2}, which identifies H\"ormander-Mikhlin multipliers as Littlewood-Paley averages of fractional Riesz transforms. It does not apply in the group algebra of $\SL$ because of the lack of orthogonal cocycles, but local transference opens the door to a twisted form in $\RR$ via the infinite-dimensional, but orthogonal cocycles naturally linked to fractional powers of Euclidean Laplacians. A duality argument then shows that it suffices to prove a square function inequality for twisted Riesz transforms. The behavior of these maps is highly asymmetric: the $\alpha$-twist affects the row index $g$, but not the column index $h$. The column case follows by a combination of harmonic analysis and operator space techniques, whereas the row case seems to be false in the full algebra $\RR$. At this point, the proof turns more technical. Roughly, we restrict to the image of our local embedding into $\RR$ and invert transference to rewrite the row square function in the group algebra. In this context, group inversion |which is locally smooth around the unit, so that H\"ormander-Mikhlin conditions are stable up to $\Omega$-constants| allows us to switch rows to columns. Then we crucially use a new form of the Littlewood-Paley theorem for Schur multipliers over $\SL$, which yields the local form of Theorem A. Our local argument holds for any compact $\Omega$, but only yields optimal estimates around the singularity. The constants grow too fast in terms of $\mathrm{diam} \hskip1pt (\Omega)$ and the alluded patching argument gives better asymptotic estimates. 

\noindent d) \hskip1pt \textit{On the rigidity theorem.} The proof of Theorem B starts with the same idea which already worked in previous rigidity theorems \cite{dLdlS,LdlS}. Given $p$ as in the statement of Theorem B, we first prove $\mathcal{C}^\alpha$-regularity for $\mathrm{SO}(n-1)$-biinvariant $S_p$-multipliers on $\mathrm{SO}(n)$. This follows in turn from related estimates in the $n-1$--dimensional sphere for the averaging operator which was introduced in \cite{LdlS}. 
The main novelty in our argument is to infer global regularity and decay estimates for radial multipliers in $\SL$, as stated in Theorem B. The key point is to factorize the radial generating function $\varphi$ as a  composition $\varphi(x) = \psi_r \circ H_r(x)$ (with $r$ depending on $x$), where $\psi_r$ is smooth as a consequence of our auxiliary estimates in $\mathrm{SO}(n)$, and the derivatives of $H_r$ decay fast enough. Then, the assertion in Theorem B follows from an application of Fa\`a di Bruno's formula. Quite surprisingly, the same technique gives \lq\lq higher order" rigidity estimates in the rank 1 groups $\mathrm{SO}(n,1)$ as explained above.

\section{\bf Local measurement of nonamenability} \label{SectLocalTransf}

Let $\G$ be a locally compact unimodular group with Haar measure $\mu$ and left regular representation $\lambda$. Let $\Omega$ be a relatively compact neighborhood of the identity in $\G$ and let $\Sigma \subset \G$ be open. If $\Sigma$ is large enough (e.g. $\overline{\Omega} \subset \Sigma$), we can always find $\phi: \G \to \R_+$ in $L_2(\G)$ and a constant $0 \le \delta_\phi(\Omega) < 1$ such that 
\begin{equation}\label{cond}
\left\{\begin{array}{l}
\displaystyle \mathrm{supp} \hskip1pt \phi \subset \Sigma \\[6pt]
\displaystyle \int_\G |\phi(h)|^2 d\mu(h) = 1\\[8pt]
\displaystyle \int_\G \big| \phi(gh) - \phi(h) \big|^2 d\mu(h) \, \le \, 2 \delta_\phi(\Omega) \textrm{ for all } g \in \Omega.
\end{array}\right.
\end{equation}
Consider the best possible constant 
\begin{eqnarray*}
\delta_\Sigma(\Omega) 
& = & 1 - \sup_{\begin{subarray}{c} \|\phi\|_2 = 1 \\ \mathrm{supp} \phi \subset \Sigma \end{subarray}} \, \inf_{g \in \Omega} \big\langle \lambda(g) \phi, \phi \big\rangle_{L_2(\G)} 
\\ & = & \inf_{\begin{subarray}{c} \|\phi\|_2 = 1 \\ \mathrm{supp} \phi \subset \Sigma \end{subarray}} \, \sup_{g \in \Omega} \, \frac12 \int_\G \big| \phi(gh) - \phi(h) \big|^2 d\mu(h).
\end{eqnarray*}
\begin{lemma}The following dichotomy holds
\begin{itemize}
\item[i)] If $\G$ is amenable, then $\delta_\G(\Omega) = 0$ for all $\Omega$. 

\item[ii)] If $\G$ is nonamenable, then $\lim_{\Omega \to \G} \delta_\G(\Omega) = 1$.
\end{itemize}
\end{lemma}
\dem
Indeed, the first assertion is a well-known characterization of
amenability. To sketch the second assertion, assume that
$\delta(\Omega) < 1 - \varepsilon$ for all compact $\Omega \subset
\G$. This means that there exists $\phi_\Omega$ in the unit ball of
$L_2(\G)$ such that $\langle \lambda(g) \phi_\Omega, \phi_\Omega
\rangle \ge \varepsilon$ for all $g \in \Omega$. Taking ultraproducts of $L_2(\G)$ over $\Omega$, we obtain an element $\Phi = (\phi_\Omega)_\Omega$ of the Hilbert space $\prod_{\mathcal U} L_2(\G)$ and we construct the
closed convex hull K of $\{ g \Phi = (\lambda(g) \phi_\Omega)_\Omega:  g \in \G\}$. The convex set K is not empty by
assumption, weakly compact and invariant under the diagonal action of
$\G$. Moreover, $0$ is not in K as $\langle g \Phi, \Phi \rangle \ge
\varepsilon$. By the Ryll-Nardzewski fixed-point theorem, there must
exist a non-zero $\Psi = (\psi_\Omega)_\Omega$ in K such that $g\Psi =
\Psi$ for all $g \in \G$. Equivalently, there exists almost invariant
vectors in $L_2(\G)$ and $\G$ is amenable.
\fin
This suggests that $\delta_\Sigma(\Omega)$ measures the \lq\lq $\Sigma$-nonamenability relative to $\Omega$\rq\rq${}$ with constants $\delta_\Sigma(\Omega) \approx 0$ for $\Omega$ small and $\delta_\Sigma(\Omega) \approx 1$ for $\Omega$ large. Along this paper, we shall not use precise estimates of these constants in terms of $(\Omega,\Sigma)$, but they can be relevant for future applications.

\begin{remark} \label{Rem-Distortion}
\emph{The Iwasawa decomposition $\SL = \mathrm{KP}$ takes place with maximal compact group $\mathrm{K} = \mathrm{SO}(n)$ and parabolic part $\mathrm{P}$ formed by the subgroup of upper triangular matrices. Let $\Delta$ be the left $\mathrm{K}$-invariant extension $\Delta(kp) = \Delta(p)$ of the left-modular function in $\mathrm{P}$. The Harish-Chandra function $\Xi: \SL \to \R_+$ is then defined as $$\Xi(g) \, = \, \int_\mathrm{K} \Delta(gk)^{-\frac12} \, dk.$$
Let us write $B_R$ for the ball of radius $R$ around the identity in the pseudometric $\log L$, as defined in the Introduction. Then, the following inequalities hold for the pair $(\Omega, \Sigma) = (B_R, B_{2 n R})$ and every $q>2$ $$\exp \Big(- \frac12 \Big[ \frac{n^2}{2} \Big] R \Big) \ \lesssim \ 1 - \delta_\Sigma(\Omega) \ \approx \ \inf_{g \in \Omega} \Xi(g) \ \le \ C_q \exp \Big( - \frac1q \Big[ \frac{n^2}{2} \Big] R \Big).$$ The lower/upper bounds arise from well-known estimates \cite{HC1} of the Harish-Chandra function $\Xi$ in terms of the modular function $\Delta$. The upper bound in the equivalence above reduces to an inequality for left $\mathrm{K}$-invariant matrix coefficients $\langle \lambda(g) \phi, \phi \rangle$ \cite{CHH,Oh}. These are the most elementary Fourier multipliers which are bounded in the Fourier and group algebras of $\SL$. The lower bound requires fairly precise estimates for F\o lner sequences in the parabolic part $\mathrm{P}$. We omit the argument since these bounds will not play a significant role in this paper.}
\end{remark}

\subsection{Matrix algebras}

Let 
$$j: \V \to \mathcal{B}(L_2(\G))$$
be the canonical embedding, it satisfies for $\widehat{f} \in \mathcal{C}_c(\G)$:
\begin{eqnarray*}
j(f) \!\!\! & = & \!\!\! j \Big( \int_\G \widehat{f}(g) \lambda(g) \, d\mu(g) \Big) \\ \!\!\! & = & \!\!\! \int_{\G \times \G} \widehat{f}(gh^{-1}) e_{gh} \, d\mu(g) d\mu(h).
\end{eqnarray*}
The $e_{gh}$ stand for infinitesimal matrix units and the last integral must be understood in the weak-$*$ sense. The map $j$ is a $*$-homomorphism, but it is not trace preserving and it fails to be $L_p$-bounded. It has been further studied for amenable actions in \cite{GP}. Given $1 \le p \le \infty$ and $0 \le \theta \le 1$, define for $f\in \V\cap L_p(\V)$ 
\begin{equation} \label{Eq-jp0}
j_{p\theta}^\phi: f \mapsto \Phi_{p,1-\theta} \hskip2pt j(f) \hskip1pt \Phi_{p,\theta}
\end{equation}
where $\Phi_{p,\theta}$ is the pointwise multiplication map by $\phi(g)^{2\theta/p}$ for some $\phi$ satisfying the  conditions \eqref{cond} with constant $\delta_\phi(\Omega)$, with the convention $\phi(g)^0=1$ if $\phi(g)>0$ and $\phi(g)^0=0$ otherwise. Without loss of generality we may also assume that $\phi\in \mathcal{C}_c(\G)$. As an operator affiliated to $\mathcal{B}(L_2(\G))$ we may think of $\Phi_{p,\theta}$ as the diagonal matrix with entries $\phi(g)^{2\theta/p}$. Consider the constants $$\varepsilon_\phi(\Omega) = \frac{\delta_\phi(\Omega) - \delta_\Sigma(\Omega)}{1 - \delta_\phi(\Omega)} \quad \mbox{so that} \quad \lim_{\delta_\phi \to \delta_\Sigma} \varepsilon_\phi(\Omega) = 0.$$ Finally, when $p \in 2\Z_+$, we also define $\Omega_p = \Omega \Omega^{-1} \Omega \Omega^{-1} \cdots \Omega^{\pm 1}$ with $p/2$ terms.

\begin{lemma} \label{LocalAmenability}
The map $j_{p\theta}^\phi$ can be extended to $L_p(\V)$ such that$\hskip1pt :$
\begin{itemize}
\item[i)] $\big\| j_{p\theta}^\phi: L_p (\V) \to S_p(L_2(\G)) \big\|_{\mathrm{cb}} \le 1$ for $p \ge 2$.

\vskip5pt

\item[ii)] If in addition $\widehat{f}\in \mathcal{C}_c(\G)$ with support in $\Omega$, then we get $$\|f\|_p \, \le_{\mathrm{cb}} \, \frac{1 + \varepsilon_\phi(\Omega_p)}{1-\delta_\Sigma(\Omega_p)} \big\| j_{p0}^\phi(f) \big\|_p \quad \mbox{for} \quad p \in 2\Z_+.$$ The same statement above holds for $j_{p1}^\phi$ instead of $j_{p0}^\phi$. 
\end{itemize}
\end{lemma}

\dem As usual, we will denote maps and their continuous extensions in the same way. The first property trivially holds for $p=\infty$. Thus, by interpolation it suffices to prove it for $p=2$. Given $0 \le \theta \le 1$, define $(2/q_0,2/q_1) = (1-\theta,\theta)$ so that every $f \in L_2(\V)\cap \V$ factorizes as $f = f_0 f_1$ with $f_i\in L_{q_i}(\V)\cap \V$ and $\|f_0\|_{q_0} \|f_1\|_{q_1} = \|f\|_2$. This yields $$\big\|j_{2\theta}^\phi(f)\big\|_2 = \big\|j_{q_00}^\phi(f_0)j_{q_11}^\phi(f_1)\big\|_2 \le \big\|j_{q_00}^\phi(f_0)\big\|_{q_0} \big\|j_{q_11}^\phi(f_1)\big\|_{q_1},$$ which reduces the problem to showing that $j_{q_00}^\phi$ and $j_{q_11}^\phi$ extend to complete contractions. By symmetry in the argument and again by interpolation, we are left to study $j_{20}^\phi$, which extends to a  complete isometry since
\begin{eqnarray*}
\big\| j_{20}^\phi(f) \big\|_2^2 & = & \mathrm{tr} \Big[ \hskip1pt \Big| \Big( \phi(g) \widehat{f}(gh^{-1}) \Big)_{\hskip-2pt g,h} \Big|^2 \Big] \\ [3pt] & = & \int_{\G \times \G} \phi(g)^2  \big| \widehat{f}(gh^{-1}) \big|^2 \, d\mu(g) d\mu(h) \\ [3pt] & = & \Big( \int_\G \phi(g)^2 d\mu(g) \Big) \Big( \int_\G |\widehat{f}(h)|^2 \, d\mu(h) \Big) \ = \ \|f\|_2^2.
\end{eqnarray*}
To prove the second property, we note that for $f,f'\in L_2(\V)$ 
\begin{eqnarray*}
\lefteqn{\hskip-8pt \big\langle j_{2,1-\theta}^\phi(f), j_{2,\theta}^\phi(f') \big\rangle} \\ [3pt] \!\!\! & = & \!\!\! \mathrm{tr} \Big[ \Big( \phi(g)^\theta \widehat{f}(gh^{-1}) \phi(h)^{1-\theta} \Big)_{\hskip-2pt g,h}^* \Big( \phi(g)^{1-\theta} \widehat{f'}(gh^{-1}) \phi(h)^\theta \Big)_{\hskip-2pt g,h} \Big] \\ \!\!\! & = & \!\!\! \int_{\G \times \G} \phi(g) \overline{\widehat{f}(gh^{-1})} \widehat{f'}(gh^{-1}) \phi(h) \, d\mu(g) d\mu(h) = \int_\G \overline{\widehat{f}(g)} \widehat{f'}(g) (1 + a_g) \, d\mu(g)  
\end{eqnarray*}
with $\displaystyle a_g = \int_\G \phi(gh) \phi(h) \, d\mu(h) - 1$ . If $g \in \Omega$
\begin{eqnarray*}
|a_g| & = & \frac12 \Big| \int_\G \big( \phi(gh)^2 + \phi(h)^2 - 2 \phi(gh) \phi(h) \big) \, d\mu(h) \Big| \ \le \ \delta_\phi(\Omega).
\end{eqnarray*}
Therefore, given $f \in L_2(\V)$ with $\widehat{f}(g)=0$ for $g \notin \Omega$, we may prove the second property as follows. Let $f'=f/\|f\|_2$ so that $\|f\|_2 = \langle f, f' \rangle$. Then we get 
\begin{eqnarray*}
\|f\|_2 & \le & \Big| \int_\G \overline{\widehat{f}(g)} \widehat{f'}(g) (1+a_g) d\mu(g) \Big| + \Big| \int_\Omega \overline{\widehat{f}(g)} \widehat{f'}(g) a_g \, d\mu(g) \Big| \\ & \le & \big| \big\langle j_{2,1-\theta}^\phi(f), j_{2,\theta}^\phi(f') \big\rangle \big| + \Big( \int_\Omega |\widehat{f'}(g)|^2 |a_g|^2 \, d\mu(g) \Big)^\frac12 \|f\|_2 \\ [5pt] & \le & \big\|j_{2,1-\theta}^\phi(f)\big\|_2 + \delta_\phi(\Omega) \|f\|_2. 
\end{eqnarray*}
This proves that 
\begin{equation} \label{Eq-ThetaL2}
\|f\|_2 \le \frac{1}{1-\delta_\phi(\Omega)} \big\|j_{2\theta}^\phi(f)\big\|_2 = \frac{1 + \varepsilon_\phi(\Omega)}{1-\delta_\Sigma(\Omega)} \big\|j_{2\theta}^\phi(f)\big\|_2 \quad \mbox{for} \quad 0 \le \theta \le 1.
\end{equation}
Now let $\widehat f\in \mathcal{C}_c(\G)$ such that $f\in L_p(\V)$ for some $p \in 2\Z_+$ greater than $2$. Consider its polar decomposition $f = u_f |f| = |f^*| u_f$ and let $q > 2$ satisfying $2/p + 2/q = 1$. Then we pick $$f' \, = \, \left\{ \begin{array}{ll} \displaystyle \frac{f^*|f^*|^{\frac{p}{2}-2}}{\|f^*|f^*|^{\frac{p}{2}-2}\|_q} & \mbox{if }\ p \in 4\Z_+, \\ [15pt] \displaystyle \frac{f^*|f^*|^{\frac{p}{2}-2} u_f}{\|f^*|f^*|^{\frac{p}{2}-2} u_f\|_q} & \mbox{if } \ p \notin 4\Z_+. \end{array} \right.$$
This gives an element of the unit ball of $L_q(\V)\cap \V$ satisfying $\|f\|_p = \|ff'\|_2$. On the other hand, note that the support of the Fourier spectrum of $ff'$ is the same as that of $|f^*|^{p/2} = (ff^*)^{p/4}$ for $p \in 4\Z_+$ and $|f^*|^{p/2-1}f = (ff^*)^{(p-2)/4}f$ when $p \in 2\Z_+ \setminus 4\Z_+$. In both cases, the Fourier spectrum is supported in the set $\Omega_p$ defined before the statement. Moreover, taking $\theta = 1-2/p$ and applying the $L_2$-inequality above we find 
\begin{eqnarray*}
\|f\|_p \ = \ \|ff'\|_2 & \le & \frac{1 + \varepsilon_\phi(\Omega_p)}{1-\delta_\Sigma(\Omega_p)} \big\|j_{2\theta}^\phi(ff')\big\|_2 \\ & = & \frac{1 + \varepsilon_\phi(\Omega_p)}{1-\delta_\Sigma(\Omega_p)} \big\|j_{p0}^\phi(f) j_{q1}^\phi(f')\big\|_2 \ \le \ \frac{1 + \varepsilon_\phi(\Omega_p)}{1-\delta_\Sigma(\Omega_p)}\big\|j_{p0}^\phi(f)\big\|_p
\end{eqnarray*}
by H\"older's inequality and the contractivity of $j_{q1}^\phi$. The cb-analogue is similar. \fin

\subsection{Local transference} 

Given $m \hskip-2pt : \G \to \C$, let $S_m$ denote the Schur multiplier associated to $m$, as defined in the Introduction. The main results in \cite{CS,NR} give that Fourier and Schur multipliers share the same cb-$L_p$-norm for amenable groups. In the nonamenable setting, Schur multiplier cb-norms are always dominated by Fourier multiplier cb-norms. The reverse inequality remains open. If we consider $\Sigma' = \Omega^{-1} \Sigma$, the distortion constants $\delta_\Sigma(\Omega_p)$ provide a local form of such an inequality between $\Omega$--supported Fourier multipliers and the corresponding Schur multiplier restricted to matrices supported in $\Sigma' \times \Sigma'$. We write
$$S_m(a) = \Big( m(gh^{-1}) a_{g,h} \Big)_{\hskip-2pt g,h \in \Sigma'} \quad \mbox{for} \quad a = \Big( a_{g,h} \Big)_{\hskip-2pt g,h \in \Sigma'}.$$ 
In all the estimates, we will always specify on which matrices $S_m$ is acting.

\begin{theorem} \label{Thm-LocalTransf}
Let $\G$ be a locally compact unimodular group. Consider a relatively compact neighborhood of the identity $\Omega$ and $p \in 2\Z_+$. Let $\Sigma$ be an arbitrary open subset in $\G$ containing $\overline{\Omega}_p$. Then, the following inequality holds for any bounded symbol $m \hskip-2pt : \G \to \C$ supported in $\Omega$ $$\big\| T_m \hskip-2pt : L_p(\V) \to L_p(\V) \big\|_{\mathrm{cb}} \, \le \, \frac{1}{1-\delta_\Sigma(\Omega_p)} \big\| S_m \hskip-2pt : S_p(L_2(\Sigma')) \to S_p(L_2(\Sigma')) \big\|_{\mathrm{cb}}.$$ 
\end{theorem}

\dem As $\mathrm{supp} \hskip1pt m \subset \Omega$, Lemma \ref{LocalAmenability} ii) gives for $f\in L_p(\mathcal L(\G))\cap \mathcal L(\G)$
\begin{equation} \label{Eq-FourierSchur} 
\|T_mf\|_p \le_{\mathrm{cb}} \frac{1+\varepsilon_\phi(\Omega_p)}{1-\delta_\Sigma(\Omega_p)} \big\| S_m(j_{p0}^\phi f) \big\|_{S_p(L_2(\Sigma'))}
\end{equation}
since $j_{p0}^\phi (T_m f) = S_m(j_{p0}^\phi f)$. Note that $\phi(g)^{2/p} m(gh^{-1})$ appears as a factor of the $(g,h)$-entry. In particular, since $\mathrm{supp} \hskip1pt (\phi \otimes m) \subset \Sigma \times \Omega$, we easily deduce that each nonvanishing entry $(g,h) \in \Sigma' \times \Sigma'$. This implies that $S_m \circ j_{p0}^\phi = S_m^{\Sigma'} \circ j_{p0}^\phi$ and Lemma \ref{LocalAmenability} i) yields the assertion for $p \ge 2$ by taking $\varepsilon_\phi$ arbitrarily small. \fin 

\begin{remark} \label{Rem-Sigma'}
\emph{The proof actually gives $$\big\| T_m \big\|_{\mathrm{cb}(L_p(\V))} \, \le \, \frac{1}{1-\delta_\Sigma(\Omega_p)} \big\| S_m \big\|_{\mathrm{cb}(S_p(L_2(\Sigma'),L_2(\Sigma)))}.$$ In other words, the same inequality holds with the $\Sigma \times \Sigma'$--truncation instead.}
\end{remark}

\begin{problem}
\emph{The above result yields the expected inequality with constant 1 for amenable groups \cite{CS,NR}, since $\delta_\G(\Omega) = 0$ for all $\Omega$. Our inequalities are new for nonamenable groups, but unfortunately leave several questions unsolved. Can we generalize Lemma \ref{LocalAmenability} and Theorem \ref{Thm-LocalTransf} to noninteger values $p \ge 2$? What is the behavior for $p<2$? On the other hand, we do not recover the $L_2$-isometry and Bo\.zejko-Fendler's $L_\infty$-isometry \cite{BF}. In other words, we should expect our constants in Theorem \ref{Thm-LocalTransf} to be close to $1$ when $p$ approaches $2$ or $\infty$. A quick review of our argument shows that this would be the case if the constants in \eqref{Eq-ThetaL2} converge fast enough to $1$ when $\theta$ approaches $0$ or $1$. Is that true? Last but not least, is there a nonlocal upper bound for nonamenable groups? Proving such an inequality or providing a counterexample would enlighten very much this relation. We have not investigated similar inequalities for nonamenable nonunimodular groups.}
\end{problem}

\subsection{Nonorthogonal cocycles}

Let $\alpha \hskip-2pt : \G \to \SLk$ be any volume-preserving continuous representation and let $\beta: \G \to \R^d$ be a continuous map satisfying the cocycle law $\alpha_g(\beta(h)) = \beta(gh) - \beta(g)$ for $g,h \in \G$. We call $\beta$ orthogonal when the action $\alpha$ is orthogonal. In contrast with \cite{JMP1,JMP2}, nonorthogonal cocycles will be admissible in what follows. A $\beta$-lifted multiplier for the symbol $m \hskip-2pt : \G \to \C$  is any function $\dot{m}: \R^d \to \C$ satisfying the identity $m = \dot{m} \circ \beta$. We need some additional notation:
\begin{itemize}
\item $\RR = L_\infty(\R^d) \bar\otimes \mathcal{B}(L_2(\Sigma))$.

\vskip5pt

\item $\dot{m}_g(\xi) = \dot{m}(\alpha_g(\xi))$ for $g \in \G$ and $\xi \in \R^d$.  

\vskip2pt

\item  $\widetilde{T}_{\dot{m}}(f)(x) = \Big( T_{\dot{m}_g}(f_{gh})(x) \Big)_{\hskip-2pt g,h \in \Sigma}$ when $f = \Big( f_{gh} \Big)_{g,h \in \Sigma}$. 

\item Characters in $\R^d$: $\exp_\xi(x) = \exp(2 \pi i \langle x, \xi \rangle)$ with $\xi \in \R^d$.
\end{itemize}
The map $\widetilde{T}_{\dot{m}}$ is a \lq\lq twisted Fourier multiplier\rq\rq${}$ which acts over $\Sigma \times \Sigma$ matrix-valued functions in $\R^d$, see Section \ref{SectTwisted} for further details. Twisted Fourier multipliers are very relevant for this paper and will be studied later on. In the following result, we adapt the arguments from \cite{CS,NR} in conjunction with de Leeuw's transference \cite{dL} to relate Schur and twisted multipliers. 

\begin{proposition} \label{PropTransference}
If $p \ge 1$, we obtain
$$\big\| S_{m} \hskip-2pt : S_p(L_2(\Sigma)) \to S_p(L_2(\Sigma)) \big\|_{\mathrm{cb}} \, \le \, \ \big\| \widetilde{T}_{\dot{m}} \hskip-2pt : L_p(\RR) \to L_p(\RR) \big\|_{\mathrm{cb}}$$ 
for any continuous Fourier symbol $\dot m \hskip-2pt : \R^d \to \C$ satisfying the identity $m = \dot{m}
\circ \beta$.
\end{proposition}

\dem Consider the unitary $u$ on $L_2(\R^d \times \Sigma)$ $$uf(\xi,g) = \exp_{\beta(g^{-1})}(\xi) f(\xi,g) = e^{2 \pi i \langle \beta(g^{-1}),\xi\rangle} f(\xi,g).$$ Define $\pi_u \hskip-2pt : S_p(L_2(\Sigma)) \to L_\infty(\R^d;S_p((L_2(\Sigma)))$ by 
\begin{eqnarray*}
\pi_u(a) \!\!\! & = & \!\!\! u^*(1\otimes a)u 
\\ \!\!\! & = & \!\!\! \Big( \exp_{-\beta(g^{-1})} a_{gh} \exp_{\beta(h^{-1})} \Big)_{\hskip-2pt g,h \in \Sigma} 
= \Big( \exp_{\alpha_{g^{-1}}(\beta (gh^{-1}))}
a_{gh} \Big)_{\hskip-2pt g,h \in \Sigma}.
\end{eqnarray*}
The last identity follows from the cocycle law for $\beta$. For every $a \in S_2(L_2(\Sigma))$, note the intertwining identity 
\begin{eqnarray*}
\widetilde{T}_{\dot{m}} (\pi_u(a)) \!\!\! & = & \!\!\! \Big( T_{\dot{m}_g}(\exp_{\alpha_{g^{-1}}(\beta (gh^{-1}))}) a_{gh}\Big)_{\hskip-2pt g,h\in   \Sigma} 
\\ \!\!\! & = & \!\!\! \Big(\dot{m}(\beta(gh^{-1}))\exp_{\alpha_{g^{-1}}(\beta   (gh^{-1}))} a_{gh}\Big)_{\hskip-2pt g,h\in \Sigma} = \pi_u(S_m(a)).
\end{eqnarray*}
Since $\pi_u$ is a normal representation of $\mathcal B(L_2(\Sigma))$ when $p= \infty$, one concludes easily by weak-$*$ density. The case $p<\infty$ requires some normalization in the spirit of de Leeuw. If $\gamma_s(x)= \exp(- {s |x|^2})$ for $s>0$, we get $$\|\gamma_{s/p}\|_{L_p(\R^d)}=\Big(\frac {\pi}{s}\Big)^{d/2p}.$$
If $\frac 1 p+ \frac 1 q=1$ and $a,b \in S_2(L_2(\Sigma))$, we claim that
\begin{equation} \label{Eq-KeydL}
\lim_{\varepsilon \to 0} \frac 1{\| \gamma_{\varepsilon/p}\|_p \| \gamma_{\varepsilon/q}\|_q} 
\Big\langle \widetilde{T}_{\dot{m}} (\gamma_{\varepsilon/p}\pi_u(a)), \gamma_{\varepsilon/q} \pi_u(b) \Big\rangle 
\, = \, \big\langle S_m(a), b \big\rangle.
\end{equation} 
Before justifying the claim, note that it implies the statement since 
\begin{equation} \label{Eq-GammaPi}
\big\| \gamma_s \pi_u(a) \big\|_{L_p(\RR)}^p = \int_{\R^d} |\gamma_s(x)|^p \big\| \pi_u(a)(x) \big\|_{S_p(L_2(\Sigma))}^p dx = \| \gamma_s \|_p^p \big\| a \big\|_{S_p(L_2(\Sigma))}^p.
\end{equation}
Therefore, as the same identity holds for $\gamma_s \pi_u(b)$ in $L_q(\RR)$, we obtain
$$\big| \big\langle S_m(a), b \big\rangle \big| \, \le \, \big\| \widetilde{T}_{\dot{m}} \hskip-2pt : L_p(\RR) \to L_p(\RR) \big\|_{\mathrm{cb}} \|a\|_p \|b\|_q.$$ This is enough to conclude by density of $S_2$ in $S_p$ and $S_q$. To prove \eqref{Eq-KeydL} we use the Plancherel formula. More precisely, by approximation we may
also assume that $a$ and $b$ have continuous kernels and
\begin{eqnarray*}
\lefteqn{\hskip-25pt \Big\langle \widetilde{T}_{\dot{m}} (\gamma_{\varepsilon/p}\pi_u(a)),
\gamma_{\varepsilon/q} \pi_u(b) \Big\rangle} \\ \hskip25pt \!\! & = & \!\! \int_{\Sigma \times \Sigma} \big\langle \dot{m}_g \widehat{[\gamma_{\varepsilon/p} \pi_u(a)_{gh}]}, \widehat{[\gamma_{\varepsilon/q}\pi_u(b)_{gh}]} \big\rangle \, d\mu(g) d\mu(h).
\end{eqnarray*}
If we set $\eta_{gh} = \alpha_{g^{-1}}(\beta(gh^{-1}))$, we note that $$\widehat{[\gamma_{\varepsilon/p} \pi_u(a)_{gh}]}(\xi) \, = \, a_{gh} \widehat{[\gamma_{\varepsilon/p}
\exp_{\eta_{gh}}]}(\xi) \, = \, a_{gh} \underbrace{\Big(\frac{p\pi}{\varepsilon}\Big)^{\frac{d}{2}} \exp\Big(-\frac
{p\pi^2 |\eta_{gh} - \xi|^2}{\varepsilon}\Big)}_{\psi_{p,\varepsilon} (\eta_{gh} - \xi)}.$$ Thus, the left hand term in \eqref{Eq-KeydL} can be rewritten as follows
$$\lim_{\varepsilon \to 0} \Big( \frac{\varepsilon}{\pi} \Big)^{\frac{d}{2}} \int_{\Sigma\times \Sigma} a_{gh}\overline{b_{gh}} \Big( \int_{\R^d} \dot{m}_g(\xi) \psi_{p,\varepsilon}(\eta_{gh} - \xi)\psi_{q,\varepsilon}(\eta_{gh} - \xi)
 d\xi \Big) d\mu(g)d\mu(h).$$
Next, using the change of variables $\xi = \varepsilon^{\frac12} \rho + \eta_{gh}$, the inner integral equals
$$\Big( \frac{pq\pi^2}{\varepsilon} \Big)^{\frac{d}{2}} \int_{\R^d} \dot{m}_g(\varepsilon^{\frac12} \rho + \eta_{gh}) \exp\big(-{\pi^2(p+q) |\rho|^2}\big) \, d\rho.$$ Altogether, the dominated convergence theorem gives the desired identity since
\begin{eqnarray*}
\lefteqn{\lim_{\varepsilon \to 0} \frac 1{\| \gamma_{\varepsilon/p}\|_p \| \gamma_{\varepsilon/q}\|_q} 
\Big\langle \widetilde{T}_{\dot{m}} (\gamma_{\varepsilon/p}\pi_u(a)), \gamma_{\varepsilon/q} \pi_u(b) \Big\rangle} 
\\ \!\!\! & = & \!\!\! \big( \pi p q \big)^{\frac{d}{2}} \Big( \int_{\R^d} e^{-\pi^2(p+q) |\rho|^2} \, d\rho \Big)  \Big( \int_{\Sigma\times \Sigma} \dot{m}_g(\eta_{gh}) a_{gh}\overline{b_{gh}} \, d\mu(g)d\mu(h) \Big)
\\ \!\!\! & = & \!\!\! \Big( \frac{p q}{p+q} \Big)^{\frac{d}{2}} \int_{\Sigma\times \Sigma} \dot{m}(\alpha_g(\eta_{gh})) a_{gh}\overline{b_{gh}} \, d\mu(g)d\mu(h) \, = \, \big\langle S_m(a), b \big\rangle. 
\end{eqnarray*}
We have used $p+q=pq$. The proof for the complete boundedness is identical. \fin

\begin{remark} \label{Rem-ShiftofTwist}
\emph{There are many variants of Proposition \ref{PropTransference}. The following will be of particular interest below. Namely, using the same argument one can prove for $\dot{m} \hskip-2pt : \R^d \to \C$ continuous that the (non-Herz) Schur multiplier with symbol $\widetilde{m}(g,h) = \dot{m}(\alpha_{g^{-1}}(\beta(gh^{-1})))$ satisfies the inequality $$\big\| S_{\widetilde{m}} \hskip-2pt : S_p(L_2(\G)) \to S_p(L_2(\G)) \big\|_{\mathrm{cb}} \, \le \, \big\| T_{\dot{m}} \hskip-2pt : L_p(\R^d) \to L_p(\R^d) \big\|_{\mathrm{cb}}.$$ This is the analogous result moving the twist from the Fourier to the Schur side. Just note that $T_{\dot m} \circ \pi_u = \pi_u \circ S_{\widetilde{m}}$ and \eqref{Eq-KeydL} still holds with the same replacements.}
\end{remark}

\begin{corollary} \label{Cor-FullTransf} 
Let $\G$ be a locally compact unimodular group. Consider a relatively compact neighborhood of the identity $\Omega$ and $p \in 2\Z_+$. Let $\Sigma$ be an arbitrary open subset containing $\overline \Omega_p$. Let $m \hskip-2pt : \G \to \C$ be a continuous symbol supported by $\Omega$ and satisfying $m = \dot{m} \circ \beta$. Then the following inequality holds $$\big\| T_m \hskip-2pt : L_p(\V) \to L_p(\V) \big\|_{\mathrm{cb}} \, \le \, \frac{1}{1-\delta_\Sigma(\Omega_p)} \big\| \widetilde{T}_{\dot{m}} \hskip-2pt : L_p(\RR) \to L_p(\RR) \big\|_{\mathrm{cb}}.$$ 
\end{corollary}

\dem The inequality in the statement with $\mathcal{R}_{\Sigma'}$ in place of $\RR$ immediately follows from Theorem \ref{Thm-LocalTransf} and Proposition \ref{PropTransference}. According to Remark \ref{Rem-Sigma'}, we may replace $S_p(L_2(\Sigma'))$ by $S_p(L_2(\Sigma'), L_2(\Sigma))$. Then, the argument in Proposition \ref{PropTransference} still applies for rectangular matrices $a,b \in S_2(L_2(\Sigma'), L_2(\Sigma))$. Indeed, the left term in \eqref{Eq-KeydL} is then written in terms of an integral over $\Sigma \times \Sigma'$ instead. This yields an upper bound which equals the cb-$L_p$-norm of the twisted Fourier multiplier over the subspace of $\Sigma \times \Sigma'$ matrix-valued functions. However, the twist of our multiplier only affects the $g$-variable and acts trivially in the $h$-variable. Thus, the cb-norm is unaffected after replacing $\Sigma'$ by $\Sigma$. This completes the proof. \fin

\begin{remark}
\emph{The lifted symbol $\dot{m}$ is called regulated when it satisfies the  Lebesgue differentiation theorem everywhere. Both Proposition \ref{PropTransference} and Corollary \ref{Cor-FullTransf} hold for regulated symbols as well, for the same reasons as in \cite{dL}.}
\end{remark}

\section{\bf Twisted multipliers} \label{SectTwisted}

Along this paper, a twisted Fourier multiplier will be any linear map sending a matrix-valued function $f$ to its Schur product $M \bullet f$ with a matrix $M$ of Fourier multipliers. Unless specified otherwise, we will work with matrices defined over a continuous parameter $(g,h)$ in $\Sigma \times \Sigma$, and the function 
\[(g,h,\xi) \in \Sigma \times \Sigma \times \R^d \mapsto M_{gh}(\xi) \in \C\] will be at least assumed to be bounded and measurable. More precisely, if $f = ( f_{gh} ) : \R^{d} \to S_2(L_2(\Sigma))$ we consider the map which sends $f$ to $$\widetilde{T}_M f = M \bullet f = \Big( T_{M_{gh}}(f_{gh}) \Big)_{\hskip-2pt g,h \in \Sigma} = \int_{\R^{d}} \Big( M_{gh}(\xi) \widehat{f}_{gh}(\xi) \Big)_{\hskip-2pt g,h \in \Sigma} \exp_\xi d\xi.$$ This determines the action of $\widetilde{T}_M$ on $L_p(\RR)$ for $p < \infty$. Typically, $M_{gh}$ will be a continuous deformation |independent of $h$, so that $(g,\xi) \mapsto M_g(\xi)$ is a bounded measurable function| of a bounded symbol $M \hskip-2pt : \R^{d} \to \C$. In that case, we may also define the twist action on $\RR$ by $(\widetilde{T}_M f) \zeta =  \widetilde{V}_M(f \zeta)$ for the map on $L_2(\Sigma\times \R^d)$ $$(\widetilde{V}_M h)_g = T_{M_g}(h_g).$$ Let $m = \dot{m} \circ \beta$ for some cocycle $\beta: \G \to \R^d$ associated to the volume-preserving action $\alpha$. In Corollary~\ref{Cor-FullTransf} when $\mathrm{supp} \, m \subset \Omega$, we have related the $L_p$-cb-boundedness of the Fourier multiplier in $L_p(\V)$ associated to $m$ with the $L_p$-cb-boundedness of the $\Sigma \times \Sigma$ twisted multiplier $$M_{gh}(\xi) = \dot{m}_g(\xi) = \dot{m}(\alpha_g(\xi)).$$ We shall work with relatively compact $\Sigma$'s in what follows. Euclidean harmonic analysis appears to be insufficient to establish $L_p$-bounds for this operator |even for small $\Sigma$ and smooth $\alpha$| and noncommutative harmonic analysis tools become necessary. A key result in \cite{JMP2} establishes that H\"ormander-Mikhlin multipliers are Littlewood-Paley averages of families of Riesz transforms associated to a fractional Laplacian. The lack of finite-dimensional orthogonal cocycles for $\SL$ forces us to express our twisted multipliers in terms of \lq\lq twisted Riesz transforms\rq\rq${}$ for fractional Laplacians. We start with some partial estimates based on these methods.

\subsection{Fractional Riesz transforms} 

The paper \cite{JMP2} revisited the classical theory of fractional Riesz transforms to see 
them as particular cases of noncommutative directional Riesz transform. We shall use their terminology, which may seem odd to readers familiar with integro-differential operators. Given $0 < \varepsilon < 1$, let 
\begin{equation} \label{Eq-Heps}
\H_\varepsilon = L_2(\R^d, \mu_\varepsilon) \quad \mbox{with} \quad d\mu_\varepsilon(x) = \frac{dx}{|x|^{d + 2\varepsilon}}.
\end{equation}
The map $b_\varepsilon: \R^d \to \H_\varepsilon$ given by $b_\varepsilon(\xi) = \exp(2 \pi i \langle \xi, \cdot \rangle) -1$ is an orthogonal cocycle with respect to the action $\alpha_{\varepsilon,\xi}(f) = \exp(2\pi i \langle \xi, \cdot \rangle) f$. The length function associated to it $\psi_\varepsilon(\xi) = \langle b_\varepsilon(\xi), b_\varepsilon(\xi) \rangle_{\H_\varepsilon}$ is the fractional Laplacian length $$\psi_{\varepsilon}(\xi) = 2 \int_{\R^d} \big( 1 - \cos (2 \pi \langle \xi, x \rangle) \big) \, \frac{dx}{|x|^{d+2\varepsilon}} = \mathrm{c}_{d,\varepsilon} |\xi|^{2 \varepsilon} \quad \mbox{with} \quad \mathrm{c}_{d,\varepsilon} \approx \frac{\pi^{d/2}}{\Gamma(d/2)} \frac{1}{\varepsilon(1-\varepsilon)},$$ see \cite[Example 1.4.A]{JMP2}. As usual, define 
\[ \mathsf{H}_\alpha^2(\R^d) = \Big\{ f \hskip-2pt : \R^d \to \C: \, \big\| (1 + | \ |^2)^{\frac{\alpha}{2}} \widehat{f} \, \big\|_{L_2(\R^d)} < \infty \Big\}.\] Define also $\mathrm{W}_{d,\varepsilon}^2(\R^d)$ as the completion of $\mathcal{C}_\mathrm{c}^\infty(\R^{d} \setminus \{0\})$ for the norm
\[ \|f\|_{\mathrm{W}_{d,\varepsilon}^2(\R^d)} := \big\| | \hskip3pt |^{\frac{d}{2} + \varepsilon} \widehat{ \sqrt{\psi_\varepsilon} f } \big\|_{L_2(\R^d)}.\]

\begin{remark} \label{Rem-DilationW}
\emph{The norm in $\mathrm{W}_{d,\varepsilon}^2(\R^d)$ is dilation invariant $$\| f \|_{\mathrm{W}_{d,\varepsilon}^2(\R^d)} = \big\| f(\lambda \hskip1pt \cdot) \big\|_{\mathrm{W}_{d,\varepsilon}^2(\R^d)} \quad \mbox{for all} \quad \lambda > 0.$$}
\end{remark}

We shall use several times along the paper Littlewood-Paley partitions of unity. Given a smooth function $\eta: \R^d \to \R_+$ with $\chi_{\mathrm{B}_1(0)} \le \eta \le \chi_{\mathrm{B}_2(0)}$, the model of such a partition of unity is 
\begin{equation} \label{Eq-LP}
\varphi_j(\xi) = \Big( \eta(2^{-j} \xi) - \eta(2^{1-j} \xi) \Big)^\frac12 \quad \mbox{with} \quad j \in \Z.
\end{equation}

\begin{lemma} \label{SobVsClassical}
If $0 < \varepsilon \le 1 - \frac{d}{2} + [\frac{d}{2}]$
\begin{eqnarray*}
\sup_{j \in \Z} \big\| \varphi_j^2 \hskip1pt M \big\|_{\mathrm{W}_{d,\varepsilon}^2(\R^d)}
& \lesssim & \sup_{j \in \Z} \big\| \varphi_0^2 \hskip1pt M(2^j \cdot ) \big\|_{\mathsf{H}_{\frac{d}{2}+\varepsilon}^2(\R^d)} \\
& \lesssim & \max_{0 \le |\gamma| \le [\frac{d}{2}] + 1} \esssup_{\xi \in \R^d} \ |\xi|^{|\gamma|} \big| \partial_\xi^\gamma M(\xi) \big|.
\end{eqnarray*} 
\end{lemma}

The second inequality in the statement is standard \cite[Theorem 5.2.7]{Gra}. The first inequality follows since $\varphi_j(\xi) = \varphi_0(2^{-j}\xi)$ and $\mathrm{W}_{d,\varepsilon}^2(\R^d)$ has a dilation invariant norm. The inequality then reduces to showing that $$f \in \mathsf{H}_{\frac{d}{2} + \varepsilon}^2(\R^d) \mapsto | \ |^\varepsilon \psi f \in \ \mathsf{H}_{\frac{d}{2}+\varepsilon}^2(\R^d)$$ is bounded for a smoothing $\psi$ of the characteristic function of $\mathrm{supp} \hskip1pt \varphi_0$ vanishing around $0$. This easily follows by complex interpolation using endpoint spaces $\mathsf{H}_\alpha^2$ with $\alpha \in 2\Z_+$. As noted in \cite{JMP2}, the main advantages of working with this new Sobolev norm (instead of the classical space $\mathsf{H}_\alpha^2(\R^d)$) come from dilation invariance in Remark \ref{Rem-DilationW} and, especially, from the crucial result below.

\begin{lemma} \label{HMLemma}
The following map is a unitary
\[ \Psi_\varepsilon: \H_\varepsilon = L_2(\R^d,\mu_\varepsilon) \ni h \mapsto \Bigg( \xi \mapsto \frac{1}{\sqrt{\psi_\varepsilon(\xi)}} \int_{\R^d} b_\varepsilon(\xi) h d\mu_{\varepsilon} \Bigg) \in \mathrm{W}_{d,\varepsilon}^2(\R^d) .\]
Moreover, the space $\mathrm{W}_{d,\varepsilon}^2(\R^d)$ consists of bounded continuous functions on $\R^d \setminus \{0\}$.
\end{lemma}

\dem A first useful observation is that $\xi \in \R^d \mapsto b_\varepsilon(\xi) \in \H_\varepsilon$ is continuous, and therefore $\xi \in \R^d \setminus\{0\} \mapsto b_\varepsilon(\xi)/\sqrt{\psi_\varepsilon(\xi)} \in \H_\varepsilon$ is bounded by $1$ and continuous. This implies that for every $h \in \H_\varepsilon$, the function $\Psi_\varepsilon(h)$ is a bounded continuous function on $\R^d \setminus \{0\}$, and $\|\Psi_\varepsilon(h)\|_\infty \leq \|h\|_{\H_\varepsilon}$.

We proceed by constructing first the inverse map. Set $u(m) := | \hskip3pt |^{d+2\varepsilon} \widehat{\sqrt{\psi_\varepsilon} m}$ for any function $m$ in $\mathcal{C}_\mathrm{c}^\infty(\R^{d} \setminus \{0\})$. Since $m$ vanishes around $0$, it turns out that $\sqrt{\psi_\varepsilon} m$ belongs to $\mathcal{C}_\mathrm{c}^\infty(\R^{d})$ and its Fourier transform is in the Schwartz class. This proves that $u(m)$ belongs to $\H_\varepsilon$ and we get
\[ \|u(m)\|_{\H_\varepsilon} = \Big\| \frac{1}{| \hskip3pt |^{\frac{d}{2} + \varepsilon}} u(m) \Big\|_{L_2(\R^d)} = \big\| | \hskip3pt |^{\frac{d}{2} + \varepsilon}  \widehat{\sqrt{\psi_\varepsilon} m} \big\|_{L_2(\R^d)} = \|m\|_{\mathrm{W}_{d,\varepsilon}^2(\R^d)}.\]
By density, $u$ extends to an isometry from $\mathrm{W}_{d,\varepsilon}^2(\R^d)$ into $\H_\varepsilon$. Moreover, we have
\[\sqrt{\psi_\varepsilon(\xi)} \Psi_\varepsilon(u(m))(\xi)  = \int b_\varepsilon(\xi)(x) \widehat{\sqrt{\psi_\varepsilon} m}(x) dx = \big\langle \widehat{b_\varepsilon(\xi)}, \sqrt{\psi_\varepsilon} m \big\rangle\]
where $\widehat{b_\varepsilon(\xi)} = \delta_\xi - \delta_0$ is the Fourier transform of the tempered distribution $b_\varepsilon(\xi)$. This gives $\sqrt{\psi_\varepsilon(\xi)} \Psi_\varepsilon(u(m))(\xi) = \sqrt{\psi_\varepsilon (\xi)}m(\xi)$ as the test function $\sqrt{\psi_\varepsilon}m$ vanishes at $0$. In particular, by the observation made at the beginning of the proof, $\mathrm{W}_{d,\varepsilon}^2(\R^d)$ consists of bounded continuous functions on $\R^d\setminus\{0\}$. By density, $\Psi_\varepsilon \circ u$ is the identity map on $\mathrm{W}_{d,\varepsilon}^2(\R^d)$.

All we are left to prove is that $u$ is surjective, or equivalently that the only element $h$ of $\H_\varepsilon$ orthogonal to $u(m)$ for every $m \in \mathcal{C}_\mathrm{c}^\infty(\R^{d} \setminus \{0\})$ is $h=0$. Given $h \in \H_\varepsilon$ and $m \in \mathcal{C}_\mathrm{c}^\infty(\R^{d} \setminus \{0\})$, we have
\[ \int_{\R^d} h u(m) \, d \mu_\varepsilon = \big\langle \widehat h, \sqrt{\psi_\varepsilon} m \big\rangle\]
for $\widehat h$ the Fourier transform of $h$ seen as a tempered distribution. If the preceding is $0$ for every $m \in \mathcal{C}_\mathrm{c}^\infty(\R^{d} \setminus \{0\})$, we obtain that the support of $\widehat h$ is contained in $0$ or equivalently that $h$ is a polynomial function. But the only polynomial function in $\H_\varepsilon$ is the zero polynomial. This proves that $u$ is surjective, as desired. \fin

Twisted forms of fractional Laplacian Riesz transforms will play a crucial role later in this paper. Let us introduce these maps for future reference. The Riesz transform for the $\varepsilon$-fractional Laplacian pointing towards $u \in L_2(\R^d; \mu_\varepsilon)$ is $$R_{\psi_\varepsilon,u} f(x) \, = \, \int_{\R^d}\underbrace{\frac{\langle b_\varepsilon(\xi), u \rangle_{\H_\varepsilon}}{\sqrt{\psi_\varepsilon(\xi)}}}_{\rho_{\varepsilon,u}(\xi)} \widehat{f}(\xi) \exp_\xi(x)\, d\xi = T_{\rho_{\varepsilon,u}} f(x).$$ The associated $\G$-twisted Riesz transforms are \vskip-12pt
\begin{eqnarray} 
\label{Eq-TwistRiesz}
\widetilde{R}_{\psi_\varepsilon,u}(f) & = & \Big( R_{\psi_\varepsilon,u}^g(f_{gh}) \Big)_{\hskip-2pt g,h \in \Sigma}, \\
\label{eq-SqTwistRiesz}
\widetilde{R}_{\psi_\varepsilon, \mathbf{u}} \Big( \sum_{j \in \Z} f_j \otimes \delta_j \Big) & = & \sum_{j \in \Z} \widetilde{R}_{\psi_\varepsilon,u_j}(f_j) \otimes \delta_j,\end{eqnarray}
for the multipliers $R_{\psi_\varepsilon,u}^g$ with symbol $\rho_{\varepsilon,u}(\alpha_g(\xi))$ and any family $\mathbf{u} = (u_j) \subset \H_\varepsilon$.

\subsection{Calder\'on-Zygmund methods} 

The theory of singular integral operators acting over matrix-valued functions has been recently developed to include $L_\infty$-BMO estimates \cite{JMP1,Pa1}. If $f$ is a function affiliated to $\RR = L_\infty(\R^d) \bar\otimes \mathcal{B}(L_2(\Sigma))$ and $k(x,y)$ is a linear operator on $L_2(\Sigma)$ for every $(x,y) \in \R^{2d} \setminus \Delta$, the $k$-singular integral acting on $f$ is formally given by $$Tf(x) = \int_{\R^d} k(x,y) \, (f(y)) \, dy \qquad \mbox{for} \qquad x \notin \mathrm{supp}_{\R^d} f.$$ The next result follows from \cite[Sect VI.4.4 Prop 2 (b)]{St}, see also \cite[Lemma 4.3]{Xiong}.

\begin{lemma} \label{KernelHM}
For $A \hskip-2pt : \R^d \to \H$ for some Hilbert space $\H$, we have
$$\|A\|_{L_\infty(\R^d; \H)} + \int_{|x| \ge 2 |w|} \big\| \widehat{A}(x-w) - \widehat{A}(x) \big\|_\H \, dx \, \le \, C_{d,\varepsilon} \sup_{j \in \Z} \big\| \varphi_0^2 A(2^j \cdot) \big\|_{\mathsf{H}_{\frac{d}{2}+ \varepsilon}^2(\R^k; \H)}.$$
\end{lemma}

In the case of twisted Fourier multipliers, Calder\'on-Zygmund conditions can be streamlined using Lemma \ref{KernelHM} and noncommutative techniques. Let us now present these conditions for arbitrary column-twisted multipliers 
$$\widetilde{T}_M(f)(x) \, = \, \left( \int_{\R^d} M_g(\xi) \widehat{f}_{gh}(\xi) e^{2\pi i \langle x, \xi \rangle} \, d\xi \right)_{\hskip-2pt g,h \in \Sigma}.$$ Our analysis requires  Hardy and BMO spaces over $\RR$ in the sense of \cite{Mei} that we briefly describe.

Let $\Q$ be the set of balls in $\R^d$ and for $Q \in \Q$ let $f_Q$ be the $Q$-average of a matrix-valued function $f$. We also use the standard notation of $\lambda Q$ for the ball concentric with $Q$ and radius multiplied by $\lambda$. The Hardy and BMO column-norms of $f$ are defined in \cite{Mei} as 
\begin{eqnarray} \label{Eq-HpBMO}
\|f\|_{H_p^c(\RR)} & = & \Big\| \Big( \int_{\R_+} \Big| \frac{\partial}{\partial t} P_t f \Big|^2 t \, dt \Big)^{\frac12} \Big\|_{L_p(\RR)}, \\ 
\|f\|_{\mathrm{BMO}_{\RR}^c} & = & \sup_{Q \in \Q} \Big\| \Big( \frac{1}{|Q|} \int_Q \big| f(x) - f_Q \big|^2 dx \Big)^{\frac12} \Big\|_{\mathcal{B}(L_2(\Sigma))}.
\end{eqnarray}
Here $(P_t)_{t > 0}$ denotes the Poisson semigroup 
$$P_t(f)(x) \, = \, e^{-t(-\Delta)^{\frac12}}f(x) \, = \, \Big( \int_{\R^d} e^{-t |\xi|} \widehat{f}_{gh}(\xi) e^{2\pi i \langle x,\xi \rangle} \, d\xi \Big)_{\hskip-2pt g,h \in \Sigma}.$$
 We consider the associated spaces $H_p^c(\RR)$ and $\mathrm{BMO}_{\RR}^c$.
 Their row versions are defined using adjoints:  the $\mathrm{BMO^r}$-norm is
 $\|f\|_{\mathrm{BMO}_{\RR}^r} = \|f^*\|_{\mathrm{BMO}_{\RR}^c}$. The full
$\mathrm{BMO}$-space is then $\mathrm{BMO}_{\RR} = \mathrm{BMO}_{\RR}^r \cap \mathrm{BMO}_{\RR}^c$.
 
 Theorem 6.2 in Mei's paper \cite{Mei} establishes for $p>2$
\begin{eqnarray*}
H_p^c(\RR) & \simeq & \big[ \mathrm{BMO}_{\RR}^c, L_2(\RR) \big]_{2/p}, \\ L_p(\RR) & \simeq & \big[ \mathrm{BMO}_{\RR}, L_2(\RR) \big]_{2/p}.
\end{eqnarray*}

Let us now consider the following regularity conditions

\begin{itemize}
\item[A)] \emph{Uniform Sobolev smoothness} $$\hskip20pt \sup_{g \in \Sigma} \hskip1pt \sup_{j \in \Z} \big\| \varphi_0^2 \hskip1pt M_g(2^j \cdot ) \big\|_{\mathsf{H}_{\frac{d}{2}+\varepsilon}^2(\R^d)} \, \le \, C_\mathrm{sob} \, < \, \infty.$$ 

\vskip5pt

\item[B)] \emph{Schur factorization of the twist} $$\hskip30pt M_g(\xi) \, = \, \big\langle A_g, B_\xi \big\rangle_\mathcal{K}$$ for some Hilbert $\mathcal{K}$ and measurable functions $A: \Sigma \to \mathcal{K}$ and $B:\R^d\to \mathcal{K}$ satisfying $\displaystyle \sup_{(g,\xi) \in \Sigma \times \R^d} \|A_g\|_\mathcal{K} \|B_\xi\|_{\mathcal{K}} \le C_{\mathrm{sch}} < \infty$.
\end{itemize}

\begin{lemma}\label{BMOA}
Under  condition \emph{A}, $\widetilde{T}_M: \RR \stackrel{\mathrm{cb}}{\longrightarrow} \mathrm{BMO}_{\RR}^c(\RR)$.
  \end{lemma}
\dem We shall only prove boundedness, since cb-boundedness follows from the same argument. Given $f \in \RR$ and $Q \in \Q$, set $f_{1Q} = f \chi_{5Q}$ and $f_{2Q} = f \chi_{\R^d \setminus 5Q}$. By the triangle and Kadison-Schwarz inequalities
\begin{eqnarray*}
\lefteqn{\hskip-30pt \Big\| \Big( \frac{1}{|Q|} \int_Q \big| \widetilde{T}_Mf(z) - (\widetilde{T}_Mf)_Q \big|^2 dz \Big)^{\frac12} \Big\|_{\mathcal{B}(L_2(\Sigma))}} \\ \hskip30pt & \le & 2 \Big\| \Big( \frac{1}{|Q|} \int_Q \big| \widetilde{T}_Mf_{1Q}(z) \big|^2 dz \Big)^{\frac12} \Big\|_{\mathcal{B}(L_2(\Sigma))} \\ & + & \Big\| \Big( \frac{1}{|Q|} \int_Q \big| \widetilde{T}_Mf_{2Q}(z) - (\widetilde{T}_Mf_{2Q})_Q \big|^2 dz \Big)^{\frac12} \Big\|_{\mathcal{B}(L_2(\Sigma))} \ = \ 2\mathrm{A}_{1Q} + \mathrm{A}_{2Q}.
\end{eqnarray*}
Let $\widetilde{V}_M$ be the $L_2(\R^d;L_2(\Sigma))$-bounded map $$(\widetilde{V}_M f)_g(x) = \int_{\R^d} M_g(\xi) \widehat{f}_g(\xi) e^{-2\pi i \langle x,\xi \rangle} \, d\xi.$$ $\widetilde{T}_M$ is cb-bounded on $L_2(\RR)$ and a \lq\lq right-module" since $(\widetilde{T}_Mf)\zeta = \widetilde{V}_M(f\zeta)$ for $\zeta \in L_2(\Sigma)$ |warning: it is not a left-module| so we may estimate $\mathrm{A}_{1Q}$ as follows
\begin{eqnarray*}
\mathrm{A}_{1Q} & = & \sup_{\|\zeta\|_{L_2(\Sigma)} = 1} \Big( \frac{1}{|Q|} \int_Q \big\| \widetilde{V}_M(f_{1Q}\zeta)(z) \big\|_{L_2(\Sigma)}^2 dz \Big)^{\frac12} \\
& \lesssim & \sup_{\|\zeta\|_{L_2(\Sigma)} = 1} \Big( \frac{1}{|Q|} \int_{\R^d} \big\| f_{1Q}(z) \zeta \big\|_{L_2(\Sigma)}^2 dz \Big)^{\frac12} \ \le \ 5^{\frac{d}{2}} \|f\|_{\RR}.
\end{eqnarray*}
On the other hand, Jensen's inequality gives 
$$\mathrm{A}_{2Q} \, \le \, \sup_{y,z \in Q} \big\| \widetilde{T}_Mf_{2Q}(z) - \widetilde{T}_Mf_{2Q}(y)\big\|_{\mathcal{B}(L_2(\Sigma))} \, =: \sup_{y,z \in Q} \, B_Q(y,z).$$
Recall that the kernel $K$ of $\widetilde{T}_M$ is the matrix $K(x-y) = \mathrm{diag} ( \widehat{M}_g(x-y))_{g \in \Sigma}$ acting by left matrix multiplication. Since we have $|x-z| \ge 2 |y-z|$ for all $(x,y,z) \in (\R^d \setminus 5Q) \times Q \times Q$, we may write for $f_{2Q}^\zeta := f_{2Q}\zeta$
\begin{eqnarray*}
B_Q(y,z) \!\!\! & \leq & \!\!\!\! \sup_{\|\zeta\|_{L_2(\Sigma)} = 1} \Big\| \int_{|x-z| \ge 2 |y-z|} \big( K(x-y) - K(x-z) \big) f_{2Q}^\zeta(x) \, dx \Big\|_{L_2(\Sigma)} \\
\!\!\! & = & \!\!\!\! \sup_{\begin{subarray}{c} \|\zeta\|_{L_2(\Sigma)} = 1 \\ \|D\|_{L_2(\Sigma)} = 1 \end{subarray}} \Big[ \int_{|x-z| \ge 2 |y-z|} \hskip-2pt \Big\langle \underbrace{ D, \big( K(x-y) - K(x-z) \big) f_{2Q}^\zeta(x)}_{\begin{array}{c} \parallel \\ \big\langle D_K(x,y,z), f_{2Q}^\zeta(y) \big\rangle \end{array}} \Big\rangle_{\hskip-2pt L_2(\Sigma)}  \, dx \Big],
\end{eqnarray*} 
where $D_K(x,y,z)$ is the vector in $L_2(\Sigma)$ given by $$\overline{D_K^g(x,y,z)} \, = \, \overline{D}_g \big( \widehat{M}_g(x-y) - \widehat{M}_g(x-z) \big).$$
This readily yields $$\sup_{Q \in \Q} B_Q(y,z) \, \le \, \sup_{\|D\|_{L_2(\Sigma)} = 1} \Big[ \underbrace{\int_{|x-z| \ge 2 |y-z|} \big\| D_K(x,y,z) \big\|_{L_2(\Sigma)} \, dx}_{D_K(y,z)} \Big] \|f\|_{\RR},$$ so that it remains to estimate the integral above. Let us consider the change of variables $(x-z, y-z) \mapsto (x,w)$ and consider $A(\xi) = (\overline{D}_g M_g(\xi))_{g \in \Sigma} \in L_2(\Sigma)$. Then it turns out that $D_K(y,z) = D_K(y-z,0) = D_K(w)$ coincides with the integral term in the statement of Lemma \ref{KernelHM}. Taking $B_{gj} = \overline{D}_g \varphi_0^2 M_g(2^j \cdot)$ and $\alpha = \frac{d}{2} + \varepsilon$, this yields the desired inequality
\begin{eqnarray*}
D_k(w) & \!\!\! \le \!\!\! & \sup_{j \in \Z} \big\| \varphi_0^2 A(2^j \cdot) \big\|_{\mathsf{H}_{\frac{d}{2}+ \varepsilon}^2(\R^d; \H)} \\ & \!\!\! = \!\!\! & \sup_{j \in \Z} \Big( \int_{\R^d} \int_\Sigma \Big| \int_{\R^d} \big( 1 + |\xi|^2 \big)^{\frac{\alpha}{2}} \widehat{B}_{gj}(\xi) e^{2\pi i \langle x, \xi \rangle} d\xi \Big|^2 d\mu(g) dx \Big)^\frac12 \\ & \!\!\! = \!\!\! & \sup_{j \in \Z} \Big( \int_\Sigma |D_g|^2 \big\| \varphi_0^2 M_g(2^j \cdot) \big\|_{\mathsf{H}^2_{\frac{d}{2}+\varepsilon}(\R^d)}^2 \, d\mu(g) \Big)^\frac12 \! \lesssim C_\mathrm{sob}.
\end{eqnarray*} \fin

\begin{lemma}\label{BMOAB}
  Under  condition \emph{A + B}, $\widetilde{T}_M: \RR \stackrel{\mathrm{cb}}{\longrightarrow} \mathrm{BMO}_{\RR}^r(\RR)$.
\end{lemma}
\dem The proof is variation of that of Lemma \ref{BMOA} considering $\widetilde{T}_M(f)^*$ instead of $\widetilde{T}_M(f)$ and using the stability of the conditions by complex conjugation.

Given $f \in \RR$ and $Q \in \Q$, we use the same decomposition as there
$f=f_{1Q}+f_{2Q}$. The given estimate for $\mathrm{A}_{2Q}$ in the column case
is adjoint invariant. Therefore, we just need to prove $$\sup_{Q \in \Q} \Big\| \Big( \frac{1}{|Q|} \int_Q \big| \widetilde{T}_M(f_{1Q})^*(z) \big|^2 dz \Big)^{\frac12} \Big\|_{\mathcal{B}(L_2(\Sigma))} \, \lesssim_{\mathrm{cb}} \, \|f\|_{\RR}.$$ 
Indeed, letting $\widetilde{T}_M^\dag(f) = \widetilde{T}_M(f^*)^*$, this will follow from   
$$\Big\| \Big( \int_{\R^d} \widetilde{T}_M^\dag f(x)^* \widetilde{T}_M^\dag f(x) \, dx \Big)^\frac12 \Big\|_{\mathcal{B}(L_2(\Sigma))} \, \lesssim_{\mathrm{cb}} \, \Big\| \Big( \int_{\R^d} f(x)^* f(x) \, dx \Big)^\frac12 \Big\|_{\mathcal{B}(L_2(\Sigma))}.$$ It is also straightforward to show that we have $$\widetilde{T}_M^\dag f \, = \, \Big( T_{M_h}^\dag(f_{gh})\Big) \, = \, \Big( T_{\overline{M}_h(- \hskip2pt \cdot)}^{\null} (f_{gh}) \Big).$$ Indeed, the $\dag$-operation on $T_{M_h}$ transforms the Fourier symbol $M_h(\xi)$ into $\overline{M}_h(-\xi)$ but Schur factorization is stable under this kind of transformation. In other words, we need to justify the $L_2$-column inequality 
\begin{equation} \label{Eq-Col(h)}
\Big\| \int_{\R^d} \Big| \Big( T_{M_h}(f_{gh}) \Big) \Big|^2 (x) \, dx \Big\|_{\mathcal{B}(L_2(\Sigma))}^\frac12 \, \lesssim_{\mathrm{cb}} \, \Big\| \int_{\R^d} \Big| \Big( f_{gh} \Big) \Big|^2 (x) \, dx \Big\|_{\mathcal{B}(L_2(\Sigma))}^\frac12
\end{equation}
According to Plancherel theorem, the left hand side can be written as 
\begin{eqnarray*}
\mathrm{LHS}_{\eqref{Eq-Col(h)}}^2 & = & \sup_{\|\zeta\|_{L_2(\Sigma)} = 1} \int_\Sigma \int_{\R^d} \Big| \int_\Sigma M_h(\xi) \widehat{f}_{gh}(\xi) \zeta(h) d \mu(h) \Big|^2 d\xi d\mu(g).
\end{eqnarray*}
Assume that $M_h(\xi) = \langle A_h, B_\xi \rangle_\mathcal{K}$ as in the statement. 
This yields 
$$\mathrm{LHS}_{\eqref{Eq-Col(h)}}^2 \, \le \, \sup_{\|\zeta\|_{L_2(\Sigma)} = 1} \int_\Sigma \int_{\R^d} \Big\| \int_\Sigma  \widehat{f}_{gh}(\xi) \zeta(h) A_h d \mu(h) \Big\|_\mathcal{K}^2 d\xi d\mu(g) \times \sup_{\xi \in \R^d} \|B_\xi\|_\mathcal{K}^2.$$
By Plancherel's formula,
\begin{align*} \mathrm{LHS}_{\eqref{Eq-Col(h)}}^2 &\lesssim \, \sup_{\|\zeta\|_{L_2(\Sigma)} = 1} \int_\Sigma\int_{\R^d} \Big\| \int_\Sigma f_{gh}(x) \zeta(h) A_h d \mu(h) \Big\|_\mathcal{K}^2 dx d\mu(g)\\
& = \, \sup_{\|\zeta\|_{L_2(\Sigma)} = 1} \Big\langle \zeta A, (R\otimes \mathrm{Id}_{\mathcal{K}})(\zeta A) \Big\rangle_{L_2(\Sigma;\mathcal K)}
\end{align*}
where $R=\displaystyle \int_{\R^d} \Big| \Big( f_{gh}(x) \Big) \Big|^2 dx$. The assertion then follows, since it is clear that $\|\zeta A\|_{L_2(\Sigma;K)}\| \le \sup_{h \in \Sigma} \|A_h\|_{\mathcal{K}}$. \fin

\begin{proposition} \label{CZProposition}
The following results hold$:$
\begin{itemize}
\item[i)] Condition \emph{A} implies $\widetilde{T}_M: L_p(\RR) \stackrel{\mathrm{cb}}{\longrightarrow} H_p^c(\RR)$ for $2 \le p < \infty$.

\item[ii)] Conditions \emph{A + B} imply $\widetilde{T}_M: L_p(\RR) \stackrel{\mathrm{cb}}{\longrightarrow} L_p(\RR)$ for $1 < p < \infty$.
\end{itemize}
\end{proposition}

\dem We rely on interpolation. Lemma \ref{KernelHM} readily implies that the twisted multiplier $\widetilde{T}_M$ is completely bounded in $L_2(\RR)$.

\noindent \textbf{i)} The missing endpoint estimate is exactly Lemma \ref{BMOA}.

\noindent \textbf{ii)} When $p>2$, similarly one just needs to use Lemma \ref{BMOA} and \ref{BMOAB}.

To reach $1<p<2$, we use duality. The identity $$\widetilde{T}_M^* =
\widetilde{T}_{\overline{M}}$$ holds for the duality pairing
$L_p(\R^d; S_p(L_2(\Sigma)))^* = L_q(\R^d;
S_q^\mathrm{op}(L_2(\Sigma)))$. Since conditions A/B are stable under
complex conjugation, the proof is over. \fin

\begin{remark} \label{RemConstantsCZ}
\emph{Compared to \cite{JMP1}, Proposition \ref{CZProposition} is also valid for twists not coming from orthogonal actions and Mikhlin regularity has been optimized to order $[\frac{d}{2}] + 1$.}
\end{remark}

\begin{remark} \label{RemRCCZ}
 \emph{It is possible to adapt Proposition \ref{CZProposition} and its proof to families of twisted multipliers. Since this will not be used later in the paper, we only provide a brief informal statement with no proof. To do so one first needs to introduce the usual analogs of square function norms \cite{P2}. Namely, given a family $f_j \in L_p(\M)$  
\begin{eqnarray} \label{eq-Cp}
&\Big\| \summ_j f_j \otimes e_{j1} \Big\|_{L_p(\M; C_p)} =  \Big\| \Big( \summ_j f_j^*f_j \Big)^\frac12 \Big\|_p\\
&\Big\| \summ_j f_j \otimes e_{1j} \Big\|_{L_p(\M; R_p)} =  \Big\| \Big( \summ_j f_jf_j^* \Big)^\frac12 \Big\|_p\\\label{eq-RCp}
& \hskip10pt \Big\| \summ_j f_j \otimes \delta_{j} \Big\|_{L_p(\M; RC_p)} =
\max \Bigg\{ \Big\| \Big( \summ_j f_j^*f_j \Big)^\frac12 \Big\|_p, \Big\| \Big( \summ_j f_j f_j^* \Big)^\frac12 \Big\|_p \Bigg\}.
\end{eqnarray}
The announced adaptation of Proposition \ref{CZProposition} provides sufficient conditions for the complete boundedness on $L_p(\RR;C_p)$ and $L_p(\RR;R_p)$ of maps $f\mapsto \summ_j \widetilde{T}_{M_j}(f)\otimes \delta_j$ for families of bounded measurable functions $(g,\xi)\mapsto M_{jg}(\xi)$. Condition A becomes an $\ell_2$-valued H\"ormander-Mikhlin condition for the symbols $(M_{jg})$ uniformly in $g \in \Sigma$. Condition B is more interesting. Consider $(g,j,\xi) \in \Sigma \times \Z \times \R^d$. The variable $g \in \Sigma$ is always a column variable. The variable $j \in \Z$ is a row/column according to whether we take values in $R_p$ or $C_p$. Finally the variable $\xi \in \R^d$ is a row/column according to whether we want a $\mathrm{BMO}_{\RR}^r$ or a $\mathrm{BMO}_{\RR}^c$ estimate respectively. Thus, we get:}
\begin{itemize}
\item[i)] \emph{For $L_p(\RR) \to L_p(\RR; R_p)$, we need, in addition to Condition A, $$M_{jg}(\xi) \, = \, \langle A_{g}, B_{\xi j} \rangle_{\mathcal{K}_1} \, = \, \langle A_{g \xi}', B_j' \rangle_{\mathcal{K}_2}$$ where $\sup_{g,j,\xi} \|A_g\|_{\mathcal{K}_1}  \|B_{\xi j}\|_{\mathcal{K}_1} + \sup_{g,j,\xi} \|A_{g \xi}'\|_{\mathcal{K}_2}  \|B_j'\|_{\mathcal{K}_2}$ is finite.}

\vskip5pt 

\item[ii)] \emph{For $L_p(\RR) \to L_p(\RR; C_p)$, we need, in addition to Condition A $$M_{jg}(\xi) \, = \, \langle A_{gj}, B_{\xi} \rangle_{\mathcal{K}_1}
  $$ where $\sup_{g,j,\xi} \|A_{gj}\|_{\mathcal{K}_1}  \|B_{\xi j}\|_{\mathcal{K}_1} 
  $ is again finite.}
\end{itemize}
\end{remark}

\subsection{The twisted column estimate} \label{SectColEst}

Twisted Riesz transforms $\widetilde{R}_{\psi_\varepsilon, u}$ and $\widetilde{R}_{\psi_\varepsilon,\mathbf{u}}$ were introduced in \eqref{Eq-TwistRiesz} and \eqref{eq-SqTwistRiesz}. Given $p>2$, in this section we shall establish the column inequality $L_p(\RR;C_p) \to L_p(\RR;C_p)$ for $\widetilde{R}_{\psi_\varepsilon, \mathbf{u}}$ and certain directions $\mathbf{u} = (u_j) = (\Psi_\varepsilon^{-1}(\varphi_j M))$ in $\H_\varepsilon$ which are determined by a symbol $M$ satisfying the H\"ormander-Mikhlin conditions. The definition of 
$L_p(\RR;C_p)$ is given in \eqref{eq-Cp} and corresponds to the norm in 
$L_p(\RR\overline \otimes \mathcal B(\ell_2))= L_p(L_\infty(\R^d)\overline \otimes \mathcal B(L_2(\Sigma)\otimes_2\ell_2))$. Recall that the symbol of $R_{\psi_\varepsilon,u_j}^g$ equals 
\begin{equation} \label{Eq-ujugj}
\rho_{\varepsilon,u_j}(\alpha_g(\xi)) = \frac{\langle b_\varepsilon(\alpha_g(\xi)), u_j \rangle_{\H_\varepsilon}}{\sqrt{\psi_\varepsilon(\xi)}} \frac{|\xi|^\varepsilon}{|\alpha_g(\xi)|^\varepsilon} = \frac{\langle b_\varepsilon(\xi), u_{gj} \rangle_{\H_\varepsilon}}{\sqrt{\psi_\varepsilon(\xi)}} \frac{|\xi|^\varepsilon}{|\alpha_g(\xi)|^\varepsilon}.
\end{equation}
Indeed, the last identity follows from 
\begin{eqnarray*}
\big\langle b_\varepsilon(\alpha_g(\xi)), u_j \big\rangle_{\H_\varepsilon} & = & \int_{\R^d} \big( \exp_{\alpha_g(\xi)}(x) - 1\big) u_j(x) \, d\mu_\varepsilon(x) \\
& = & \int_{\R^d} \big( \exp_{\xi}(s) - 1\big) \underbrace{u_j(\alpha_g'(s)) \, \frac{|s|^{d+2\varepsilon}}{|\alpha_g'(s)|^{d+2\varepsilon}}}_{u_{gj}(s)} \, d\mu_\varepsilon(s)
\end{eqnarray*}
with $\alpha_g' = (\alpha_{g}^*)^{-1}$. In particular, $\widetilde{R}_{\psi_\varepsilon,u_j} = R_{\psi_\varepsilon,\widetilde{u}_j} \circ \widetilde{H}$ with:
\begin{itemize}
\item A \lq\lq homogeneous twisted" multiplier
\begin{eqnarray*}
\widetilde{H}(f)(x) \!\! & = & \!\! \Big( \int_{\R^d} \frac{|\xi|^\varepsilon}{|\alpha_g(\xi)|^\varepsilon} \widehat{f}_{gh}(\xi) \exp_\xi(x) \, d\xi \Big)_{\hskip-2pt g,h \in \Sigma}.
\end{eqnarray*}

\item A \lq\lq multidirectional" Riesz transform
\begin{eqnarray*}
\hskip6pt R_{\psi_\varepsilon, \widetilde{u}_j}(f)(x) \!\! & = & \!\! \Big( R_{\psi_\varepsilon,u_{gj}}(f_{gh}) \Big)_{\hskip-2pt g,h \in \Sigma} \\ \!\! & = & \!\!  \Big( \int_{\R^d} \frac{\langle b_\varepsilon(\xi), u_{gj} \rangle_{\H_\varepsilon}}{\sqrt{\psi_\varepsilon(\xi)}} \widehat{f}_{gh}(\xi) \exp_\xi(x) \, d\xi \Big)_{\hskip-2pt g,h \in \Sigma}.
\end{eqnarray*}
\end{itemize}
Therefore, it suffices to prove that 
\begin{eqnarray*}
\| \widetilde{H}(f) \|_{L_p(\RR)} & \le_{\mathrm{cb}} & C_{p,d}(\Sigma) \, \| f \|_{L_p(\RR)}, \\ [3pt]
\Big\| \sum_{j \in \Z} R_{\psi_\varepsilon,\widetilde{u}_j} (f_j) \otimes e_{j1} \Big\|_{L_p(\RR; C_p)} & \le_{\mathrm{cb}} & C_{p,d}(\Sigma) \, \Big\| \sum_{j \in \Z} f_j \otimes e_{j1} \Big\|_{L_p(\RR; C_p)}.
\end{eqnarray*}

\subsubsection{Multidirectional Riesz transforms} \label{SubsectRieszCol}

\begin{proposition} \label{PropColEst}
Let $p \ge 2$ and $\Sigma$ relatively compact, then \\ [-3pt]
$$\Big\| \sum_{j \in \Z} R_{\psi_\varepsilon,\widetilde{u}_j} (f_j) \otimes e_{j1} \Big\|_{L_p(\RR; C_p)} \, \le_{\mathrm{cb}} \, C_{p,d}(\Sigma) \Big\| \sum_{j \in \Z} f_j \otimes e_{j1} \Big\|_{L_p(\RR; C_p)}$$ provided $\displaystyle \frac{\langle b_\varepsilon(\xi), u_j \rangle}{\sqrt{\psi_\varepsilon(\xi)}} = \varphi_jM$ for some $M$ satisfying that $\displaystyle \sup_{j \in \Z} \big\| \varphi_0^2 M(2^j \cdot) \big\|_{\mathsf{H}_{\frac{d}{2} + \varepsilon}} < \infty$.
\end{proposition}

\dem Given an orthonormal basis $e_1, e_2, \ldots$ of $\H_\varepsilon$, $u_{gj} = \sum_k \langle u_{gj}, e_k \rangle e_k$. This gives $R_{\psi_\varepsilon, u_{gj}} = \sum_k \langle u_{gj}, e_k \rangle R_{\psi_\varepsilon, e_k}$. By \cite[Theorem A1 and Remark 1.8]{JMP2} we obtain
$$\Big\| \sum_{j \in \Z} \sum_{k \ge 1} R_{\psi_\varepsilon,e_k} (f_j) \otimes e_{j1} \otimes e_{k1} \Big\|_{L_p(\RR; C_p(\N \times \Z))} \, \le_{\mathrm{cb}} \, C_p \Big\| \sum_{j \in \Z} f_j \otimes e_{j1} \Big\|_{L_p(\RR; C_p(\Z))}.$$ It therefore suffices to prove that
\begin{equation} \label{Eq-Lambdau}
\Lambda_u: \sum_{j \in \Z} \sum_{k \ge 1} \Big( a_{gh}^{jk} \Big) \otimes e_{j1} \otimes e_{k1} \mapsto \sum_{j \in \Z} \Big( \sum_{k \ge 1} \langle u_{gj}, e_k \rangle a_{gh}^{jk} \Big) \otimes e_{j1}
\end{equation}
is a cb-bounded map $\Lambda_u \hskip-2pt : S_p(L_2(\Sigma); C_p(\N \times \Z)) \to S_p(L_2(\Sigma); C_p(\Z))$ for $p \ge 2$. Factorization of $S_p = C_p \otimes_\mathrm{h} R_p$ in terms of the Haagerup tensor product \cite{P1} yields 
\begin{eqnarray*}
S_p(L_2(\Sigma); C_p(\Z)) & = & C_p(L_2(\Sigma) \otimes_2 \ell_2(\Z)) \otimes_\mathrm{h} R_p(L_2(\Sigma)), \\
S_p(L_2(\Sigma); C_p(\N \times \Z)) & = & C_p(L_2(\Sigma) \otimes_2 \ell_2(\N \times \Z)) \otimes_\mathrm{h} R_p(L_2(\Sigma)). 
\end{eqnarray*}
The map $\Lambda_u$ acts trivially on $R_p(L_2(\Sigma))$ |the $h$-variable| so that cb-boundedness is equivalent to that of the map $\Lambda_u: C_p(L_2(\Sigma) \otimes_2 \ell_2(\N \times \Z)) \to C_p(L_2(\Sigma) \otimes_2 \ell_2(\Z))$. On the other hand, $C_p$-spaces are homogeneous Hilbertian operator spaces \cite{P2}, that means that a bounded map on $C_p$ is automatically completely bounded with the same norm. Moreover, since one can tensorize completely bounded maps with the
Haagerup tensor product, we are reduced to simply proving boundedness. In formulas, for every linear map $V: C_p \to C_p$,
$$\| V \otimes Id_{R_p}: C_p\otimes_h R_p \to C_p\otimes_h R_p \|_{\rm cb}=\|V : C_p\to C_p\|_{\rm cb}=\|V: C_p\to C_p\|.$$
Now the Hilbertian nature of the spaces make then isomorphic for all values of $p$ and it suffices to prove the boundedness of $\Lambda_u$ in $S_2[C_2]$. We have
\begin{eqnarray*}
\|\Lambda_u(a)\|_{S_2[C_2]}^2 & = & \int_{\Sigma \times \Sigma}  \sum_{j \in \Z} \Big| \sum_{k \ge 1} \langle u_{gj}, e_k \rangle a_{gh}^{jk} \Big|^2 d\mu(g)d\mu(h) 
\\ & \le & \sup_{\begin{subarray}{c} j \in \Z \\ g \in \Sigma \end{subarray}} \hskip2pt \sum_{k \ge 1} \big| \langle u_{gj}, e_k \rangle \big|^2 \, \|a\|_{S_2[C_2]}^2
\ = \ \sup_{\begin{subarray}{c} j \in \Z \\ g \in \Sigma \end{subarray}} \hskip2pt \|u_{gj}\|_{\H_\varepsilon}^2 \, \|a\|_{S_2[C_2]}^2.
\end{eqnarray*}
By compactness of $\Sigma \subset \G$ 
\begin{eqnarray} \label{Eq-Normugjuj}
\hskip20pt \sup_{g \in \Sigma} \|u_{gj}\|_{\H_\varepsilon}^2 & = & \sup_{g \in \Sigma} \int_{\R^d} |u_j(\alpha_g'(s))|^2 \frac{|s|^{2d + 4\varepsilon}}{|\alpha_g'(s)|^{2d + 4\varepsilon}} \, \frac{ds}{|s|^{d + 2 \varepsilon}} \\ \nonumber
& = & \sup_{g \in \Sigma} \int_{\R^d} |u_j(s)|^2 \frac{|\alpha_g^*(s)|^{d + 2\varepsilon}}{|s|^{d + 2\varepsilon}} \, d\mu_\varepsilon(s) \ \le \ C_d(\Sigma) \|u_j\|_{\H_\varepsilon}^2.
\end{eqnarray}
The uniform bound for $u_j = \Psi_\varepsilon^{-1}(\varphi_j M)$ follows from Lemmas \ref{SobVsClassical} and \ref{HMLemma}. \fin 

\begin{remark} \label{RemFailureRow}
\emph{If $p > 2$, it turns out that 
$$\big\| \widetilde{T}_{\dot m} \hskip-2pt : L_p(\RR) \to L_p(\RR) \big\|_{\mathrm{cb}} \, \le \, C_{p,d}(\Sigma) \big\| \widetilde{R}_{\psi_\varepsilon, \mathbf{u}} \hskip-2pt : L_p(\RR;RC_p) \to L_p(\RR;RC_p) \big\|_{\mathrm{cb}}$$
for some Riesz directions $\mathbf{u} = (u_j)$ dictated by $\dot m$. It follows from Propositions \ref{PropColEst} and \ref{PropHomogeneous} below that the column estimate holds for the term in the right-hand side. However, the first factorization in Remark \ref{RemRCCZ} i) fails for twisted Riesz transforms, and everything indicates that the row estimate does not hold in general. This is partly explained from the asymmetric nature of our twist. This difficulty will be sorted out in the next section using a new Littlewood-Paley type inequality for Schur multipliers, together with a local inversion trick which allows us to write the row estimates in terms of column ones.}
\end{remark}

\subsubsection{The homogeneous twisted multiplier} \label{SubsectHTwisted}

We now study the map $\widetilde{H}$ using the above Calder\'on-Zygmund methods. In addition, we need a standard Sobolev bound for Schur multipliers which we now recall. Given a pair of cubes $\mathrm{Q}_j \subset \R^{d_j}$, let $S: \mathrm{Q}_1 \times \mathrm{Q}_2 \to \C$. The following result establishes a sufficient condition on $S$ to be a Schur multiplier in $\mathcal{B}(S_\infty(L_2(\mathrm{Q}_1), L_2(\mathrm{Q}_2)))$. 

\begin{lemma} \label{LemmaSobolev}
If $\gamma = (1,1,\ldots, 1) \in \R^{d_1+ d_2}$ 
\begin{eqnarray*}
\lefteqn{\hskip-30pt \Big\| \Big( S(x,y) A(x,y) \Big) \Big\|_{S_\infty(L_2(\mathrm{Q}_1),L_2(\mathrm{Q}_2))}} \\
\hskip20pt & \le_{\mathrm{cb}} & \sum_{\rho \le \gamma} \big\| \partial^\rho S \big\|_{L_2(\mathrm{Q}_1 \times \mathrm{Q}_2)} \Big\| \Big( A(x,y) \Big) \Big\|_{S_\infty(L_2(\mathrm{Q}_1),L_2(\mathrm{Q}_2))}.
\end{eqnarray*}
\end{lemma}

\dem If $\ell_j = \mbox{length}(\mathrm{Q}_j)$ and $\mathbf{Z}_j = \ell_j^{-1} \Z^{d_j}$
$$S(x,y) = \sum_{(p,q) \in \mathbf{Z}_1 \times \mathbf{Z}_2} \widehat{S}(p,q) \underbrace{\exp(2\pi i \langle x,p \rangle)}_{u_x(p)} \hskip1pt \underbrace{\exp(2\pi i \langle y, q \rangle)}_{u_y(q)}.$$ By Grothendieck's characterization \cite[Proposition 1.1]{PisAJM} (or \cite[Theorem 1.7]{LdlS} for a formulation with continuous parameters), it suffices to factorize the symbol $S(x,y) = \langle A_x,B_y \rangle_\mathcal{K}$ for some Hilbert space $\mathcal{K}$ with uniformly bounded vectors $A_x, B_y \in \mathcal{K}$. Using the Fourier expansion above, this will be the case when the Fourier coefficients are summable, since we may pick 
\begin{eqnarray*}
A_x & = & \sum_{(p,q) \in \mathbf{Z}_1 \times \mathbf{Z}_2} |\widehat{S}(p,q)|^{\frac12} u_x(p) \otimes \delta_{p,q}, \\
B_y & = & \sum_{(p,q) \in \mathbf{Z}_1 \times \mathbf{Z}_2} |\widehat{S}(p,q)|^{\frac12} \mbox{sgn}(\widehat{S}(p,q)) u_y(q) \otimes \delta_{p,q}.
\end{eqnarray*}  
Elementary integration by parts and Plancherel theorem give $$\sum_{(p,q) \in \mathbf{Z}_1 \times \mathbf{Z}_2} |\widehat{S}(p,q)| \, \lesssim \, \sum_{\rho \le \gamma} \|\partial^\rho S\|_{L_2(\mathrm{Q}_1 \times \mathrm{Q}_2)}.$$ Multi-indices of order $j$ are used for $(p,q)$ with $d_1 + d_2 - j$ vanishing entries. \fin

\begin{proposition} \label{PropHomogeneous}
If $p \ge 2$ 
$$\| \widetilde{H}(f) \|_{L_p(\RR)} \, \le_{\mathrm{cb}} \, C_{p,d}(\Sigma) \| f \|_{L_p(\RR)}.$$
\end{proposition}

\dem According to Proposition \ref{CZProposition}, it suffices to show that the twisted symbol is uniformly Sobolev-smooth and admits Schur factorization. The stronger Mikhlin smoothness condition $\sup_g | \partial_\xi^\gamma M_g(\xi) | \le C_\Sigma |\xi|^{-|\gamma|}$ for $M_g(\xi) = |\xi|^\varepsilon |\alpha_g(\xi)|^{-\varepsilon}$ holds for all multi-index $\gamma$. This is a simple exercise which follows from the compactness of $\Sigma$. It remains to prove Schur's factorization. Equivalently $$\frac{|\xi|^\varepsilon}{|\alpha_g(\xi)|^\varepsilon} = \frac{1}{\big|\alpha_g(\frac{\xi}{|\xi|}) \big|^\varepsilon} \in \mathcal{B}\big(S_\infty(L_2(\Sigma), L_2(\R^d))\big)$$ as a Schur multiplier in $(g,\xi) \in \Sigma \times \R^d$. Of course, Schur factorization is stable by composition, so we may replace $(g,\xi)$ by $(\alpha_g, \xi) \in \Pi \times \R^d$ for some compact set $\Pi \subset \SLk$. Moreover, since boundedness and cb-boundedness are equivalent for this class of Schur multipliers, we may use the homogeneity of our symbol to reduce it to $(\alpha_g,\xi) \in \Pi \times \mathbf{S}^{d-1}$. The set $\Pi \times \mathbf{S}^{d-1}$ is a compact manifold in $\R^{d^2} \times \R^d$ and our symbol $H(g,\xi) = H_g(\xi)$ admits a smooth extension |still denoted by $H$| to an open neighborhood of it. In particular, since row/column restriction is a continuous operation for Schur multipliers, we may cover that open (relatively compact) set by a finite number $Q_{1j} \times Q_{2j} \subset \R^{d^2} \times \R^d$ of pairs of cubes. This gives 
$$\|H\|_{\mathcal{B}(S_\infty(L_2(\Sigma), L_2(\mathbf{S}^{d-1})))} \le \sum_{j=1}^\mathrm{N} \|H\|_{\mathcal{B}(S_\infty(L_2(Q_{1j}), L_2(Q_{2j})))}.$$ The assertion follows from Lemma \ref{LemmaSobolev} and the smoothness of our symbol. \fin 

\section{\bf Proof of Theorem A} \label{SectSmoothness}

As explained in the Introduction, the main challenge in the proof of Theorem A comes from the local behavior of the multiplier around the singularity, whereas the asymptotic behavior will follow at the end from an elementary patching argument due to the exponential nature of the metric. Thus, we assume momentarily that $m \hskip-2pt : \G \to \C$ is supported in a compact neighborhood around the identity.
 
\subsection{Local inversion} \label{SubSectLocInversion}

Set $$\mathcal{I}(A) = (A + e)^{-1} - e \quad \mbox{for} \quad A \in \GL - e.$$
Let $\mathrm{K}$ be a compact set in $\GL - e$ containing $0$. In the following we shall use local stability properties of $\mathcal{I}$ relative to $\mathrm{K}$, which we now collect. Let us fix $d=n^2$ for the rest of this section and let $(\varphi_j)_{j \in \Z}$ be the Littlewood-Paley partition of unity \eqref{Eq-LP} in $\R^d$. Construct the partition of unity 
\begin{equation} \label{eq-PUsigmaj}
\sigma_j = \frac{1}{2\mathrm{N_K}+1} \sum_{k=j-\mathrm{N_K}}^{j+\mathrm{N_K}} \varphi_k^2 \quad \mbox{for some} \quad \mathrm{N_K} \in \N.
\end{equation}
Let us recall that $(2\mathrm{N_K}+1) \sigma_j \equiv 1$ in the set $\big\{ \xi \in \R^d: 2^{j-\mathrm{N_K}} \le |\xi| \le 2^{j+\mathrm{N_K}} \big\}$.
   
\begin{lemma} \label{LemmaSobolevInversion}
If $\mathrm{supp} \hskip1pt m \subset \mathcal{I}(\mathrm{K})$ and $\mu_j = (\sigma_j m) \circ \mathcal{I}$ $$\sup_{j \in \Z} \| \mu_j\|_{\mathrm{W}_{d,\varepsilon}^2(\R^d)} \, \le \, C_{\mathrm{K},d} \sup_{j \in \Z} \big\| \varphi_0^2 m(2^j \cdot) \big\|_{\mathsf{H}_{\frac{d}{2} + \varepsilon}^2(\R^d)}.$$
\end{lemma}

\dem If $A \in \mathrm{K}$ and $|\cdot|$ denotes the Hilbert-Schmidt norm, we first observe that $|\mathcal{I}(A)| \approx |A|$ up to a constant $C_\mathrm{K}$ depending only of $\mathrm{K}$. Indeed, matrix inversion is operator Lipschitz on any compact $\mathrm{K} + e$ in $\mathrm{GL}_n(\R)$, so we get the upper estimate over $\mathrm{K}$ with some constant $c_\mathrm{K}$. The lower estimate follows since $A = \mathcal{I}((A+e)^{-1}-e)$ and $\mathcal{I}(\mathrm{K})$ is another compact set in $\mathrm{GL}_n(\R) - e$. 

We shall use as well a standard fact on Sobolev norms of composition by smooth functions \cite{T}. Assume that $\mathrm{supp} \hskip1pt f \subset \mathrm{K} \subset \Delta \subset \overline{\Delta } \subset \Lambda$ for some open domains $\Delta$ and $\Lambda$. Given $\Phi \in \mathcal{C}^\infty(\Lambda)$ and $\Psi: \Lambda \to \Lambda$ a diffeomorphism satisfying $\Psi^2 = id$, let us set $$\mathrm{a} = \max_{|\gamma| \le [\alpha]+1} \|\partial^\gamma \hskip-1pt \Psi\|_{L_\infty(\Psi(\Delta))}.$$
Then, the following inequality holds 
\begin{equation}  \label{Eq-Sob2}
\|f \circ \Psi\|_{\mathsf{H}_\alpha^2(\R^d)} \, \le \, C_{\mathrm{K}, \mathrm{a}, \alpha} \|f\|_{\mathsf{H}_\alpha^2(\R^d)}.
\end{equation}
Fix a relatively compact open set $\Delta$ in $\GL -e$ containing $\mathrm K$ as above, on which the derivatives up to some order of $\mathcal I$ over $\mathcal{I}(\Delta)$ can be controlled by a constant determined by $\mathrm K$. Applying the triangular inequality, the dilation invariance in $\mathrm{W}_{d,\varepsilon}^2(\R^d)$ and Lemma \ref{SobVsClassical} (for a fixed $j$), we get
$$\| \mu_j\|_{\mathrm{W}_{d,\varepsilon}^2(\R^d)} \, 
\le \,  \sup_{\ell \in \Z}  \big\| \big[(\varphi_\ell^2 m) \circ \mathcal{I} \big](2^\ell \cdot)
\big\|_{\mathsf{H}_{\frac{d}{2}+\varepsilon}^2(\R^d)}.$$ 
As in the proof of Lemma \ref{SobVsClassical}, this requires to check that $[(\varphi_\ell^2 m) \circ \mathcal{I}] (2^\ell \cdot)$ is supported in certain corona around the unit sphere of $\R^d$ which is independent of the value of $\ell \in \Z$. Indeed, we know that $\mathrm{supp} \hskip1pt m \circ \mathcal{I}(2^\ell \cdot) \subset 2^{-\ell} \mathrm{K}$. Given $\xi \in 2^{-\ell} \mathrm{K}$ and according to our first observation, we get $|\mathcal{I}(2^\ell \xi)| \approx |2^\ell \xi|$ up to a constant determined by $\mathrm{K}$ but independent of $\ell \in \Z$. This proves the required condition. Let us now define $\mathcal{I}_\ell(\xi) = 2^{-\ell} \mathcal I(2^\ell\xi)$, so that
$[(\varphi_\ell^2 m) \circ \mathcal{I}](2^\ell \cdot)=(\varphi_0^2 m(2^\ell\cdot) )\circ \mathcal{I}_\ell = f_\ell \circ \Psi_\ell$
 where $f_\ell =\varphi_0^2  m(2^\ell \cdot)$ is supported by
$2^{-\ell} \mathcal{I}(\mathrm{K})\cap \{\frac 1 2\leq |\xi|\leq 2\}$ and $\Psi_\ell = \mathcal{I}_\ell$ is defined in $\Lambda_\ell = 2^{-\ell} (\GL - e)$.
It is clear that $f_\ell = 0$ when $2^\ell > 2 \mathrm{diam} (\mathcal{I}(\mathrm{K}))$. Using our observation one more time, this implies that $2^\ell \le 2 C_\mathrm{K} \mathrm{diam}(\mathrm{K})$. Moreover, given $\Delta_\ell = 2^{-\ell} \Delta$ around $2^{-\ell} \mathrm{K}$, we get by construction
$$\max_{0 \le |\gamma| \le [\frac{d}{2}+\varepsilon]+1} \| \partial^\gamma \hskip-1pt \Psi_\ell \|_{L_\infty(\Psi_\ell(\Delta_\ell))} \le 2^{\ell |\gamma|-1} C_K'
\le  C_K'' \mathrm{diam} (\mathrm{K})^{|\gamma|}.$$ Therefore, the assertion in the statement follows from inequality \eqref{Eq-Sob2} above. \fin

\subsection{The local theorem}

We have now all the ingredients to prove the local form of Theorem A. Recall that we work in dimension $d = n^2$ and the natural cocycle is given by $\beta(g) = g-e$ with cocycle action $\alpha_g(\xi) = g \xi$. We begin with an $S_p$-column inequality for Schur multipliers in $\SL$ which satisfy a uniform Mikhlin bound. 

\begin{proposition}\label{Prop:SpCp_HM} Let $\mathfrak S=(m_j)_{j \in \Z}$ be a sequence of functions $ \SL \to \C$ and let $\Sigma \subset \SL$ be a relatively compact subset. Assume that $m_j(g) = \dot m_j(g-e)$ for certain $\dot m_j \colon \R^d \to \C$ satisfying 
\[ C_\mathrm{hm} (\mathfrak S) := \sup_{j \in \Z} \| \dot m_j\|_{\mathrm{W}_{d,\varepsilon}^2(\R^d)} <\infty \quad \mbox{for some} \quad 0<\varepsilon<1.\]
Then, for any $1<p<\infty$ and any sequence $A_j \in S_p(L_2 (\Sigma))$, we have
\[ \Big\| \sum_{j \in \Z} S_{m_j}(A_j) \otimes \delta_j\Big\|_{S_p[C_p]} \leq C_{p,d}(\Sigma) C_{\mathrm{hm}} (\mathfrak S) \Big\| \sum_{j \in \Z} A_j \otimes \delta_j\Big\|_{S_p[C_p]}.\]
\end{proposition}

\dem By a straightforward adaptation of the proof of Proposition \ref{PropTransference}, the norm on $S_p[C_p]$ of the map in the statement is bounded above by the norm on $L_p(\RR;C_p)$ of the map $$\widetilde{T}_{\mathfrak S}: \sum_{j \in \Z} a_j \otimes \delta_j \mapsto \sum_{j \in \Z}  \widetilde{T}_{\dot{m_j}}(a_j) \otimes \delta_j.$$ Moreover, Lemma \ref{HMLemma} yields \[\widetilde{T}_{\dot m_j} = \widetilde{R}_{\psi_{\varepsilon}, u_j}\] for a sequence $(u_j)_j$ in $\H_\varepsilon = L_2(\R^d,\mu_\varepsilon)$ which is uniformly bounded in norm by the constant $C_{\mathrm{hm}}((m_j)_j)$. Next, we follow the exact same argument as in Section \ref{SectColEst} and factor \[ \widetilde R_{\psi_{\varepsilon,u_j}} = R_{\psi_{\varepsilon}, \widetilde{u}_j} \circ \widetilde{H}.\]
Therefore, $\widetilde{T}_{\mathfrak S}$ is bounded above by the product of the completely bounded norm of $\widetilde H$ with the norm on $L_p(\RR;C_p)$ of $\sum_j a_j \otimes \delta_j \mapsto \sum_j R_{\psi_{\varepsilon}, \widetilde{u}_j}(a_j) \otimes \delta_j$. The assertion then follows by applying Proposition \ref{PropColEst} and Proposition \ref{PropHomogeneous}. \fin

\begin{theorem} \label{Prop:HM_Schur}
Let $\Sigma$ be a relatively compact subset of $\G=\SL$ and $m \colon \G \to \C$. Assume $m(g) = \dot{m}(g-e)$ for some $\dot{m} \colon \R^{d} \to \C$ satisfying the Mikhlin type condition
\[\big| \partial_\xi^\gamma \dot{m}(\xi) \big| \, \le \,  |\xi|^{-|\gamma|} \quad \mbox{for all} \quad 0 \le |\gamma| \le \Big[ \frac{n^2}{2}\Big] + 1.\]
Then, the following inequality holds for each $1 < p < \infty$
\[\big\| S_m \colon S_p(L_2(\Sigma)) \to S_p(L_2(\Sigma)) \big\|_{\mathrm{cb}} \, \le \, C_p(\Sigma).\]
\end{theorem}

\dem
We shall prove that $\|S_m(A)\|_p \leq C_p(\Sigma) \|A\|_p$ for $A \in S_p(L_2(\Sigma))$. The completely bounded norm is also dominated by the same constant with the same argument. Alternatively, we can use \cite[Theorem 1.19]{LdlS}, which implies that, if $\Sigma$ is open (which we can always assume by enlarging it), the norm and completely bounded norm of the Schur multiplier $S_m$ coincide for $m$ continuous. Next, by duality we may restrict to the case $p\geq 2$. Moreover, by replacing $\dot m$ by $\psi \dot m$ for a suitable compactly supported smooth function $\psi$ equal to $1$ on $\Sigma\Sigma^{-1}$, the Schur multiplier $S_m$ is not affected while the Mikhlin condition for $\dot{m}$ holds up to a constant $C_\Sigma$ determined by $\Sigma$. Therefore, we may and will assume that there is a compact subset $\mathrm{K}$ in $\GL-e$ such that $\dot m=0$ outside of $\mathrm{K}$.

\noindent \textbf{A. Reduction to an $RC_p$ inequality.} By Remark \ref{Rem-ShiftofTwist}, the classical Littlewood theorem can be transfered to twisted Schur multipliers. More precisely, consider $\varphi_j$ as in \eqref{Eq-LP} and the Schur multiplier $S_j$ on $L_2(\G)$ with symbol $$(g,h) \mapsto \varphi_j^2 \big( \alpha_{g^{-1}}(\beta(gh^{-1})) \big) = \varphi_j^2 \big( \beta(h^{-1}) - \beta(g^{-1}) \big) \quad  \mbox{for} \quad \beta(g) = g-e.$$ This gives an unconditional decomposition of the identity 
\begin{equation} \label{Eq-LPTwist1}
\Big\| \sum_{j \in \Z} S_j(A) \otimes \delta_j \Big\|_{S_p[RC_p]} \, \approx\|A\|_{S_p(L_2(\G))},
\end{equation}
see \eqref{eq-RCp} for the definition of $S_p[RC_p]$. Next, we construct a partition of unity of the form \eqref{eq-PUsigmaj}. More precisely, pick a positive integer $\mathrm{N}_\Sigma$ (to be fixed) and consider 
$$\sigma_j = \frac{1}{2\mathrm{N}_\Sigma+1} \sum_{k=j-\mathrm{N}_\Sigma}^{j+\mathrm{N}_\Sigma} \varphi_k^2.$$
It is clear that $\sum_j \sigma_j = 1$ and
\begin{itemize}
\item $\mbox{supp} \hskip1pt \varphi_j \subset \big\{ \xi \in \R^d: 2^{j-1} \le |\xi| \le 2^{j+1} \big\}$,

\vskip5pt

\item $(2\mathrm{N}_\Sigma+1) \sigma_j \equiv 1$ in $\big\{ \xi \in \R^d: 2^{j-\mathrm{N}_\Sigma} \le |\xi| \le 2^{j+\mathrm{N}_\Sigma} \big\}$.
\end{itemize}
Thus, since $\Sigma$ is relatively compact and the action $\alpha_g(h) = gh$ is continuous, there must exist $\mathrm{N}_\Sigma \in \Z_+$ determined by $\Sigma$ such that $(2\mathrm{N}_\Sigma + 1) \sigma_j(\alpha_g(\xi)) \equiv 1$ for all $(j,g,\xi) \in \Z \times \Sigma \times \mbox{supp} \hskip1pt \varphi_j$. This yields for $(g,h)\in \Sigma \times \Sigma$
\[ \varphi_j^2 \big( \alpha_{g^{-1}}(\beta(gh^{-1})) \big)m(gh^{-1}) = (2 \mathrm{N}_\Sigma + 1) \varphi_j^2 \big( \alpha_{g^{-1}}(\beta(gh^{-1})) \big) \sigma_j(\beta(gh^{-1})) m(gh^{-1}).\]
This means that \[S_j(S_m(A)) = (2 \mathrm{N}_\Sigma +1)S_{m_j}(S_j(A))\] for $m_j = (\sigma_j \circ \beta) m$. Therefore, by \eqref{Eq-LPTwist1} the assertion is equivalent to 
\begin{equation}\label{eq:Smj_on_RCp} \Big\| \sum_{j \in \Z} S_{m_j}(A_j) \otimes \delta_j\Big\|_{S_p[RC_p]} \leq C_p(\Sigma) \Big\| \sum_{j \in \Z} A_j \otimes \delta_j\Big\|_{S_p[RC_p]}\end{equation}
for $A_j = S_j(A)$. We shall prove \eqref{eq:Smj_on_RCp} for every family $(A_j)_{j \in \Z} \in S_p[RC_p]$.

\noindent \textbf{B. Rows to columns by local inversion.} We have 
\[S_{m_j}(A^*)^* = S_{\mu_j}(A) \quad \mbox{for} \quad \mu_j(g) = \overline{m_j(g^{-1})}.\] 
This allows us then to write the row term in \eqref{eq:Smj_on_RCp} as follows
\[\Big\| \sum_{j \in \Z} S_{m_j} (A_j) \otimes \delta_j \Big\|_{S_p[R_p]}  = \Big\| \sum_{j \in \Z} S_{\mu_j} (A_j^*) \otimes \delta_j \Big\|_{S_p[C_p]}.\]
Therefore, the $RC_p$-inequality \eqref{eq:Smj_on_RCp} follows from the following inequalities 
\begin{equation}\label{eq:Smj_on_Cp} 
\Big\| \sum_{j \in \Z} S_{\bullet_j}(A_j) \otimes \delta_j\Big\|_{S_p[C_p]} \leq C_p(\Sigma) \Big\| \sum_{j \in \Z} A_j \otimes \delta_j\Big\|_{S_p[C_p]}
\end{equation}
for $\bullet = m$ and $\bullet=\mu$. The case $\bullet = m$ is the content of Proposition \ref{Prop:SpCp_HM}, since $\dot{m}_j = \sigma_j \dot{m}$ is uniformly bounded in $\mathrm{W}_{d,\varepsilon}^2(\R^d)$ by virtue of the triangular inequality and Lemma \ref{SobVsClassical}. On the other hand, 
group inversion $g \mapsto g^{-1}$ is smooth in any relatively compact set of $\G$. Thus $\mu_j$ satisfies (up to constants) the same Mikhlin conditions as $m_j$. To be more precise, using the map $\mathcal{I}$ from Section \ref{SubSectLocInversion}  
\[\mu_j(g) = \dot{\mu}_j(g-e) \textrm{ with } \dot{\mu}_j(\xi) = \overline{\dot{m}_j \big( (\xi + e)^{-1} - e)} = \overline{(\sigma_j \dot{m}) \circ \mathcal{I}(\xi)}.\]
By Lemmas \ref{SobVsClassical} and \ref{LemmaSobolevInversion} (the norm in $\mathrm{W}_{d,\varepsilon}^2(\R^d)$ is conjugate invariant) we get
\begin{equation} \label{Eq-InversionControl}
\sup_{j \in \Z} \|\dot{\mu}_j\|_{\mathrm{W}_{d,\varepsilon}^2(\R^d)} \, \le \, C_\Sigma \max_{|\gamma| \le [\frac{d}{2}] + 1} \big\| |\xi|^{|\gamma|} \partial_\xi^\gamma \dot{m}(\xi) \big\|_\infty < \infty.
\end{equation}
So \eqref{eq:Smj_on_Cp} for $\bullet = \mu$ also follows from Proposition \ref{Prop:SpCp_HM}. This proves the assertion. \fin

\begin{TheoLocal}
Assume that $m \colon \SL \to \C$ is supported by a relatively compact neighborhood $\Omega$ of the identity and satisfies \eqref{Eq-ThmA}. Then, $T_m$ is completely $L_p$-bounded for $1  <  p <  \infty$ by $C_p(\Omega)$, with $C_p(\Omega) \approx C_p C_\mathrm{hm}$ for $\Omega$ small.
\end{TheoLocal}

\dem
By duality and interpolation we may assume that $p \in 2\Z_+$. Fix a smooth function $\varphi \colon \R \to \R_+$ which is $0$ outside of $[1/2,2]$ and equal to $1$ at $1$. Extend the multiplier $m \colon \SL \to \C$ to a function $M\colon \R^d = M_n(\R) \to \C$ by the formula $M(A) = \varphi(\det A) m(A/(\det A)^{\frac 1 n})$ if $\det(A)>0$ and $M(A)=0$ otherwise. Let $\dot m(\xi) = M(\xi + e)$, so that $m(g) = \dot m(g-e)$. It follows from \eqref{Eq-ThmA} that 
\begin{equation} \label{eq-EuclideanHM}
\big| \partial_\xi^\gamma \dot{m}(\xi) \big| \, \le \, C(\Omega) C_{\rm hm} |\xi|^{-|\gamma|} \quad \mbox{for all} \quad 0 \le |\gamma| \le \Big[ \frac{n^2}{2}\Big] + 1.
\end{equation}
Choosing $\Sigma$ large enough, the local theorem follows from Theorems \ref{Thm-LocalTransf} and \ref{Prop:HM_Schur}.
\fin
 
\begin{remark}
\emph{The proof above gives in fact a stronger result. Namely, by the left invariant nature of our Lie differential operators, it suffices to prove \eqref{eq-EuclideanHM} in a small neighborhood of the identity. Then, the Mikhlin constants of $\dot m$ are dominated by a subfamily of Lie derivatives of $m$. More precisely, it suffices to assume \eqref{Eq-ThmA} in Theorem A for a family $\Gamma_0$ of multi-indices satisfying that every other multi-index $\gamma$ with $|\gamma| \le [n^2/2]+1$ is the permutation of an element in $\Gamma_0$.}
\end{remark}

\begin{remark} \label{Rem-SobolevLocal}
\emph{Given any $\varepsilon > 0$, the Sobolev condition $$\sup_{j \in \Z} \big\| \varphi_0^2 \dot{m}(2^j \cdot) \big\|_{\mathsf{H}_{\frac{d}{2} + \varepsilon}^2(\R^d)} \, < \, \infty$$ also suffices for the local form of Theorem A. Indeed, we just need to observe that both Lemmas \ref{SobVsClassical} and \ref{LemmaSobolevInversion} do provide upper bounds of the $\mathrm{W}_{d,\varepsilon}^2(\R^d)$-norms in terms of the smaller Sobolev norms above. Our Calder\'on-Zygmund type estimates in Proposition \ref{CZProposition} |needed to bound the homogeneous twisted multiplier| are also given in terms of Sobolev norms. This Sobolev condition is standard in Euclidean harmonic analysis and less demanding than the Mikhlin one.} 
\end{remark}

\begin{remark}
\emph{We also deduce the following result of independent interest. Let $(\varphi_j)_{j \in \Z}$ be a Littlewood-Paley partition of unity in $\R^{d}$ of the form \eqref{Eq-LP}. Set $\psi_j(g) = \varphi_j(g-e)$ for $g \in \G = \SL$. Let us write $\Psi_j$ for the Fourier multiplier associated to the symbol $\psi_j$. Let $M \in \Z$. Then, the following holds for $1 < p < \infty$
\begin{eqnarray*}
\Big\| \sum_{j \le M} \Psi_j(f) \otimes \delta_j \Big\|_{L_p(\V; RC_p)} \, \le_{\mathrm{cb}} \, C_p(M) \hskip1pt \| f \|_{L_p(\V)}. 
\end{eqnarray*}
This becomes a cb-norm equivalence when $\mathrm{supp} \hskip1pt \widehat{f} \subset \Omega = \bigcup_{j \le M-1} \mathrm{supp} \hskip1pt \psi_j$. In other words, a \emph{local form of the Littlewood-Paley theorem} in the group algebra of $\SL$.}
\end{remark}

\dem According to the noncommutative Khintchine inequality \cite{L,LP}, the $RC_p$ norm/square functions in the statement can be linearized and rewritten as follows 
\begin{eqnarray*}
\Big\| \sum_{j \le \mathrm{M}_\Omega} \Psi_j(f) \otimes \delta_j \Big\|_{L_p[RC_p]} \!\!\! & \simeq_{\mathrm{cb}} & \!\!\! \mathbf{E}_\varepsilon \Big\| \sum_{j \le \mathrm{M}_\Omega} \varepsilon_j \Psi_j(f) \Big\|_p \\ \!\!\! & \le_{\mathrm{cb}} & \!\!\! \sup_{\varepsilon_j = \pm 1} \Big\| \sum_{j \le \mathrm{M}_\Omega} \varepsilon_j \Psi_j(f) \Big\|_p \le_{\mathrm{cb}} 2 \sup_{A \subset \Z_\Omega} \big\| \Psi_A(f) \big\|_p,
\end{eqnarray*}
where $\Z_\Omega = \{j \in \Z: j \le \mathrm{M}_\Omega\}$ and $\Psi_A$ has symbol $\psi_A = \sum_{j \in A} \psi_j$. By the local form of Theorem A, the upper estimate will follow if the Euclidean symbols $\psi_A(\xi+e)$ satisfy the Mikhlin regularity imposed there with HM-constants uniformly bounded in $A$. This is however standard for Littlewood-Paley radial decompositions and follows by construction. Therefore, it  remains to justify the lower estimate for $f$ with frequency support inside $\Omega$. Assume first that $f \in L_2(\V)$ and consider any other $f' \in L_2(\V)$. We obtain for $\Psi = \sum_j \Psi_j \otimes \delta_j$
$$\langle \Psi(f), \Psi(f') \rangle \, = \, \sum_{j \le \mathrm{M}} \langle \Psi_j(f), \Psi_j(f') \rangle \, = \, \sum_{j \le \mathrm{M}} \int_{\G} \psi_j^2(g) \widehat{f}(g) \overline{\widehat{f}'(g)} d\mu(g) \, = \, \langle f, f' \rangle$$
since $\sum_{j \le \mathrm{M}} \psi_j^2$ is identically $1$ in the support of $\widehat{f}$ by hypothesis. Now, by density we may assume that $f \in L_p \cap L_2$  and that its norm is nearly attained by duality against a norm $1$ element $f' \in L_q\cap L_2$. Altogether, we get the lower estimate since $\|f\|_p \approx_{\mathrm{cb}} \langle f, f' \rangle = \langle \Psi(f), \Psi(f') \rangle \lesssim_{\mathrm{cb}} \|\Psi(f)\|_p.$ This completes the proof.
\fin 

\subsection{The asymptotic condition}

The proof of Theorem A will be completed with a simple patching argument from its local form. The key point is to observe that condition \eqref{Eq-ThmA} in Theorem A implies that the symbol $m \hskip-2pt : \SL \to \C$ differs from $L_1(\SL)$ by a constant function. This follows in turn from the exponential nature of the metric and Weyl's integration formula. 

\begin{lemma} \label{Lem-Exponential}
Given $\phi \in \mathcal{C}^1(\SL \setminus \{e\})$ and $\beta>2$
$$\sup_{\mathrm{X} \in \mathfrak{sl}_n(\R)} L(g)^{\beta} \big| \partial_\mathrm{X} \phi (g) \big| \, \le \, 1 \ \Rightarrow \ L(g)^\beta |\phi(g) - \alpha| \le C_\beta$$ for some  $\alpha \in \C$. The supremum runs over all unit vectors in the Lie algebra. 
\end{lemma}

\dem We claim that $\phi - \alpha \in \mathcal{C}_0(\SL)$ for some $\alpha \in \C$. This claim gives the statement. Indeed, every $g \in \SL$ factorizes as $g = u \exp(s \mathrm{X})$ for some unit vector $\mathrm{X} \in \mathfrak{sl}_n(\R)$ and $u \in \mathrm{SO}(n)$. By assumption and K-biinvariance of $L$, we obtain
\begin{eqnarray*}
| \phi(g) - \alpha| \!\!\!\! & = & \!\!\!\! \Big| \sum_{k \ge 1} \phi \big( g \exp((k-1)\mathrm{X}) \big) - \phi \big( g \exp(k \mathrm{X}) \big) \Big| = \Big| \sum_{k \ge 1} \partial_\mathrm{X} \phi \big( g \exp(s_k \mathrm{X}) \big) \Big| \\
\!\!\!\! & \le & \!\!\!\! \sum_{k \ge 1} L \big(u \exp((s+s_k) \mathrm{X}) \big)^{-\beta} \hskip-3pt = \hskip-1pt \sum_{k \ge 1} e^{-\beta (s+s_k)} \lesssim L(\exp(s \mathrm{X}))^{-\beta} \hskip-2pt = \hskip-1pt L(g)^{-\beta}
\end{eqnarray*}
for some $s_k \in (k-1,k)$. Let us now justify the claim. Every $g \in \SL$ factorizes as $g = u_1 \exp(\mathrm{Z}) u_2$ with $\mathrm{Z}$ a diagonal matrix in $\mathfrak{sl}_n(\R)$ and $u_1, u_2 \in \mathrm{SO}(n)$. By the surjectivity of the exponential map onto $\mathrm{SO}(n)$, we get that $u_j = \exp (\mathrm{A}_j)$ for some skew-symmetric $\mathrm{A}_j \in \mathfrak{so}_n$ with $\|\mathrm{A}_j\| \le 2\pi$. Under this factorization, we have $L(g) = \exp (\|\mathrm{Z}\|)$ and we conclude that
\begin{eqnarray} \label{eq-A2}
\lefteqn{\hskip-40pt \Big| \phi \big( \exp(\mathrm{A}_1) \exp(\mathrm{Z}) \exp(\mathrm{A}_2) \big) - \phi \big( \exp (\mathrm{A}_1) \exp(\mathrm{Z}) \big) \Big|} \\ \nonumber
\hskip70pt & = & \|\mathrm{A}_2\| \, \Big| \partial_{\hskip-2pt \frac{\mathrm{A}_2}{\|\mathrm{A}_2\|}} \phi \big( \exp(\mathrm{A}_1) \exp(\mathrm{Z}) \exp (r \mathrm{A}_2) \big) \Big| \, \le \, 2\pi \exp(-\beta \|\mathrm{Z}\|)
\end{eqnarray}
for some $0 < r < 1$. Similarly, let us note that $\exp(\mathrm{A}_1) \exp(\mathrm{Z}) = \exp(\mathrm{Z}) \mathrm{w}$ for $\mathrm{w} = \exp(- \mathrm{Z}) \exp(\mathrm{A}_1) \exp(\mathrm{Z}) = \exp (\mathrm{Y})$ where $\mathrm{Y} = \exp(- \mathrm{Z}) \mathrm{A}_1 \exp(\mathrm{Z})$ belongs to $\mathfrak{sl}_n(\R)$. Therefore, the following identity holds for some $r \in (0,1)$
\begin{equation} \label{eq-A1}
\Big| \phi \big( \exp(\mathrm{A}_1) \exp(\mathrm{Z}) \big) - \phi \big(\exp(\mathrm{Z}) \big) \Big| = \|\mathrm{Y}\| \, \Big| \partial_{\hskip-2pt \frac{\mathrm{Y}}{\|\mathrm{Y}\|}} \phi \big( \exp(\mathrm{Z}) \exp (r \mathrm{Y}) \big) \Big|.
\end{equation}
Since $\|\mathrm{Y}\| \le 2\pi \exp(2\|\mathrm{A}\|)$ and $L(\exp(\mathrm{Z}) \exp (r \mathrm{Y})) = L(\exp(r \mathrm{A}_1) \exp(\mathrm{Z})) = e^{\|\mathrm{Z}\|}$, the above quantity is bounded by $2\pi \exp(-(\beta-2)\|\mathrm{Z}\|)$, which decreases to $0$ for any $\beta > 2$ as $\mathrm{Z} \to \infty$. According to \eqref{eq-A2} and \eqref{eq-A1}, it suffices to prove that $\phi - \alpha \in \mathcal{C}_0$ when restricted to diagonal matrices $\exp \mathrm{Z}$. To prove it, we identify diagonal matrices in $\mathfrak{sl}_n(\R)$ with $\R^{n-1}$ as follows 
$$\mathrm{Z} = \mathrm{diag}(z_1, z_2, \ldots, z_n) \stackrel{\Lambda}{\longmapsto} (z_1, z_2, \ldots, z_{n-1}) = z.$$
Consider the function $\rho(z) = \phi(\exp \mathrm{Z})$. If $\Lambda(\mathrm{U}) = u$, we get  
\begin{eqnarray*}
\partial_u \rho (z) & = & \lim_{s \to 0} \frac{\rho(z + su) - \rho(z)}{s} \\
& = & \lim_{s \to 0} \frac{\phi(\exp(\mathrm{Z} \exp(s \mathrm{U})) - \phi(\exp(\mathrm{Z}))}{s} \, = \,  \partial_\mathrm{U} \phi(\exp(\mathrm{Z})).
\end{eqnarray*} 
Since $\Lambda$ is a contraction, we deduce the following inequality for $\rho$
$$\sup_{\|u\| = 1}  \big| \partial_u \rho (z) \big| \le \sup_{\|\mathrm{X}\| = 1}  \big| \partial_\mathrm{X} \phi (\exp(\mathrm{Z})) \big| \le L(\exp(\mathrm{Z}))^{-\beta} \le \exp(-\beta \|z\|).$$ This readily implies that $\rho$ has a limit $\alpha$ at infinity and the same holds for $\phi$. \fin

\begin{remark} \label{Rem-Linear}
\emph{By linearity of Lie differentiation $\partial_{\mathrm{X}_1 + \mathrm{X}_2} = \partial_{\mathrm{X}_1} + \partial_{\mathrm{X}_2}$, the supremum in Lemma \ref{Lem-Exponential} may be replaced by a maximum over norm $1$ matrices $\mathrm{X} \in \mathfrak{sl}_n(\R)$ in the directions of a fixed  orthogonal basis of the Lie algebra. In particular, if condition \eqref{Eq-ThmA} holds, we may apply Lemma \ref{Lem-Exponential} to $\phi = d^\gamma m$ with $|\gamma| = [n^2 \hskip-2pt /2]$. By assumption, $\partial_\mathrm{X} \phi(g)$ decays as $L(g)^{-\beta}$ for $\beta = [n^2 \hskip-2pt /2]+1 > 2$. In addition, $\phi \in \mathcal{C}_0$ so that $\alpha=0$ and $L(g)^{|\gamma|+1} |d^\gamma m(g)| \le C_\gamma$. In other words, the asymptotic decay of derivatives of order $[n^2 \hskip-2pt /2]$ is the same as those of order $[n^2 \hskip-2pt /2]+1$. Iterating this argument, Lemma \ref{Lem-Exponential} finally applies to $(\phi,\beta) = (m, [n^2 \hskip-2pt /2] +1)$ up to a constant depending on $n$. Condition \eqref{Eq-ThmA} does not assume $m \in \mathcal{C}_0$ and the lemma gives a constant $\alpha \in \C$ with 
$$L(g)^{\sigma_n + 1} |m(g) - \alpha| \, \le \, C_n \quad \mbox{for} \quad \sigma_n = \Big[\frac{n^2}{2} \Big].$$ It explains how \eqref{Eq-ThmA} incorporates Lafforgue/de la Salle rigidity \cite{LdlS} in Theorem A.} 
\end{remark}

\begin{remark} \label{eq-Integrability}
\emph{Let $\sigma_n = [n^2/2]$. By Weyl's integration formula 
$$\int_{\SL} f(g) \, d\mu(g) \, = \, \int_{\mathrm{SO}(n) \times \mathfrak{a}_+ \times \mathrm{SO}(n)} f(k_1 \! \exp \mathrm{Z} \, k_2) \prod_{j > k} \sinh (\mathrm{Z}_j - \mathrm{Z}_k) \, dk_1 d\mathrm{Z} \hskip1pt dk_2,$$ we get $\mu(B_R) \approx \exp(\sigma_n R)$ for $$B_R = \{g \in \SL: \log L(g) \le R\}.$$ Thus, $\sigma_n$ is the critical integrability index for the metric $L$: $L^{-\sigma_n-\varepsilon} \in L_1(\SL)$.}
\end{remark}

\demAGral By Remark \ref{Rem-Linear}, we may assume $m \in \mathcal{C}_0(\SL \setminus \{e\})$, equivalently $\alpha = 0$. Let $\Gamma$ be a cocompact lattice in $\SL$ with fundamental domain $\Delta$. Consider a relatively compact neighborhood of the identity $\Omega$ containing the closure of $\Delta$. Given $\phi \in \mathcal{C}_c^{\infty,+}(\SL)$ supported in $\Omega$ and identically $1$ over $\Delta$, define \vskip-10pt $$\Phi_\gamma(g) = \frac{\phi(\gamma g)}{\sum_{\rho \in \Gamma} \phi(\rho g)} \quad \mbox{for each} \quad \gamma \in \Gamma.$$ It is clear by construction that the denominator above is greater or equal than 1 and the $\Phi_\gamma$'s form a smooth partition of unity in $\SL$ indexed by $\Gamma$. Let us decompose the symbol $m = \sum_\gamma m_\gamma$ accordingly. By the triangle inequality, it suffices to prove that
$$A_p(m) = \sum_{\gamma \in \Gamma} \big\| T_{m_\gamma} \hskip-2pt : L_p(\mathcal{L}(\SL)) \to L_p(\mathcal{L}(\SL)) \big\|_{\mathrm{cb}} \, \le \, C_p C_{\mathrm{hm}}.$$ 
\vskip-8pt \noindent Using a conjugation with the translation by $\gamma$, we may clearly replace $m_\gamma$ by its left translate $M_\gamma(g) = m(\gamma^{-1} g) \Phi_e(g) = m_\gamma(\gamma^{-1} g)$. Then, the local form of Theorem A yields for $\sigma_n = [n^2 \hskip-2pt / 2]$
\begin{eqnarray*}
A_p(m) \!\! & \le & \!\! C_p \max_{|\beta| \le \sigma_n+1} \sum_{\gamma \in \Gamma} \, \sup_{g \in \Omega} \hskip-3pt \weight{g}^{\hskip-3pt |\beta|} \hskip-2pt \big| d_g^\beta M_\gamma(g) \big| \\
\!\! & \lesssim & \!\! C_p \max_{|\beta| \le \sigma_n+1} \Big( \sup_{g \in \Omega} \hskip-3pt \weight{g}^{\hskip-3pt |\beta|} \hskip-2pt \big| d_g^\beta M_e(g) \big| + \sum_{\gamma \neq e} \, \sup_{g \in \Omega} \big| d_g^\beta M_\gamma(g) \big| \Big).
\end{eqnarray*} 
\vskip-8pt \noindent Next, Leibnitz rule and left invariance of Lie differentiation give 
$$\big| d_g^\beta M_\gamma(g) \big| \, \le \, \sum_{\rho \le \beta} \big| d_g^{\rho} m(\gamma^{-1} g)  d_g^{\beta - \rho} \Phi_e (g) \big| \, \lesssim \, \sum_{\rho \le \beta} \big| d_g^{\rho}m (\gamma^{-1} g) \big|.$$
\vskip-8pt \noindent In particular, combining the above estimates we get the expected inequality 
$$A_p(m) \, \le \, C_p \max_{|\beta| \le \sigma_n+1} \Big( \sup_{g \in \Omega} \hskip-3pt \weight{g}^{\hskip-3pt |\beta|} \hskip-2pt \big| d_g^\beta m(g) \big| + \sum_{\gamma \neq e} \, \sup_{g \in \Omega} \big| d_g^\beta m(\gamma g) \big| \Big) \, \le \, C_p C_{\mathrm{hm}}.$$
The last inequality follows from Remark \ref{eq-Integrability}. Namely, the proof there gives that the Lie derivatives $|d_g^\beta m(\gamma g)|$ are all dominated by $C_{\mathrm{hm}} L(\gamma g)^{-(\sigma_n + 1)}$. Using relative compactness of $\Omega$ and Weyl's integration formula as we did before, the above sum is dominated (up to absolute constants) by $C_{\mathrm{hm}}$. This completes the proof. \fin    

\section{\bf Proof of Theorem B} \label{SectNecessary}

The first rigidity theorems for Fourier multipliers in $\SL$ which are relevant for this paper can be traced back to \cite{CHH,H2,Oh}. Haagerup's paper \cite{H2} is particularly relevant since it proves that $\SL$ and $\mathrm{SL}_n(\Z)$ fail to be weakly amenable for $n \ge 3$. This was strengthened in \cite{dLdlS,LdlS} by disproving the CBAP for the noncommutative $L_p$ spaces over the group algebra of any lattice in $\SL$ when $|1/p - 1/2|$ is large enough in terms of the rank. In this section we prove Theorem B as stated in the Introduction, which strengthens in turn \cite{dLdlS,LdlS} with more demanding rigidity conditions for $L_p$-multipliers. 

The proof relies, as in \cite{dLdlS,LdlS} (or \cite{HdL,Laat} for other higher rank Lie groups), on the idea developped in Lafforgue's work \cite{Laff} on strong property (T): to prove first local H\"older rigidity on the level of compact Lie groups, and then to combine these local estimates to obtain estimates on the whole noncompact Lie group. The local H\"older rigidity results for the compact group $\mathrm{SO}(n)$ are the content of Proposition \ref{prop:SOn_coeff_biinvariant}, and their combination to explore the whole group are derived in Theorem  \ref{thm:multipliers_SO(n,1)} for $\mathrm{SO}(n,1)$ and at the end of the section for $\SL$.

\subsection{Composition of H\"older continuous functions}

If $J$ is an interval of $\R$ and $\alpha>0$, we denote by $\mathcal{C}^\alpha(J)$ the space of all functions $f \colon J \to \R$ which are $[\alpha]$ times differantiable, and whose $[\alpha]$-th derivative is $(\alpha - [\alpha])$-H\"older-continuous on every compact subset of $J$ where $[ \, \cdot \, ]$ is the integer part function. We also write ${\lceil \cdot \rceil}$ for the ceiling function. If $X$ is a Banach space we denote $\mathcal{C}^\alpha(J;X)$ the space of such functions with values in $X$. When $J$ is compact, $\mathcal{C}^\alpha(J;X)$ is a Banach space for the norm
\[ \|f\|_{\mathcal{C}^\alpha(J;X)} = \max \left\{ \max_{k \leq [\alpha]} \sup_{x \in J} \big\| \partial^k \hskip-1pt f(x) \big\|, \sup_{\begin{subarray}{c} x \neq y \\ x,y \in J \end{subarray}} \frac{\| \partial^{[\alpha]} \hskip-1pt f(x) - \partial^{[\alpha]} \hskip-1pt f(x)\|}{|x-y|^{\alpha - [\alpha]}}\right\}.\]
When no confusion is possible, we simply write $\|f\|_{\mathcal{C}^\alpha}$ for this norm. In general the family of seminorms $\|\cdot \|_{\mathcal{C}^\alpha(K;X)}$ for $K$ a compact subinterval of $J$ turn $\mathcal{C}^\alpha(J;X)$ into a Fr\'echet space. It is clear that the spaces $\mathcal{C}^\alpha(J; X)$ are invariant under precomposition by a sufficiently smooth function. When $\alpha < 1$, $f \in \mathcal{C}^\alpha(J;X)$ and $\varphi \colon I \to J$ is a $\mathcal{C}^1$ function, the following inequality holds for every $x,y \in I$
\begin{equation}\label{eq:hoelder_composition} 
\big\| f \circ \varphi(x) - f \circ \varphi(y) \big\| \le \|f\|_{\mathcal{C}^\alpha} \|\varphi'\|_\infty^\alpha |x-y|^\alpha.
\end{equation}
We shall need quantitative estimates for higher derivatives, that we now collect.

\begin{proposition}\label{prop:Calpha_under_precomposition}
Given $\alpha > 0$ and two compact intervals $I,J$ of $\R$, assume that $\varphi \colon I \to J$ is a function of class $\mathcal{C}^{\lceil \alpha \rceil}$. Then, the map $f \mapsto f \circ \varphi$ maps $\mathcal{C}^{\alpha}(J;X)$ into $\mathcal{C}^{\alpha}(I;X)$. More precisely, if $\alpha\geq 1$, there is a constant $C_\alpha>0$ such that, for every such $I,J,f,\varphi$, every $1 \leq k \leq [\alpha]$ and every $x \in I$, we have
\begin{equation}\label{eq:derivative_from_FaaDiBruno} 
\big\| \partial^k (f \circ \varphi)(x) \big\| \, \leq \, C_\alpha \|f\|_{\mathcal{C}^k} \max_{1 \leq j \leq k} \big| \partial^j \varphi (x) \big|^{\frac k j}.
\end{equation}
Moreover, given $x,y \in I$ we also get the following inequality
\begin{multline} \label{eq:hoelder_from_FaaDiBruno} 
\big\| \partial^{[\alpha]} (f \circ \varphi) (x) - \partial^{[\alpha]} (f \circ \varphi)(y) \big\| \\ \null \hskip30pt \leq \, C_\alpha \|f\|_{\mathcal{C}^\alpha} \left( \|\varphi'\|_\infty^{\alpha} |x-y|^{\alpha - [\alpha]} + \max_{1 \leq j \leq \lceil \alpha \rceil} \Big( \|\partial^j \varphi \|_\infty^{\lceil \alpha \rceil/j} \Big) |x-y|\right).
\end{multline}
In particular, if $\alpha \notin \Z$, the $(\alpha - [\alpha ])$-H\"older constant of $\partial^{[\alpha]}(f \circ \varphi)$ at $x$ 
\begin{equation}\label{eq:hoelder2_from_FaaDiBruno} 
\limsup_{y \to x} \frac{\| \partial^{[\alpha]}(f \circ \varphi)(x) - \partial^{[\alpha]}(f \circ \varphi)(y)\|}{|x-y|^{\alpha - [\alpha]}} \le C_\alpha \|f\|_{\mathcal{C}^\alpha} |\varphi'(x)|^{\alpha}.
\end{equation}
\end{proposition}

\dem Fa\`a di Bruno's formula asserts for $k \leq [ \alpha ]$ that
\[ \partial^k (f \circ \varphi)(x) = \sum_{j=1}^k B_{k,j} \big( \varphi'(x),\dots, \partial^{k-j+1} \varphi(x) \big) (\partial^j \hskip-1pt f \circ \varphi)(x). \]
The coefficients in the above sum are given by the Bell polynomials
\[ B_{k,j}(z_1,z_2,\dots,z_{k-j+1}) \, = \, \sum \frac{k!}{i_1!i_2!\cdots i_{k-j+1}!} \prod_{s=1}^{k-j+1} \left(\frac{z_s}{s!} \right)^{i_s}, \]
where the sum is over all sequences $i_1,\dots,i_{k-j+1}$ of non-negative integers such that $i_1+i_2+\dots+i_{k-j+1} = j$ and $i_1+2 i_2+\dots+(k-j+1)i_{k-j+1} = k$. Elementary computations provide a constant $C_k$ such that
\[ \big| B_{k,j}(z_1,\dots, z_{k-j+1}) \big| \leq C_k \max \Big\{ |z_1|^k, |z_2|^{\frac k 2}, \dots, |z_{k-j+1}|^{\frac{k}{k-j+1}} \Big\}.\]
This immediately implies the inequality \eqref{eq:derivative_from_FaaDiBruno} in the statement.

Inequality \eqref{eq:hoelder2_from_FaaDiBruno} follows from \eqref{eq:hoelder_composition} or \eqref{eq:hoelder_from_FaaDiBruno}, according to the value of $\alpha$. To prove \eqref{eq:hoelder_from_FaaDiBruno}, we consider Fa\`a di Bruno's formula for $k = [ \alpha ]$. The term $j=[\alpha]$ is $\varphi'(x)^{[\alpha]} \partial^{[\alpha]} f (\varphi(x))$. In particular, the difference $\partial^{[\alpha]} (f \circ \varphi)(x) - \partial^{[\alpha]}(f \circ \varphi)(y)$ yields a term $j=[\alpha]$ given by
$$\varphi'(x)^{[\alpha]} \Big( \partial^{[\alpha]} f (\varphi(x)) - \partial^{[\alpha]} f(\varphi(y)) \Big) + \Big( \varphi'(x)^{[\alpha]} - \varphi'(y)^{[\alpha]} \Big) \partial^{[\alpha]} f (\varphi(y)) = \mathrm{A} + \mathrm{B}.$$ These two terms are respectively bounded above as follows
\begin{eqnarray*}
\mathrm{A} \!\!\! & \le & \!\!\! | \varphi'(x) |^{[\alpha]} \|f\|_{\mathcal{C}^\alpha} |\varphi(x)-\varphi(y)|^{\alpha - [\alpha]} \leq \|f\|_{\mathcal{C}^\alpha} \|\varphi'\|_\infty^{\alpha} |x-y|^{\alpha - [\alpha]}, \\
\mathrm{B} \!\!\! & \le & \!\!\! \| ((\varphi')^{[\alpha]})'\|_\infty \|f\|_{\mathcal{C}^\alpha} |x-y| \le C_\alpha \|f\|_{\mathcal{C}^\alpha} \|\varphi'\|_\infty^{[\alpha]-1} \|\varphi''\|_\infty |x-y|.
\end{eqnarray*}
The other terms in Fa\`a di Bruno's formula for $1 \le j < [\alpha]$ can be split into two terms $\mathrm{A}_j + \mathrm{B}_j$ as above. The terms $\mathrm{A}_j$ can be uniformly bounded, as for \eqref{eq:derivative_from_FaaDiBruno}, by
\[ C_\alpha \|f\|_{\mathcal{C}^\alpha} \max_{1 \leq j \leq [\alpha]} \| \partial^j \varphi \|_\infty^{\frac {[\alpha]} j} |x-y|.\]
Letting $\beta_{k,j}^\varphi(x) = B_{k,j}(\varphi'(x), \ldots, \partial^{k-j+1} \varphi (x))$, we have $$|\mathrm{B}_j| = \big| \big( \beta_{[\alpha],j}^\varphi(x) - \beta_{[\alpha],j}^\varphi(y) \big) ( \partial^j \hskip-1pt f \circ \varphi)(x)\big| \le \|f\|_{\mathcal{C}^\alpha} \big\| (\beta_{[\alpha],j}^\varphi)' \big\|_\infty |x-y|.$$ It is an straightforward exercise to show by direct calculation that $$\big| (\beta_{k,j}^\varphi)'(x) \big| \le C_\alpha B_{k+1,j} \big( |\varphi'(x)|, \ldots, | \partial^{k-j+2} \varphi (x)|\big).$$ Thus we get $|\mathrm{B}_j| \le \displaystyle C_\alpha \|f\|_{\mathcal{C}^\alpha} \max_{1 \leq j \leq \lceil \alpha \rceil} \| \partial^j \varphi \|_\infty^{\frac {\lceil \alpha \rceil} j} |x-y|$. 
This proves \eqref{eq:hoelder_from_FaaDiBruno}. \fin

\subsection{$S_p$-multipliers: Estimates on $\mathrm{SO}(n)$}

Given a locally compact group $\G$ and $1 \leq p \leq \infty$, we say that a bounded measurable function $m \colon \G \to \C$ is an \emph{$S_p$-multiplier} (resp. $S_p$-$S_q$-multiplier) if the map
\[ S_m\colon \Big( a_{g,h} \Big)_{g,h \in \G} \mapsto \Big( m(gh^{-1}) a_{g,h} \Big)_{g,h \in \G} \]
is bounded on $S_p(L_2(\G))$ (resp. from $S_p(L_2(\G))$ to $S_q(L_2(\G))$).

Let us identify $\mathrm{SO}(n-1)$ with the subgroup of $\mathrm{SO}(n)$ fixing the first coordinate vector $e_1$ of $\R^n$. Then, the double quotient $\mathrm{SO}(n-1) \backslash \mathrm{SO}(n) / \mathrm{SO}(n-1)$ identifies with $[-1,1]$ through $$\mathrm{SO}(n-1) k \mathrm{SO}(n-1)\mapsto k_{1,1}.$$ Therefore, to an $\mathrm{SO}(n-1)$-biinvariant function $\varphi \colon \mathrm{SO}(n) \to \C$ corresponds a unique function $\widetilde \varphi \colon [-1,1] \to \C$ satisfying $\varphi(k) = \widetilde \varphi(k_{1,1})$. In what follows, we shall fix $n \ge 3$, $p>2+\frac{2}{n-2}$ and set 
\begin{equation} \label{eq-alpha}
\alpha_0 = \frac{n-2}{2} -\frac{n-1}{p}>0 \quad \mbox{and} \quad \alpha = \begin{cases} \alpha_0&\textrm{if }{\alpha_0 \notin \Z}\\ \alpha_0 - \varepsilon& \textrm{if }\alpha_0 \in \Z\end{cases}
\end{equation}
for an arbitrarily small $\varepsilon$. The following is a strengthening of \cite[Proposition 3.1]{dLdlS}.

\begin{proposition}\label{prop:SOn_coeff_biinvariant} 
Let $p$ and $\alpha$ be as above for any $\varepsilon$. Assume that $\varphi \colon \mathrm{SO}(n) \to \C$ is an $\mathrm{SO}(n-1)$-biinvariant $S_p$-multiplier in $\mathrm{SO}(n)$. Then, we get $\widetilde \varphi \in \mathcal{C}^{\alpha}((-1,1))$.
\end{proposition}

\begin{remark}
\emph{The proof will actually show more:} 
\begin{itemize}
 \item \emph{If $\alpha_0 \in \Z$, $\widetilde \varphi$ is $\alpha_0-1$ times differentiable and
  \[ \big| \partial^{\alpha_0-1} \widetilde \varphi (x) - \partial^{\alpha_0-1} \widetilde \varphi (y) \big| \lesssim |x-y| \hskip1pt \big| \log |x-y| \big|^{\frac 1 p}\]
  holds uniformly on every compact subset of the interval $(-1,1)$.}

\vskip3pt

\item \emph{The conclusion holds if one merely assumes that $\varphi$ is an $S_p$-$S_\infty$-multiplier.}
\end{itemize}
\end{remark}

We shall prove a dual statement. Let $\sphere^{n-1}$ denote the unit sphere in $\R^n$ equipped with the Lebesgue probability measure. For $\delta \in [-1,1]$, let $T_\delta$ be the (densely defined) operator on $L_2(\sphere^{n-1})$ given by
\[T_\delta f(x)=\textrm{the average of }f\textrm{ on } \big\{ y \in \sphere^{n-1} \mid \langle x,y\rangle = \delta \big\}.\]
Equivalently, using the identification $\sphere^{n-1} \cong \mathrm{SO}(n-1) \backslash \mathrm{SO}(n)$ through the map $\mathrm{SO}(n-1) g \mapsto g^{-1} e_1$, we can consider $L_2(\sphere^{n-1})$ as a subspace of $L_2(\mathrm{SO}(n))$. Then $T_\delta$ is the operator on $L_2(\mathrm{SO}(n))$ equal to 
\begin{equation}\label{eq=def_Tdelta} \int_{\mathrm{SO}(n-1) \times \mathrm{SO}(n-1)} \lambda( u g u') du du' \in \mathcal{B}(L_2(\mathrm{SO}(n)))\end{equation}
for $g \in \mathrm{SO}(n)$ satisfying $g_{11} = \delta$. Here, $\lambda$ denotes the left-regular representation.

\begin{proposition} 
The map $\delta \in (-1,1)\mapsto T_\delta$ belongs to $\mathcal{C}^{\alpha}((-1,1); S_p(L_2(\mathrm{SO}(n))))$.
\end{proposition}

\begin{remark}
\emph{This implies Proposition \ref{prop:SOn_coeff_biinvariant} because for an $S_p$-multiplier $\varphi$ (or more generally an $S_p$-$S_\infty$ multiplier), we have $S_\varphi(T_\delta) = \widetilde \varphi(\delta) T_\delta$ and in particular $\widetilde \varphi(\delta) = \langle S_\varphi(T_\delta) \xi,\xi\rangle$ where $\xi  \in L_2(\mathrm{SO}(n))$ is the constant function equal to $1$. So the function $\widetilde \varphi$, which is the composition of $(\delta \mapsto T_\delta)$ with the continuous linear map $T \in S_p \mapsto \langle S_\varphi(T) \xi,\xi\rangle$, is at least as regular as $(\delta \mapsto T_\delta)$.}
\end{remark}

\dem As explained in \cite[Lemma 3.2]{dLdlS}, there is an orthonormal basis in which the operators $T_\delta$ are all diagonal, and in which the eigenvalue sequence is $(\varphi_k(\delta))_{k \geq 0}$ with multiplicity \[ m_k = \frac{(n+k-3)!(n+2k-2)}{(n-2)! k!},\]
  where
  \[ \varphi_k(x) = c_n \int_0^\pi \big( x+i\sqrt{1-x^2} \cos \theta \big)^k (\sin \theta)^{n-3} \, d\theta \qquad \mbox{and} \qquad c_n = \frac{\Gamma(\frac{n-1}{2})}{\sqrt \pi \Gamma(\frac{n-2}{2})}. \] 
By derivating in the integral, we obtain  
  \[ \big| \partial^r \hskip-1pt \varphi_k (x) \big| \leq C(n,r) \frac{(1+k)^r}{\big(\sqrt{1-x^2} \hskip2pt \big)^{2r-1}} \int_0^\pi \big| x + i \sqrt{1-x^2} \cos \theta \big|^{k-r} (\sin \theta)^{n-3} \, d \theta\]
for any nonegative integer $r \le k$. When $r > k$, we must replace the exponent $k-r$ inside the integral by $0$. As (3.2) in \cite{dLdlS}, there is a constant $C(n,r)$ (depending on $n$ and $r$) such that this is less than
  \[ \frac{C(n,r)}{\big( (1+k)(1-x^2) \big)^{\frac{n-2}{2}}} \frac{(1+k)^r}{\big(\sqrt{1-x^2} \hskip2pt \big)^{2r-1}}.\]
    In particular, if $0 < c < 1$, there exists $C'(n,r)$ such that
    \[ | \partial^r \hskip-1pt \varphi_k (x) | \leq C'(n,r) (1+k)^{r+1-\frac n 2} \quad \mbox{for every} \quad x \in [-c,c].\] 
    And so, bounding $m_k \leq A(n) (1+k)^{n-2}$, we obtain that for every such $x$
    \[ \sum_{k \ge 0} m_k | \partial^r \hskip-1pt \varphi_k (x)|^p \leq AC' \sum_{k \ge 0} (1+k)^{n-2+p(r+1-\frac n 2)}  = AC' \sum_{k \ge 0} (1+k)^{p(r-\alpha_0)-1},\]
which converges if $r<\alpha_0$. Since $r \in \Z$, this holds iff $r < \alpha$. Taking $r=[\alpha]$, we deduce that $\delta \mapsto T_\delta$ belongs to $\mathcal{C}^{[\alpha]}([-c,c]; S_p(L_2(\mathrm{SO}(n))))$ for all $c<1$. Its $[\alpha]$-th derivative $\partial^{[\alpha]} T_x$ is the operator which, in the basis as above, is diagonal with eigenvalues $( \partial^{[\alpha]} \hskip-1pt \varphi_k (x))_k$ and multiplicities $(m_k)_k$. We get
\[ \big\| \partial^{[\alpha]} T_x - \partial^{[\alpha]} T_{y} \big\|_p^p = \sum_{k \ge 0} m_k \big| \partial^{[\alpha]} \hskip-1pt \varphi_k(x) - \partial^{[\alpha]} \hskip-1pt \varphi_k(y) \big|^p.\]
As above, we bound $m_k \leq A(n) (1+k)^{n-2}$. In addition, when $x,y \in [-c,c]$, we estimate the difference $| \partial^{[\alpha]} \hskip-1pt \varphi_k(x) - \partial^{[\alpha]} \hskip-1pt \varphi_k(y)|$ with two bounds. If $k |x-y| \leq 1$, we use 
\[ \big| \partial^{[\alpha]} \hskip-1pt \varphi_k(x) - \partial^{[\alpha]} \hskip-1pt \varphi_k(y) \big| \le \max_{|z| \le c} |\partial^{\lceil \alpha \rceil} \varphi_k(z)| \hskip1pt |x-y| \leq C' (1+k)^{\lceil \alpha \rceil+1-\frac n 2} |x-y|.\]
If $k |x-y|>1$, we use 
\[ \big| \partial^{[\alpha]} \hskip-1pt \varphi_k(x) - \partial^{[\alpha]} \hskip-1pt \varphi_k(y) \big| \le 2 \max_{|z| \le c} |\partial^{[\alpha]} \varphi_k(z)| \leq 2C' (1+k)^{[\alpha]+1-\frac n 2}.\]
When $\alpha_0 \notin \Z$, we get $\alpha = \alpha_0$ and obtain 
\begin{eqnarray*}
\big\| \partial^{[\alpha]} T_x - \partial^{[\alpha]} T_{y} \big\|_p^p \!\! & \lesssim & \!\! \sum_{k \leq \frac{1}{|x-y|}} (1+k)^{p(\lceil \alpha \rceil - \alpha_0)-1} |x-y|^p \\ 
\!\! & + & \!\! \sum_{k >\frac{1}{|x-y|}} (1+k)^{p([\alpha] - \alpha_0)-1} \lesssim |x-y|^{p(\alpha - [\alpha])}.
\end{eqnarray*}
Additionally, when $\alpha_0 \in \Z$, we get $[\alpha] = \alpha_0-1$ and the same estimate gives  
\[ \big\| \partial^{[\alpha]} T_x - \partial^{[\alpha]} T_y \big\|_p \lesssim |x-y| \hskip1pt |\log |x-y||^{\frac 1 p}.\]
Therefore, $\partial^{[\alpha]} T_x$ is $(\alpha - [\alpha]$)-H\"older continuous on every compact subinterval. \fin

\subsection{Rigidity for $\mathrm{K}$-biinvariant $S_p$-multipliers on $\mathrm{SO}(n,1)$}

Theorem B will be deduced from Propositions \ref{prop:Calpha_under_precomposition} and \ref{prop:SOn_coeff_biinvariant} in the next subsection. Before that we explain, in a simpler situation, how the same idea as for Theorem B allows us to prove some rigidity for $\mathrm{K}$-biinvariant $S_p$-multipliers on the rank $1$ simple Lie group $\mathrm{SO}(n,1)$ which contain $\mathrm{SO}(n)$ as a subgroup. Of course, the results cannot be as strong as for $\SL$ in the sense that they cannot prescribe any rate of convergence at infinity as $\mathrm{SO}(n,1)$ is weakly amenable. But it turns out that there are some ``higher order'' rigidity, which only appears at the level of the derivatives. To our knowledge, Theorem \ref{thm:multipliers_SO(n,1)} and its particular case Remark \ref{rem:coeff_SO(n,1)_rep} for $\mathrm{K}$-biinvariant matrix coefficients of $\mathrm{SO}(n,1)$ is the first application to rank $1$ groups of the ideas around (the proof of) strong property (T) originating in \cite{Laff}.

Recall that $\mathrm{SO}(n,1)$ the group of $(n+1) \times (n+1)$ matrices of determinant one and preserving the symmetric bilinear form
\[ \big[ (x_1,\dots,x_{n+1}),(y_1,\dots,y_{n+1}) \big] = \sum_{i=1}^n x_i y_i - x_{n+1}y_{n+1}.\]
Denote by $\mathrm{K} \simeq \mathrm{O}(n)$ the maximal compact subgroup of $\mathrm{SO}(n,1)$ given by
\[ \mathrm{K} = \left\{ \begin{pmatrix} U & 0 \\ 0 & \det(U)\end{pmatrix} \colon \ U \in \mathrm{O}(n) \right\},\]
and by $\mathrm{A} \simeq \R$ the group
\[ \mathrm{A} = \left\{ D(s) = \begin{pmatrix} \cosh(s) & 0 & \sinh(s) \\ 0 & 1_{n-1} &0 \\ \sinh(s) & 0 & \cosh(s)\end{pmatrix} \colon \ s \in \R\right\}.\]
Denote also by $\mathrm{A}_+$ the subset corresponding to $s \geq 0$. The polar decomposition in $\mathrm{SO}(n,1)$ reads as $\mathrm{SO}(n,1) = \mathrm{KA_+K}$. That is, every element of $\mathrm{SO}(n,1)$ can be written as $g=kak'$ for $k,k' \in \mathrm{K}$ and $a \in \mathrm{A}_+$. Moreover, $a=D(s)$ is uniquely determined by $\|g\| = \|a\|$, that is $\|g\|=e^s$. Alternatively, $a=D(s)$ is uniquely determined by $\mathrm{tr}(g^* g) = \mathrm{tr}(a^*a) = n-1+2 \cosh(2s)$.

In particular, every $\mathrm{K}$-biinvariant function $m \colon \mathrm{SO}(n,1) \to \C$ can be written as $m(g) = \varphi(\mathrm{tr}(g^*g))$ for a function $\varphi \colon [n+1,\infty) \to \C$. So the next result gives regularity properties for $\mathrm{K}$-biinvariant multipliers of $\mathrm{SO}(n,1)$.

\begin{theorem}\label{thm:multipliers_SO(n,1)}
Let $\alpha = \alpha(\varepsilon)$ be defined as in \eqref{eq-alpha} for some $n\geq 3$ and $p>2+\frac{2}{n-2}$. Then, every $\mathrm{K}$-biinvariant $S_p$-$S_\infty$ Schur multiplier of $\G = \mathrm{SO}(n,1)$ is of class $\mathcal{C}^{\alpha}$. 

More precisely, let $m(g) = \varphi(\mathrm{tr}(g^* g))$ for a function $\varphi \colon [n+1,\infty) \to \C$. Assume that the Schur multiplier $S_m(g,h) = m(gh^{-1})$ is $S_p$-$S_\infty$ bounded. Then $\varphi$ is of class $\mathcal{C}^{\alpha}(n+1,\infty)$, and the following local/asymptotic estimates hold:
\begin{itemize}
\item[i)] Given $x > n+1$ and an integer $1 \le k \le [\alpha]$ $$\big| \partial^k \varphi (x) \big| \le C_{p,n}^\varepsilon \frac{\|S_m\|_{\mathcal{B}(S_p(L_2(\mathrm{G})),S_\infty(L_2(\mathrm{G}))}}{(x - n-1)^k}.$$

\item[ii)] The H\"older constants in a neighborhood of $x$ 
  $$\limsup_{y \to x} \frac{|\partial^{[\alpha]} \varphi(x) - \partial^{[\alpha]} \varphi(y)|}{|x-y|^{\alpha - [\alpha]}} \, \le \, C_{p,n}^\varepsilon \frac{\|S_m\|_{\mathcal{B}(S_p(L_2(\mathrm{G})),S_\infty(L_2(\mathrm{G}) )}}{(x-n-1)^{\alpha}},$$
\end{itemize}
 \end{theorem}

\begin{remark} \label{rem:coeff_SO(n,1)_rep}
\emph{Coefficients of uniformly bounded representations on Hilbert spaces are particular cases of $S_\infty$-multipliers. So for $p=\infty$, the previous theorem has as a consequence that $\mathrm{K}$-biinvariant coefficients of uniformly bounded representations of $\mathrm{SO}(n,1)$ are of class $\mathcal{C}^{\frac n 2-1}$ if $n$ is odd (and of class $\mathcal{C}^{\frac n 2-1-\varepsilon}$ for every $\varepsilon>0$ if $n$ is even), and the derivatives and H\"older constants of their restrictions to $\mathrm{A}$ are explicitly controlled by the bounds of Theorem \ref{thm:multipliers_SO(n,1)}. We are not aware of any such result in the literature, even for unitary representations. We recall however that $\mathrm{K}$-biinvariant coefficients of irreducible unitary representations (and more generally $\mathrm{K}$-finite coefficients of admissible representations) are $\mathcal{C}^\infty$ \cite{HC0}, and that unitary representations are direct integrals of irreducible representations. But Harish-Chandra's estimates depend on the representation and therefore do not provide higher order regularity estimates for arbitrary unitary representations.}
\end{remark}

\demC For $r>0$, the function $$k \in \mathrm{SO}(n) \mapsto m \Big(D(r) \begin{pmatrix} k & 0\\0&1\end{pmatrix} D(r) \Big)$$ is an $S_p$-$S_\infty$-multiplier of norm $\leq 1$ by restriction, and is $\mathrm{SO}(n-1)$-biinvariant because $D = D(r)$ commutes with the image of $\mathrm{SO}(n-1)$ in $\mathrm{K}$.

In particular, if we consider the rotation matrix $k_\delta$ of angle $\arccos \delta$ in the space spanned by the first two coordinate vectors
\[k_\delta = \begin{pmatrix} \delta & -\sqrt{1-\delta^2} &0 \\ \sqrt{1-\delta^2} & \delta & 0 \\ 0 & 0 & 1_{n-2} \end{pmatrix} \quad \mbox{for} \quad \delta \in [0,1],\]
Proposition \ref{prop:SOn_coeff_biinvariant} gives that the function $\psi_r \colon \delta \mapsto m( D k_\delta D)$ is of class $\mathcal{C}^\alpha$ (uniformly in $r$ and $\varphi$) and we obtain
\begin{equation}\label{eq:psir_in_Calpha_SOn1}\sup_{r>0} \|\psi_r\|_{C^\alpha([0,1-\eta])}<\infty\textrm{ for every }\eta>0.\end{equation}
We can compute
\[ \mathrm{tr} \big( (D k_\delta D)^* (D k_\delta D) \big) = 4 (\sinh r)^4 \delta^2+2 (\sinh 2r)^2\delta +n-3+4 (\cosh r)^4  = a_r\delta^2+ b_r \delta+c_r.\]
Denote by $g_r$ the inverse of $\delta \in [0,1] \mapsto \mathrm{tr}((D(r) k_\delta D(r))^* (D(r) k_\delta D(r)))$. That is
\[g_r(x) = \frac{-b_r}{2a_r} + \sqrt{\frac{b_r^2}{4 a_r^2} + \frac{x-c_r}{a_r}}\]
for every $x \in [c_r,a_r+b_r+c_r]$.
Observe that
\[ \frac{b_r}{a_r} = \frac{(\sinh 2r)^2}{2 (\sinh r)^4} \geq 2\]
Therefore, for every integer $k\geq 1$, we have
\[ \max_{1 \leq j \leq k} \big\| \partial^j g_r \big\|_\infty^{\frac k j} \leq C_k b_r^{-k}.\]
Let $r>0$ be determined by $x = c_r$. Then Proposition \ref{prop:Calpha_under_precomposition} gives,
\[ \big| \partial^k \varphi(x) \big| \leq C b_r^{-k} \quad \mbox{for every
$x \in [c_r,a_r+b_r+c_r]$ and every integer $k \leq [\alpha]$}\]
and the $(\alpha - [\alpha])$-H\"older constant of $\varphi$ at $x$ is $\leq C b_r^{-\alpha}$. This proves the theorem as there is a constant $C'$ such that $1/C' b_r \leq  c_r - n -1 \leq C' b_r$. \fin

\begin{remark}
\emph{Theorem \ref{thm:multipliers_SO(n,1)} implies a cheap and weaker form of Theorem B. Indeed, when $n \geq 4$ and $m$ and $\varphi$ are given as in Theorem B, then by restriction $m$ defines a $\mathrm{K}$-biinvariant multiplier on $\mathrm{SO}(n-1,1)$. In particular, Theorem \ref{thm:multipliers_SO(n,1)} implies that $\varphi$ is $\mathcal{C}^{\beta-}$ for $\beta = \frac{n-3}{2} - \frac{n-2}{p}$ and gives explicit estimates on the derivatives of the function $\varphi$. These estimates are weaker than the conclusion of Theorem B, as the regularity is lower and there are no additional factors $c_k$.}
\end{remark}

\subsection{Rigidity for radial $S_p$-multipliers on $\SL$}

Now we use our results so far to prove Theorem B. The idea is the same of for Theorem \ref{thm:multipliers_SO(n,1)}, but the details are more technical. In fact, we shall prove a form of Theorem B, which is slightly stronger in two senses:
\begin{itemize}
\item[i)] We shall allow Schur multipliers which are radial either in the normalized Hilbert-Schmidt norm $|g|^2 = \frac{1}{n} \mathrm{tr}(g^*g)$ or the operator norm $\|\cdot\|$ on $\mathrm{SL}_n(\R)$
\[ \hskip20pt \|g\| = \sup \Big\{ \sum_{i,j=1}^n g_{i,j} x_i y_j \, \big| \, x,y \in \R^n, \sum_i x_i^2 = \sum_j y_j^2=1 \Big\}. \] 
\item[ii)] The assumption that $S_m$ is $S_p$-bounded by the weaker assumption that $S_m$ is $S_p-S_\infty$ bounded, and with $\|S_m\|_{\mathcal{B}(S_p(L_2(\mathrm{G})))}$ replaced by the smaller quantity $\|S_m\|_{\mathcal{B}(S_p(L_2(\mathrm{G})),B(L_2(\mathrm{G})))}$.
\end{itemize}
Thus, let $\varphi \colon (0,\infty) \to \C$ be a function satisfying that $S_{\varphi(|\cdot|)}$ or $S_{\varphi(\|\cdot\|)}$ maps $S_p(L_2(G))$ to $B(L_2(G))$ with norm $1$. We start with the following crucial lemma where the letter $C$ stands for a constant depending on $p,n$ only.

\begin{lemma}\label{lemma:phi_C_alpha_using_SOn} 
The symbol $\varphi$ is of class $\mathcal{C}^{\alpha}$ on $(1,\infty)$. Moreover$\hskip1pt :$
  \begin{itemize}
  \item[i)] If $\alpha \leq 1$, then
    \[ | \varphi(x) -\varphi(y)| \leq \sup_{x \le z \le y} \frac{C}{\big( (z-1)z^{\frac{n}{n-2}} \big)^{\alpha}} |x-y|^\alpha\]
    for every pair $x,y \in (1,\infty)$ satisfying that $x \leq y \leq x^{ 1 + \frac{n}{n-2}}$.
  
\vskip3pt

  \item[ii)] If $\alpha >1$, then for every $1 \leq k \leq [ \alpha ]$,
    \[ | \partial^k \varphi(x)| \leq \frac{C}{(x-1)^k x^{\frac{n}{n-2}}}.\]

   \item[iii)] The $(\alpha - [ \alpha ])$-H\"older constant of $\partial^{[\alpha]} \varphi$ at $x$ satisfies $$\limsup_{y \to x} \frac{|\partial^{[\alpha]} \varphi(x) - \partial^{[\alpha]} \varphi (y)|}{|x-y|^{\alpha - [ \alpha ]}} \le \frac{C}{\big( (x-1)x^{\frac{n}{n-2}} \big)^{\alpha}}.$$
 \end{itemize}
\end{lemma}

\dem For $r>0$, define $s= - \frac{r}{n-1}$ and $D=\mathrm{diag}(e^r,e^s,\dots,e^s) \in \mathrm{SL}_n(\R)$. The function $k \in \mathrm{SO}(n) \mapsto \varphi (\| D k D \|)$ is an $S_p$-$S_\infty$-multiplier of norm $\leq 1$ by restriction to submatrices with entries $(g,h)$ in $D \hskip1pt \mathrm{SO}(n) \times D^{-1} \mathrm{SO}(n)$. Let us also recall that it is $\mathrm{SO}(n-1)$-biinvariant, because $D$ commutes with $\mathrm{SO}(n-1)$. Note additionally that the same properties hold for the normalized Hilbert-Schmidt norm $|\cdot|$.

In particular, if we consider the rotation matrix $k_\delta$ of angle $\arccos \delta$ in the space spanned by the first two coordinate vectors
  \[k_\delta = \begin{pmatrix} \delta & -\sqrt{1-\delta^2} &0 \\ \sqrt{1-\delta^2} & \delta & 0 \\ 0 & 0 & 1_{n-2} \end{pmatrix} \quad \mbox{for} \quad \delta \in [0,1],\] Proposition \ref{prop:SOn_coeff_biinvariant} gives that the function $\psi_r \colon \delta \mapsto \varphi (\| D k_\delta D \|)$ is of class $\mathcal{C}^\alpha$ uniformly in $r$ and $\varphi$ (the same holds one more time replacing the operator norm by the Hilbert-Schmidt norm) and we obtain
\begin{equation}\label{eq:psir_in_Calpha}\sup_{r>0} \|\psi_r\|_{\mathcal{C}^\alpha([0,1-\eta])}<\infty \textrm{ for every }\eta>0.\end{equation}

We can compute
\[ D k_\delta D = e^{r+s} \begin{pmatrix} e^{r-s}\delta & -\sqrt{1-\delta^2} &0 \\ \sqrt{1-\delta^2} & e^{s-r} \delta & 0 \\ 0 & 0 & e^{s-r} 1_{n-2} \end{pmatrix}.\]
We first consider radial multipliers in the operator norm $\|\cdot\|$. The matrix 
\[ \begin{pmatrix} e^{r-s}\delta & -\sqrt{1-\delta^2} \\ \sqrt{1-\delta^2} & e^{s-r} \delta \end{pmatrix}\]
has determinant $1$ and Hilbert-Schmidt norm $(2+4 \delta^2\sinh^2(r-s))^{1/2}$, so its norm is equal to $g(\delta \sinh(r-s))$ where $g(x) = (1+2x^2 + 2\sqrt{x^2+x^4})^{1/2}$. Therefore, $D k_\delta D$ has norm $e^{r+s} g \big(\delta \sinh(r-s) \big)$. We conclude that $$\psi_r(\delta) = \varphi \big( e^{r+s} g(\delta \sinh(r-s) \big).$$
Taking $\delta=1$ in the computation of the norm of $D k_\delta D$, it follows that $g(\sinh u) = e^u$ for every $u \geq 0$. In other words, the inverse is $g^{-1}(x) = \sinh(\log x) = \frac{1}{2} \left(x- \frac 1 x\right)$. So if we define $H_r \colon [e^{r+s},e^{2r}] \to [0,1]$ by
\begin{equation}\label{eq:def_Hr} H_r(x) = \frac{1}{2 \sinh(r-s)} \Big(\frac{x}{e^{r+s}} - \frac{e^{r+s}}{x}\Big) = \frac{\frac{x}{e^{r+s}} - \frac{e^{r+s}}{x}}{e^{r-s} - e^{s-r}},\end{equation}
we obtain that $e^{r+s} g(H_r(x) \sinh(r-s)) = x$ and 
\begin{equation}\label{eq:phi=psioH}
  \varphi = \psi_r \circ H_r \textrm{ on }[e^{r+s},e^{2r}].
\end{equation}

We are now in position to apply Proposition \ref{prop:Calpha_under_precomposition}, which gives us bounds on the derivatives of $\varphi$ and their H\"older constants in terms of the derivatives of $H_r$. So we compute, for $x \in [e^{r+s},e^{2r}]$,
\[ H_r'(x) = \frac{1}{e^{2r} - e^{2s}} \Big(1+\frac{e^{2r+2s}}{x^2}\Big) \in \Big[ \frac{1}{e^{2r} - e^{2s}},\frac{2}{e^{2r} - e^{2s}} \Big] \]
and
\[ \partial^j H_r(x) = \frac{(-1)^{j-1} j!}{(e^{-2s} - e^{-2r})x^{j+1}} \quad \mbox{for} \quad j \geq 2. \]

For a fixed $x$, $| \partial^j H_r(x)|$ is a decreasing function of $r$, so the bounds we get will be optimal when $r$ is maximal. In other words, when $x=e^{r+s}$. This determines the value of $r$ and $s$ as a function of $x$, so that $$e^{2r} = x^{1+\frac{n}{n-2}} \quad \mbox{and} \quad e^{2s}=x^{1-\frac{n}{n-2}}.$$ This yields that $|\partial^j \hskip-1pt H_r(x)| \approx \frac{1}{x-1}$ when $x$ is close to one and $|\partial^j \hskip-1pt H_r(x)| \approx x^{-(j + \frac{n}{n-2})}$ when $x$ is large. Written concisely, we get the following estimates for the derivatives of $H_r$ 
\begin{equation} \label{eq-derivatives Hr}
\big| \partial^j \hskip-1pt H_r(x) \big| \leq \frac{C}{(x-1)x^{j-1+\frac{n}{n-2}}} \ \Rightarrow \ \max_{1 \leq j \leq k} \big| \partial^j \hskip-1pt H_r(x) \big|^{\frac k j} \leq \frac{C}{(x-1)^k x^{\frac{n}{n-2}}}.
\end{equation}
Proposition \ref{prop:Calpha_under_precomposition} together with \eqref{eq:psir_in_Calpha}, \eqref{eq:phi=psioH} and \eqref{eq-derivatives Hr} implies that $\varphi$ is of class $\mathcal{C}^\alpha$ at $x$ and estimates ii) and iii) in the statement follow. When $\alpha < 1$, we use \eqref{eq:hoelder_composition} instead of Proposition \ref{prop:Calpha_under_precomposition} to deduce i). The condition $x \leq y \leq x^{ 1 + \frac{n}{n-2}}$ comes imposed by the domain of $H_r$. This proves the lemma for the operator norm. 

Let us now consider radial Schur multipliers in the normalized Hilbert-Schmidt norm $|\cdot|$. In that case, the norm computations are straightforward and lead to the decomposition
\begin{equation}\label{eq:def_Hrtilde} 
\varphi = \psi_r \circ \widetilde H_r \textrm{ on } [A_r,B_r]
\end{equation}
with $A_r^2 = \frac{1}{n} ((n-2) e^{4s} + 2 e^{2r+2s})$, $B_r^2 = \frac{1}{n} ((n-1) e^{4s} +  e^{4r})$ and 
  \[\widetilde H_r(x) = \sqrt{\frac{x^2/A_r^2-1}{B_r^2 /A_r^2-1}} = \frac{A_r}{\sqrt{B_r^2-A_r^2}} H(x/A_r).\]
  The function $H\colon x \mapsto \sqrt{x^2-1}$ has first derivative
  $x (x^2-1)^{-\frac 1 2}$ and (by induction) $j$-th derivative of the
  form $P_j(x) (x^2-1)^{\frac 1 2 - j}$ for certain polynomial $P_j$. This gives
\[ \partial^j \hskip-1pt \widetilde H_r(x)  =\frac {A_r^j}{\sqrt{(x+A_r)^{2j-1}}} P_j\Big(\frac x {A_r}\Big) \sqrt{\frac{1}{(B_r^2-A_r^2)(x-A_r)^{2j-1}}}.\]
To bound it when $x$ is close to $1$, we choose $r>0$ so that
$B_r=x$. We have that $x-1\approx c r^2$ and $B_r^2-A_r^2=dr^2$ for some $c,d>0$.
We clearly obtain for $x$ close to $1$
 \[ \max_{1 \leq j \leq k} \big| \partial^j H_r(x) \big|^{\frac k j} \leq \frac{C}{(x-1)^k}.\]
For $x$ large, we choose $r>0$ so that $x=e^{r+s}$. It follows that $x\sim c A_r$ and $B_r\sim d x^{2\frac {n-1}{n-2}}$ for some $c>1$ and $d>0$. Thus we get
 \[ \max_{1 \leq j \leq k} \big| H_r^{(j)}(x) \big|^{\frac k j} \leq \frac{C}{x^{k + \frac{n}{n-2}}},\]for $x$ big enough. Assertions i), ii) and iii) for the Hilbert-Schmidt norm then follow as above from \eqref{eq:hoelder_composition} and Proposition \ref{prop:Calpha_under_precomposition}. This completes the proof. \fin 
 
The previous lemma looks a lot like the conclusion of Theorem B, except that the precise exponents are not correct. The correct exponents are obtained by applying Proposition \ref{prop:SOn_coeff_biinvariant} for $\mathrm{SO}(m)$ for various $m \leq n$. 

\begin{lemma}\label{lemma:phi_C_alpha_using_SOm} 
Let $p > 2 + \frac{2}{n-2}$ and let $m \leq n$ be an integer such that $\frac{m-2}{2} - \frac{m-1}{p}$ is strictly positive. Consider $\beta \notin \N$ such that $\beta \leq \frac{m-2}{2} - \frac{m-1}{p}$. Then, the following estimates hold$\hskip1pt :$
  \begin{itemize}
  \item If $\beta >1$, then for every integer $1 \leq k \leq \beta$,
    \[ \big| \partial^k \hskip-1pt \varphi(x) \big| \leq \frac{C}{(x-1)^k x^{\frac{n}{m-2}}}.\]

  \item If $\beta < 1$, then
    \[ \big| \varphi(x) -\varphi(y) \big| \leq \sup_{x \le z \le y} \frac{C}{((z-1)z^{\frac{n}{m-2}})^\beta} |x-y|^\beta\]
    for every pair $x,y \in (1,\infty)$ satisfying that $x \leq y \leq x^{ 1 + \frac{n}{m-2}}$.
\end{itemize}
\end{lemma}

\dem The argument is the same as in Lemma \ref{lemma:phi_C_alpha_using_SOn}, except that $D$ is replaced by the diagonal matrix in $\mathrm{SL}_n(\R)$ with eigenvalues $e^r$ with multiplicity $1$, $e^s$ with multiplicity $m-1$ and $e^t$ with multiplicity $n-m$, where $r>0$ and $s,t$ are determined by $r$ as follows $$(s,t) = \Big( - \frac{n-m+2}{n+m-2}r, \frac{m-2}{n+m-2}r \Big).$$ As above, it turns out that $k \in \mathrm{SO}(m) \mapsto \varphi(\|D k D\|)$ is an $S_p$-$S_\infty$-multiplier by restriction. Moreover, the fact that $D$ has $m-1$ equal eigenvalues ensures that it is $\mathrm{SO}(m-1)$-biinvariant. In particular, we can use Proposition \ref{prop:SOn_coeff_biinvariant} for $\mathrm{SO}(m)$ combined with Proposition \ref{prop:Calpha_under_precomposition}. 

Any choice of $r,s,t$ with $r+(m-1)s+(n-m)t=0$ ensures that $\det D = 1$ and that $\| D k_\delta D \| = \max \{ e^{2s}, e^{2t},e^{r+s}g(\delta \sinh(r-s)) \}$. Our particular choice gives in addition that the operator norm of $D k_\delta D$ is equal to $e^{r+s}g(\delta \sinh(r-s))$. The only difference is that if $x,r,s$ are related by $x=e^{r+s}$ and $s = - \frac{n-m+2}{n+m-2}r$, then one gets $$e^{2r} = x^{1+\frac{n}{m-2}} \quad \mbox{and} \quad e^{2s} = x^{1-\frac{n}{m-2}}.$$ So if $H_r$ is still defined by \eqref{eq:def_Hr} and $x,r,s$ are related as above, then the estimates on the derivatives of $H_r$ become
\[ \big| \partial^j \hskip-1pt H_r(x) \big| \leq \frac{C}{(x-1)x^{j-1+\frac{n}{m-2}}} \ \Rightarrow \ \max_{1 \leq j \leq k} \big| \partial^j \hskip-1pt  H_r(x) \big|^{\frac k j} \leq \frac{C}{(x-1)^k x^{\frac{n}{m-2}}}.\]
Thus, the conclusion for the operator norm $\|\cdot\|$ is the same as for Lemma \ref{lemma:phi_C_alpha_using_SOn}.

On the other hand, the normalized Hilbert-Schmidt norm of $D k_\delta D$ has the form $\sqrt{A_r^2+\delta^2(B_r^2-A_r^2)}$ with
\begin{eqnarray*}
B_r^2 \!\! & = & \!\! \frac{1}{n} \left(e^{4r} + (m-1)e^{4s} +(n-m) e^{4t}\right), \\
A_r^2 \!\! & = & \!\! \frac{1}{n} \left(2 e^{2r+2s} + (m-2) e^{4s} +(n-m) e^{4t}\right).
\end{eqnarray*}
Therefore, when using radial Schur multipliers in the normalized Hilbert-Schmidt norm $|\cdot|$, if we define $\widetilde H_r$ by \eqref{eq:def_Hrtilde} with $x,r,s$ related as above, the same analysis for the derivatives of $\widetilde H_r$ can be applied and proves the lemma. 
\fin

\demB The conclusion of Theorem B follows from the preceding two lemmas. The fact that $\varphi$ is of class $\mathcal{C}^\alpha$ and the estimate on the local $\alpha - [ \alpha ]$ H\"older constant is contained in Lemma \ref{lemma:phi_C_alpha_using_SOn}. It remains to justify the pointwise estimates on $\varphi$ and its derivatives.

We start with the derivatives, for which the argument is direct. Fix an integer $1 \leq k <\alpha$. Consider $m$, the smallest integer such that $\beta := \frac{m-2}{2} - \frac{m-1}{p}$ is strictly greater than $k$. A small computation shows that it satisfies $$m-2= \left[ \frac{2k+1}{1 - \frac 2 p}\right ] \ge 3.$$ Note that $k<\beta \leq k+ \frac 1 2 - \frac 1 p$, so $\beta \notin \Z$ and Lemma \ref{lemma:phi_C_alpha_using_SOm} gives the expected estimate 
\begin{equation} \label{eq:derivativeThmB}
\big| \partial^k \hskip-1pt \varphi(x) \big| \leq \frac{C}{(x-1)^k x^{\frac{n}{m-2}}} = \frac{C}{(x-1)^k x^{c_k}} \quad \mbox{for every} \quad x \in (1,\infty).
\end{equation}

  To obtain pointwise estimates on $\varphi$, one more argument is needed. First observe that the case $x \leq 2$ is trivial because the norm of an $S_p$-multiplier is always bounded below by the $L_\infty$ norm of its symbol. So we can consider the case $x\geq 2$. When $\alpha<1$, Lemma \ref{lemma:phi_C_alpha_using_SOn} implies in particular that 
  \[ \big| \varphi(x) - \varphi(y) \big| \leq \frac{C}{x^{\alpha+\frac{\alpha n}{n-2}}} |x-y|^\alpha \leq \frac{C}{x^{\frac{\alpha n}{n-2}}} \quad \mbox{when} \quad 2 < x \le y \le 2x < \infty.\]
Thus $\varphi$ satisfies the Cauchy criterion, has a limit $\varphi_\infty$ and
  \[ |\varphi(x)-\varphi_\infty| \leq \sum_{i \geq 0} |\varphi(2^i x) - \varphi(2^{i+1}x)| \leq \frac{C}{x^{\frac{\alpha n}{n-2}}}.\]
This gives the assertion for $\alpha < 1$. When $\alpha > 1$, the inequality
  \[ \big| \varphi(x)-\varphi_\infty \big| \leq \int_x^\infty |\varphi'(y)| \, dy \leq \frac{C}{x^{c_1}}\]
is immediate from \eqref{eq:derivativeThmB}. This concludes the proof of Theorem B. \fin

\begin{remark}
\emph{Compared to the best known result \cite{dLdlS}, Theorem B gives:} 
\begin{itemize}
\item[i)] \emph{Automatic regularity of class $\mathcal{C}^\alpha$.} 

\item[ii)] \emph{More accurate decay rates for asymptotic rigidity of $\varphi$.}

\item[iii)] \emph{Local Mikhlin type conditions and asymptotic \lq\lq higher order\rq\rq${}$ rigidity.} 

\item[iv)] \emph{A larger range of $p$'s (in terms of $n$) for which the rigidity results hold.}
\end{itemize}
\end{remark}

\section{\bf Final comments} 

\noindent \textbf{A. Rank one.} According to \cite{CH}, we know that $\mathrm{SL}_2(\R)$ is weakly amenable. In particular, the rigidity theorems that shaped our statement in Theorem A do not apply to it. More precisely, there are some particularly well-behaved multipliers in $\mathrm{SL}_2(\R)$ |completely $L_p$-bounded by $1$ with Fourier symbols converging to $1$ uniformly on compact sets| which strongly break the decay implicit in Theorem A (as described in Remark \ref{Rem-Linear}). Is there a substantial improvement of Theorem A for $\mathrm{SL}_2(\R)$? The group $\mathrm{SL}_2(\R)$ does not admit finite-dimensional orthogonal cocycles and, consequently, the Mikhlin type conditions in \cite{JMP1,JMP2} do not apply. On the contrary, $\mathrm{SL}_2(\R)$ enjoys Haagerup property since it admits infinite-dimensional proper cocycles which lead to K-biinvariant associated length functions \cite[Chapter IV]{EF}. In this respect, we may construct noncommutative Riesz transforms for any such cocycle $\beta \hskip-2pt : \mathrm{SL}_2(\R) \to \mathcal{H}$ and any Riesz direction $u$ in the cocycle Hilbert space 
$$\widehat{R_uf}(g) \, = \, \frac{\langle \beta(g), u \rangle_\H}{\|\beta(g)\|_\H} \widehat{f}(g).$$ According to \cite{JMP2}, the maps $R_u$ are completely $L_p$-bounded and also satisfy more involved dimension free estimates. An optimal formulation of the Mikhlin condition in $\mathrm{SL}_2(\R)$ should include this natural class of multipliers. A quick inspection of the cocycle $\beta$ from \cite{EF} gives an asymptotic decay of order $(\log L)^{1/2}$, which is much less rigid than the behavior imposed by Theorem A. This indicates that there might be room for improvement in the rank 1 case. This construction is not possible in higher rank for the lack of such cocycles, due to Kahzdan property (T).  

\vskip3pt

\noindent \textbf{B. Twisted multipliers.} We have claimed in the Introduction that classical harmonic analysis methods are not efficient to give $L_p$-bounds of twisted Fourier multipliers. The first illustration of that was given in \cite{PRo}. There it was proved that twisted forms of $u$-directional Hilbert transforms are $L_p$-unbounded for all $p \neq 2$ unless the $\G$-orbit of $u$ is a finite set. The hidden tool here is Fefferman's Kakeya type construction for his ball multiplier theorem \cite{Fef}. On the contrary, the twisted one-dimensional Hilbert transforms are $L_p$-bounded for orthogonal actions \cite{JMP1}. This evidences that the twist is not stable under tensor product extensions! In addition, asymptotic Calder\'on-Zygmund methods worked in \cite{JMP1} for orthogonal actions, but become much less efficient for nonorthogonal ones due to the distortion effect of volume-preserving transformations. Namely,  combining Proposition \ref{CZProposition} i) with Junge's $H_p^c \to L_p$ inequality as in the proof of \cite[Theorem A]{JMP1}, we get a sufficient condition for $L_p$-boundedness $(1 < p < \infty)$ of twisted multipliers which vanish around $0$ $$\sup_{g \in \Sigma} |\xi|^{|\gamma|} \big| \partial_\xi^\gamma M_{\pm g}(\xi) \big| \, \lesssim \, 1 \ \Rightarrow \ \widetilde{T}_{\dot{m}} \hskip-2pt : L_p(\RR) \stackrel{\mathrm{cb}}{\longrightarrow} L_p(\RR)$$ for $M_{\pm g}(\xi) = |\xi|^{\pm \delta} \dot{m}(\alpha_g(\xi))$. This is effective for orthogonal actions. In the nonorthogonal case, multipliers with lower decay than in Theorem A satisfy the above condition for first-order derivatives, but fail it for higher orders (even for K-biinvariant multipliers in $\mathrm{SL}_2(\R)$) due to the distortion produced by $\alpha_g$.

\vskip3pt

\noindent \textbf{C. Calder\'on-Torchinsky theorem.} Is it true that 
\begin{equation} \label{Eq-CT} \tag{CT}
\weight{g}^{\hskip-3pt |\gamma|} \hskip-2pt \big| d_g^\gamma m(g) \big| \, \le \, C_{\mathrm{hm}} \quad \mbox{for all} \quad |\gamma| \le [s]+1
\end{equation}
suffices for the complete $L_p$-boundedness of the Fourier multiplier $T_m$ in the group algebra of $\SL$ whenever $|1/p - 1/2| < s/n^2$? This would be a natural analogue of Calder\'on-Torchinsky refinement of HM-condition \cite{CT,GHHN}, originally formulated using Sobolev spaces for fractional derivatives in the spirit of Remark \ref{Rem-SobolevLocal}. It is not hard to show that \eqref{Eq-CT} suffices locally. Indeed, the proof of the local form of Theorem A can not be directly modified since we make crucial use of Riesz transforms, for which we need full regularity. However, the proof gives an upper bound in terms of the Mikhlin condition for the lift $\dot{m}$. This bound still holds for the Sobolev condition in $\dot{m}$. In particular, the interpolation argument in \cite{GHHN} still applies. Unfortunately it seems much harder to interpolate the asymptotic behavior of the multiplier and we have no results in this direction.

The validity of \eqref{Eq-CT} beyond compactly supported symbols would be especially relevant in our context. Namely, Remark \ref{Rem-Linear} shows that the Fourier symbol decays as the largest Lie derivative in \eqref{Eq-CT}. In particular, \eqref{Eq-CT} imposes (as expected) less and less decay when $p$ approaches $2$. Moreover, working with fractional derivatives we might replace $[s]+1$ by $s+\varepsilon$, which can be arbitrarily close to $0$. Consequently, it is especially interesting to find a Sobolev formulation of Theorem A and the corresponding interpolated \lq CT-condition\rq${}$ to give room for $L_p$ multipliers with arbitrarily mild decay as $p \to 2$. Using left invariant Lie derivatives 
$$\lambda \big( \partial_\mathrm{X} m \big) = \int_\G \frac{d}{ds}_{\mid_{s=0}} \hskip-5pt m \big( g \exp(s \mathrm{X}) \big) \lambda(g) \, d\mu(g) = \lambda(m) \frac{d}{ds}_{\mid_{s=0}} \hskip-5pt \lambda(\exp(-s \mathrm{X})) = \lambda(m) a_\mathrm{X}$$
for $m$ regular enough. In particular, letting $A = - \sum_j a_{\mathrm{X}_j}^2$ for certain ONB $\mathrm{X}_1, \mathrm{X}_2, \ldots, \mathrm{X}_{\dim \G}$ of the Lie algebra $\mathfrak{g}$, we could define the Sobolev space $\mathsf{H}_{q,s}(\G)$ with \lq $s$-derivatives in $L_q$\rq${}$ as follows $$\|m\|_{\mathsf{H}_{q,s}(\G)} \, = \, \big\| \widehat{\lambda(m) (1+ A)^{\frac{s}{2}}} \big\|_{L_q(\G)}.$$ Compared to Remark \ref{Rem-SobolevLocal} this definition is intrinsic to $\G$. Nevertheless, a Sobolev form of the H\"ormander-Mikhlin condition requires to find a \lq dilation map\rq${}$ in the group, so that the resulting condition recovers the Mikhlin one for Lie derivatives as $q \to \infty$ for $s \in \Z_+$. We have no results in this direction. 

\vskip3pt

\noindent \textbf{D. Other Lie groups.} The local form of Theorem A generalizes to every real linear Lie group. Namely, every such group $\G$ admits a smooth embedding into $\SL$ for some large enough $n$. In particular, given a compactly supported symbol $m \hskip-2pt : \G \to \C$ satisfying H\"ormander-Mikhlin conditions up to order $[n^2/2]+1$, \hskip-1pt we may easily extend it to another symbol $M \hskip-2pt : \SL \to \C$ satisfying the same assumptions in $\SL$. By Theorem A, this implies that $T_M$ is an $L_p$-bounded Fourier multiplier in $\mathcal{L}(\SL)$, and hence also an $S_p$-bounded Schur multiplier in $\mathcal{B}(L_2(\SL))$. By the good restriction properties of Schur multipliers \cite{LdlS}, this implies $S_m \hskip-2pt : S_p(L_2(\G)) \to S_p(L_2(\G))$ is completely bounded for $1 < p < \infty$. Then local transference from Theorem \ref{Thm-LocalTransf} and interpolation give the following result. 

\begin{theorem}
Let $\G$ be a $\R$-linear Lie group and let $n$ be the minimal integer for which $\G$ embeds in $\SL$. Let $\mathrm{dist}$ be the distance associated with any left $\G$-invariant metric on $\G$ and assume that $m: \G \to \C$ is a compactly supported symbol in $\mathcal{C}^{k_n}(\G \setminus \{e\})$ for $k_n = [\frac{n^2}{2}]+1$ and satisfying 
$$\sup_{g \in \G \setminus \{e\}} \mathrm{dist}(g,e)^{|\gamma|} \big| d_g^\gamma m(g) \big| \, < \, \infty \quad \mbox{for all} \quad |\gamma| \le \Big[\frac{n^2}{2} \Big] + 1.$$
\noindent Then, the Fourier multiplier $T_m$ is completely $L_p$-bounded for all $1 < p < \infty$.
\end{theorem}

\vskip-2pt 

It would be very interesting to sharpen the above statement by lowering the differentiation order to $[\dim \G/2]+1$. As noticed in the Introduction, this also affects Theorem A since $\dim (\SL) = n^1-1 < n^2$. We leave it as an open problem for the interested reader. On the other hand, the new techniques in this paper are beyond the scope of \cite{GJP,JMP1,JMP2}, notably since we include nonorthogonal cocycles. This opens a door to investigate regularity conditions for $L_p$-multipliers in many other unimodular Lie groups. 

\noindent \textbf{Acknowledgement.} We thank Adri\'an M. Gonz\'alez-P\'erez for numerous comments and fruitful discussions. J. Parcet was partially supported by Europa Excelencia Grant MTM2016-81700-ERC, CSIC Grant PIE-201650E030, Spain Grant PID2019-107914GB-I00 and ICMAT Severo Ochoa Grant CEX2019-000904-S (Spain). M. de la Salle was partially supported by ANR grants GAMME and AGIRA. M. de la Salle and \'E. Ricard were partially supported by the ANR grant ANCG ANR-19-CE40-0002. Part of this work was carried out during a long term visit of J. Parcet to the Laboratoire de Math\'ematiques Nicolas Oresme at Universit\'e de Caen Normandie. The first-named author would like to express his gratitude to the members of the Laboratoire for their hospitality. 

\bibliographystyle{amsplain}

\vskip10pt

\hfill \noindent \textbf{Javier Parcet} \\
\null \hfill Instituto de Ciencias Matem{\'a}ticas 
\\ \null \hfill Consejo Superior de
Investigaciones Cient{\'\i}ficas \\ \null \hfill C/ Nicol\'as Cabrera 13-15.
28049, Madrid. Spain \\ \null \hfill\texttt{javier.parcet@icmat.es}

\vskip2pt

\hfill \noindent \textbf{\'Eric Ricard} \\
\null \hfill Normandie Univ\\ \null \hfill
\null \hfill UNICAEN, CNRS \\
\null \hfill LMNO
\\ \null \hfill 14000 Caen, France \\ \null \hfill\texttt{eric.ricard@unicaen.fr}

\vskip2pt

\hfill \noindent \textbf{Mikael de la Salle} \\
\null \hfill UMPA \\ \null \hfill
\null \hfill CNRS -- ENS de Lyon 
\\ \null \hfill 69364 Lyon Cedex 7, France \\ \null \hfill\texttt{mikael.de.la.salle@ens-lyon.fr}
\end{document}